\pdfoutput=1
\RequirePackage{silence}
\documentclass[a4paper,11pt]{amsart}
\usepackage[hmarginratio={1:1},vmarginratio={1:1},lmargin=60.0pt,tmargin=60.0pt]{geometry}

\usepackage[numbers]{natbib}

\allowdisplaybreaks

\usepackage[utf8]{inputenc}

\usepackage{amssymb,amsmath,amsthm,amsfonts,mathrsfs,bbm}
\usepackage[table]{xcolor}
\usepackage{graphicx}
\usepackage{mathtools}

\definecolor{spinach}{RGB}{46,139,87}
\definecolor{tomato}{RGB}{255,99,71}
\definecolor{pumpkin}{RGB}{224,180,80}

\definecolor{orchid}{RGB}{143,40,194}
\definecolor{lava}{RGB}{207,16,32}
\definecolor{mydarkblue}{RGB}{10,10,150}

\usepackage[yyyymmdd,hhmmss]{datetime}
\usepackage{todonotes}

\usepackage{xparse}

\usepackage{enumitem}
\setlist[enumerate]{itemsep=0.15cm,label=\emph{\upshape(\alph*)}}
\setlist[enumerate,2]{itemsep=0.15cm,label=\emph{\upshape(\roman*)}}

\let\emph\relax
\DeclareTextFontCommand{\emph}{\em}

\usepackage{array}
\newcolumntype{C}{>{$}c<{$}}

\usepackage[vcentermath]{youngtab}

\usepackage[all]{xy}
\usepackage{tikz}
\usetikzlibrary{cd}
\usetikzlibrary{decorations}
\usetikzlibrary{decorations.markings}
\usetikzlibrary{decorations.pathreplacing}
\usetikzlibrary{decorations.pathmorphing}
\usetikzlibrary{arrows.meta,shapes,positioning,matrix,calc}
\usetikzlibrary{shapes.callouts}
\tikzset{anchorbase/.style={baseline={([yshift=-0.5ex]current bounding box.center)}},
tinynodes/.style={font=\tiny, text height=0.25ex, text depth=0.05ex},
smallnodes/.style={font=\scriptsize, text height=0.75ex, text depth=0.15ex},
crossline/.style={preaction={draw=white,line width=5.0pt,-},preaction={draw=black,line width=0.9pt,-}},
usual/.style={line width=1.0,color=black},
dot/.style = {
decoration={markings,
post length=0.25mm,
pre length=0.25mm,
mark=at position #1 with {\node[circle,radius=0.2cm,inner sep=-1.5pt,color=black,fill=black]{};}
},
postaction={decorate}
},
dot/.default=1,
}

\newcommand{\ie}{\text{i.e.}}
\newcommand{\eg}{\text{e.g.}}
\newcommand{\cf}{\text{cf.}}
\newcommand{\etc}{\text{etc.}}
\newcommand{\aka}{\text{a.k.a.}}
\newcommand{\ver}{\text{verbatim}}

\newcommand{\muta}{\text{mutatis mutandis}}

\newcommand{\acts}{\centerdot}
\renewcommand{\dots}{\text{...}}

\newcommand{\placeholder}{{}_{-}}

\newcommand{\lra}{\longrightarrow}
\newcommand{\floor}[1]{\lfloor #1 \rfloor}

\newcommand{\C}{\mathbb{C}}
\newcommand{\Z}{\mathbb{Z}}
\newcommand{\R}{\mathbb{R}}
\newcommand{\Q}{\mathbb{Q}}
\newcommand{\F}{\mathbb{F}}
\newcommand{\K}{\mathbb{K}}
\newcommand{\LL}{\mathbb{L}}

\newcommand{\N}{\mathbb{Z}_{\geq 0}}
\newcommand{\cchar}[1][\K]{\mathrm{char}(#1)}

\newcommand{\Aut}{\mathrm{Aut}}
\newcommand{\End}{\mathrm{End}}
\newcommand{\Hom}{\mathrm{Hom}}
\newcommand{\Ext}{\mathrm{Ext}}
\newcommand{\HH}{\mathrm{H}}

\newcommand{\id}{id}

\newcommand{\monoid}[1][S]{\mathcal{#1}}
\newcommand{\onemon}{\mathcal{S}_{1}}
\newcommand{\xmon}[1][X]{\mathcal{#1}}
\newcommand{\category}[1][S]{\mathbf{#1}}
\newcommand{\group}[1][G]{\mathcal{#1}}
\newcommand{\msg}{\monoid\setminus\group}

\newcommand{\dimk}[1][\K]{\mathrm{dim}_{#1}}

\newcommand{\module}[1][M]{{#1}}
\newcommand{\one}{\mathbbm{1}}
\newcommand{\oneb}{\mathbbm{1}_{b}}
\newcommand{\onet}{\mathbbm{1}_{t}}
\newcommand{\onebt}{\mathbbm{1}_{bt}}
\newcommand{\simple}[1][K]{L_{#1}}
\newcommand{\fmodule}[1][F]{{#1}}

\newcommand{\gap}[2][\K]{\mathrm{gap}_{#1}(#2)}
\newcommand{\faith}[2][\K]{\mathrm{faith}_{#1}(#2)}
\newcommand{\faithb}[2][\K]{\mathrm{faith}_{#1}\big(#2\big)}
\newcommand{\gratio}[2][\K]{\mathrm{gapr}_{#1}(#2)}

\newcommand{\fratio}[2][\K]{\mathrm{faithr}_{#1}(#2)}

\newcommand{\brgr}[1][n]{\mathcal{B}\mathrm{r}_{#1}}
\newcommand{\tmon}[1][n]{\mathcal{T}_{#1}}
\newcommand{\sym}[1][n]{\mathcal{S}_{#1}}
\newcommand{\cyclic}[1][n]{\mathcal{C}_{#1}}
\newcommand{\alt}[1][n]{\mathcal{A}_{#1}}
\newcommand{\GL}[2][n]{\mathcal{GL}_{#1}(#2)}
\newcommand{\PSL}[2][n]{\mathcal{PSL}_{#1}(#2)}

\newcommand{\lcell}{\mathcal{L}}
\newcommand{\rcell}{\mathcal{R}}
\newcommand{\jcell}{\mathcal{J}}
\newcommand{\hcell}{\mathcal{H}}

\newcommand{\jb}{\mathcal{J}_{b}}
\newcommand{\jm}{\mathcal{J}_{m}}
\newcommand{\jt}{\mathcal{J}_{t}}

\newcommand{\kcell}{\mathcal{K}}
\newcommand{\mcell}{\mathcal{M}}

\newcommand{\lmod}[1][\lcell]{\Delta_{#1}}
\newcommand{\rmod}[1][\rcell]{{}_{#1}\Delta}

\newcommand{\hd}{\mathrm{Hd}}
\newcommand{\ind}{\mathrm{Ind}}
\newcommand{\rad}[1][\jcell]{\mathrm{Rad}}
\newcommand{\radd}[1][\jcell]{\mathrm{Rad}_{#1}}

\newcommand{\ssdimk}[1][\K]{\mathrm{ssdim}_{\K}}
\newcommand{\ssgap}[2][\K]{\mathrm{ssgap}_{#1}(#2)}
\newcommand{\ssratio}[2][\K]{\mathrm{ssgapr}_{#1}(#2)}

\newcommand{\tlcat}[1][\delta]{\mathbf{TL}^{lin}(#1)}
\newcommand{\tlalg}[2][n]{\mathcal{TL}_{#1}^{lin}(#2)}
\newcommand{\tlset}{\mathbf{TL}}
\newcommand{\tlmon}[1][n]{\mathcal{TL}_{#1}}

\newcommand{\tlord}{\vartriangleleft}
\newcommand{\tlordsecond}{\tlord^{\prime}}

\newcommand{\tltru}[2][n]{\mathcal{TL}_{#1}^{\leq{#2}}}
\newcommand{\protru}[3][n]{\mathrm{p}\mathcal{R}\mathrm{o}_{#1}^{\leq{#2},<{#3}}}
\newcommand{\motru}[2][n]{\mathcal{M}\mathrm{o}_{#1}^{\leq{#2}}}
\newcommand{\ppatru}[2][n]{\mathrm{p}\mathcal{P}\mathrm{a}_{#1}^{\leq{#2}}}

\newcommand{\nobottom}[2]{\mathrm{B}_{#1}^{#2}}
\newcommand{\vertical}[2]{\mathrm{Vert}_{#1}^{#2}}
\newcommand{\wvertical}[2]{\mathrm{WVert}_{#1}^{#2}}
\newcommand{\tlgraph}[2]{\Gamma_{#1}^{#2}}
\newcommand{\tlgraphtwo}[2]{\Delta_{#1}^{#2}}

\newcommand{\brcat}[1][\delta]{\mathbf{Br}^{lin}(#1)}
\newcommand{\bralg}[2][n]{\mathcal{B}\mathrm{r}_{#1}^{lin}(#2)}
\newcommand{\brset}{\mathbf{Br}}
\newcommand{\brmon}[1][n]{\mathcal{B}\mathrm{r}_{#1}}

\newcommand{\brtru}[2][n]{\mathcal{B}\mathrm{r}_{#1}^{\leq{#2}}}
\newcommand{\rotru}[3][n]{\mathcal{R}\mathrm{o}_{#1}^{\leq{#2},<{#3}}}
\newcommand{\robrtru}[2][n]{\mathcal{R}\mathrm{o}\mathcal{B}\mathrm{r}_{#1}^{\leq{#2}}}
\newcommand{\patru}[2][n]{\mathcal{P}\mathrm{a}_{#1}^{\leq{#2}}}

\newcommand{\pamon}[1][n]{\mathcal{P}\mathrm{a}_{#1}}
\newcommand{\robrmon}[1][n]{\mathcal{R}\mathrm{o}\mathcal{B}\mathrm{r}_{#1}}
\newcommand{\romon}[1][n]{\mathcal{R}\mathrm{o}_{#1}}
\newcommand{\ppamon}[1][n]{\mathrm{p}\mathcal{P}\mathrm{a}_{#1}}
\newcommand{\probrmon}[1][n]{\mathrm{p}\mathcal{R}\mathrm{o}\mathcal{B}\mathrm{r}_{#1}}
\newcommand{\momon}[1][n]{\mathcal{M}\mathrm{o}_{#1}}
\newcommand{\promon}[1][n]{\mathrm{p}\mathcal{R}\mathrm{o}_{#1}}
\newcommand{\pbrmon}[1][n]{\mathrm{p}\mathcal{B}\mathrm{r}_{#1}}
\newcommand{\psym}[1][n]{\mathrm{p}\mathcal{S}_{#1}}

\DeclarePairedDelimiterX{\inner}[2]{\langle}{\rangle}{#1, #2}

\DeclareMathOperator{\rank}{rank}

\usepackage{aliascnt}
\def\NewTheorem#1{%
\newaliascnt{#1}{equation}%
\newtheorem{#1}[#1]{#1}%
\aliascntresetthe{#1}%
\expandafter\def\csname #1autorefname\endcsname{#1}%
}
\def\equationautorefname~#1\null{(#1)\null}

\numberwithin{equation}{subsection}

\NewTheorem{Proposition}
\NewTheorem{Theorem}
\NewTheorem{Corollary}
\AtEndEnvironment{Corollary}{\null\hfill$\square$}%
\NewTheorem{Lemma}
\theoremstyle{definition}
\NewTheorem{Definition}
\NewTheorem{Notation}
\NewTheorem{Example}
\AtEndEnvironment{Example}{\null\hfill$\Diamond$}%
\theoremstyle{remark}
\NewTheorem{Remark}
\NewTheorem{Assumption}
\NewTheorem{Hypothesis}
\NewTheorem{Conjecture}
\NewTheorem{Question}
\NewTheorem{Observation}
\NewTheorem{Fact}
\NewTheorem{Conclusion}
\NewTheorem{lemma}

\usepackage[hypertexnames=false]{hyperref}
\usepackage{bookmark}
\hypersetup{
pdftoolbar=true,
pdfmenubar=true,
pdffitwindow=false,
pdfstartview={FitH},
pdftitle={Monoidal categories, representation gap and cryptography},
pdfauthor={Mikhail Khovanov, Maithreya Sitaraman and Daniel Tubbenhauer},
pdfsubject={},
pdfcreator={Mikhail Khovanov, Maithreya Sitaraman and Daniel Tubbenhauer},
pdfproducer={Mikhail Khovanov, Maithreya Sitaraman and Daniel Tubbenhauer},
pdfkeywords={},
pdfnewwindow=true,
colorlinks=true,
linkcolor=mydarkblue,
citecolor=teal,
filecolor=magenta,
urlcolor=orchid,
linkbordercolor=lava,
citebordercolor=teal,
urlbordercolor=orchid,
linktocpage=true
}

\newcommand{\nnfootnote}[1]{%
\begin{NoHyper}
\renewcommand\thefootnote{}\footnote{#1}%
\addtocounter{footnote}{-1}%
\end{NoHyper}
}

\def\makeautorefname#1#2{\csdef{#1autorefname}{#2}}

\makeautorefname{section}{Section}%
\makeautorefname{subsection}{Section}%
\makeautorefname{subsubsection}{Section}%

\begin{document}
\title[Monoidal categories...and cryptography]
{Monoidal categories, representation gap and cryptography}
\author[M. Khovanov, M. Sitaraman and D. Tubbenhauer]{Mikhail Khovanov, Maithreya Sitaraman and Daniel Tubbenhauer}

\address{M.K.: Department of Mathematics, Columbia University, New York, NY 10027, United States,\\ \href{https://www.math.columbia.edu/~khovanov/}{www.math.columbia.edu/~khovanov/}, \href{https://orcid.org/0009-0006-5296-6044}{ORCID 0009-0006-5296-6044}} 
\email{khovanov@math.columbia.edu}

\address{M.S.: Department of Mathematics, Columbia University, New York, NY 10027, United States}
\email{maithreya@math.columbia.edu}

\address{D.T.: The University of Sydney, School of Mathematics and Statistics F07, Office Carslaw 827, NSW 2006, Australia, \href{http://www.dtubbenhauer.com}{www.dtubbenhauer.com}, \href{https://orcid.org/0000-0001-7265-5047}{ORCID 0000-0001-7265-5047}}
\email{daniel.tubbenhauer@sydney.edu.au}

\begin{abstract}
The linear decomposition attack provides a serious obstacle to direct applications of noncommutative groups and monoids (or semigroups) in cryptography. To overcome this issue we propose to look
at monoids with only big representations, in the sense made precise in the paper, and undertake a systematic study of such monoids. 
One of our main tools is Green's theory of cells (Green's relations).

A large supply of monoids is delivered by monoidal categories. We consider simple examples of monoidal categories of diagrammatic origin, including the Temperley--Lieb, the Brauer and partition categories, and discuss lower bounds for their representations.
\end{abstract}

\nnfootnote{\textit{Mathematics Subject Classification 2020.} Primary: 18M05, 20M30; Secondary: 94A60.}
\nnfootnote{\textit{Keywords.} Monoidal categories, diagram categories, monoid/semigroup representations, representation gap, cryptography.}

\maketitle

\tableofcontents

\arrayrulewidth=0.5mm
\setlength{\arrayrulewidth}{0.5mm}


\section{Introduction}\label{S:Intro}

The main goal of this paper is to start  
connecting monoidal categories and cryptography.


\subsection{Protocols and platform groups}\label{SS:IntroProtocols}

Some of 
the most important cryptographic protocols in use today are based on \emph{commutative groups} and deliver a gold standard for cryptography (modulo the fear of quantum computers). On the other hand, \emph{noncommutative group-based} and \emph{monoid-based} (or \emph{semigroup-based}, but we will stay with monoids in this paper) protocols seem to be less understood and in many cases admit efficient attacks.

Exceptionally successful Diffie--Hellman (DH), Rivest--Shamir--Adleman (RSA) and elliptic curve cryptography algorithms, see {\eg} \cite{Ko-algebraic-cryptography}, \cite{Wa-elliptic-curves},
are based on the commutative group $(\Z/n\Z)^{\ast}$ of invertible residues modulo $n$ and on the group of points on an elliptic curve $E$ over a finite field $\F_{q}$, respectively. Here one usually wants these groups to contain a subgroup of large prime order and small index. 
For example, in the classical DH protocol the prime $p$ as well as 
a generator $g\in(\Z/p\Z)^{\ast}$ of the multiplicative group are public. 
Then party A 
chooses privately $a\in\Z$ and party B chooses privately $b\in\Z$. 
Party A communicates $g^{a}$, B sends $g^{b}$ and the 
common secret is $(g^{b})^{a}=g^{ab}=(g^{a})^{b}$.
A third party C has access to $n$, $g$, $g^{a}$ and $g^{b}$, but finding 
$g^{ab}$ from the known data is difficult as long as $p-1$ contains a large prime among its factors. 

There has been many ideas and there is an extensive literature on constructing cryptographic protocols from noncommutative groups and monoids (monoids generalize groups and we switch to saying monoids from now on), 
see {\eg} \cite{MyShUs-group-cryptography}, \cite{MyShUs-noncom-cryptography} and references therein. Examples of 
such are Magyarik--Wagner public key protocol \cite{WaMa-cryptosystem-word-problem},
Anshel--Anshel--Goldfeld key exchange \cite{AnAnGo-algebraic-cryptography}, 
Ko--Lee et al. key exchange protocol \cite{KoLeChHaKaPa-crypto-braids} and Shpilrain--Zapata public key protocols \cite{ShZa-groups-public-key-cryptography}.

In the literature the monoid $\monoid$ used in protocols is 
often called the platform group/monoid. In \cite[Section 4]{MyRo-linear-attack} 
there is a big list of various protocols and platform monoids, 
including but not limited to the ones named above.
Sometimes these
restrict to groups or matrix groups, sometimes general monoids can be used.
A prototypical example for this paper is the Shpilrain--Ushakov (SU) key exchange protocol, see {\eg} \cite[Section 4.2.1]{MyShUs-group-cryptography}, which works as follows. The public data is a monoid $\monoid$, and 
two sets $A,B\subset\monoid$ of commuting elements and $g\in\monoid$. Party A 
chooses privately $a,a^{\prime}\in A$ and party B chooses privately $b,b^{\prime}\in A$. Party $A$ communicates $aga^{\prime}$, B sends $bgb^{\prime}$ and the common secret is $abgb^{\prime}a^{\prime}=baga^{\prime}b^{\prime}$. Another example that does not use commuting elements 
is Stickel's secret key exchange (St) \cite{St-new-key-exchange}. Here $g,h\in\monoid$ with $gh\neq hg$ are public, party A picks $a,a^{\prime}\in\N$, party 
B picks $a,a^{\prime}\in\N$, A sends $g^{a}h^{a^{\prime}}$, B sends $g^{b}h^{b^{\prime}}$, and the common secret is 
$g^{a}g^{b}h^{b^{\prime}}h^{a^{\prime}}=g^{b}g^{a}h^{a^{\prime}}h^{b^{\prime}}$.
Note that $\monoid$ can be an arbitrary monoid in these protocols. The 
complexity of $\monoid$ determines how difficult it is to find the common secret from the public data.

As shown by Myasnikov and Roman'kov \cite{MyRo-linear-attack} and also based on earlier work, the SU and St protocols and others in this spirit, the ones named two paragraphs above included,
can be successfully attacked if $\monoid$ admits small nontrivial representations. This is called a \emph{linear decomposition attack} or \emph{linear attack}, for short.

One of the consequences of linear attacks is 
that finite noncommutative groups may not be suited for 
cryptographic purposes as they admit nontrivial 
representations of moderate size. For a toy example, the symmetric group $\sym$ has $n!$ elements, but admits a faithful $(n-1)$-dimensional representation. The dimension of this representation is smaller than logarithmic in the size of the group, and the symmetric group would be a poor choice for various standard noncommutative group protocols. 
Likewise, finite simple groups of Lie type often admit representations of (exponentially) small dimension compared to their size. 
With few exceptions, including cyclic groups of prime order, which are related to the classical and well-understood protocols, the same is true for other finite simple groups. That is, these groups admit nontrivial representations of small dimension relative to their order. Since any finite group $G$ surjects onto some finite simple group, reducing the problem of bounding representations of $G$ from below to that of the simple quotient, linear attacks rule out many finite noncommutative groups. 

Hence, it is not surprising that some platform groups proposed in the literature are infinite, 
{\eg} Artin--Tits, Thompson or Grigorchuk groups, see 
\cite[Chapter 5]{MyShUs-group-cryptography}.
\medskip 

This paper explores finite monoids (mostly coming from monoidal categories) 
instead of infinite groups.
The questions we address are:
\begin{enumerate}[label=$\bullet$]

\item What are (numerical) measures to determine whether 
a monoid can resist linear attacks?

\item How to find a good supply of finite monoids for cryptographic 
use?

\end{enumerate}


\subsection{Linear attacks, representation gap and faithfulness}\label{SS:IntroObservations}

The following observations 
regarding monoid-based 
cryptography are our starting points:  

\begin{enumerate}
\item As explained above, monoid-based protocols such as SU or St and many others
often admit efficient attacks based on linear algebra \cite{MyRo-linear-attack}.

\item A natural solution to this problem is to restrict to monoids that have nontrivial representations only starting from a suitably big dimension. We call the smallest dimension of a nontrivial $\monoid$-representation the \emph{representation gap} of $\monoid$. Alternatively and weaker, 
we also ask for the dimension of the smallest faithful
$\monoid$-representation to be big, and we call this measure the 
\emph{faithfulness} of $\monoid$. We elaborate on these in \autoref{S:RepGap}.
\end{enumerate} 

\begin{Remark}\label{R:IntroReferences}
Various monoid invariants similar 
to the representation gap and its companions 
have appeared in the literature and we give some references 
in the main body of the paper. However, the motivations to study these invariants in the literature are very different 
from ours, and it would be very interesting 
to make a connection to cryptography starting from these works.
\end{Remark}

It is thus essential to find monoids  
that have big representation gaps or with faithful representations of big dimension only. Suitably defined, a big representation gap 
or big faithfulness seem to be necessary, but not sufficient, 
conditions for a monoid to be potentially useful in cryptography, however.
Moreover, one problem not discussed here is potential information loss: multiplication by an element of a monoid may not be invertible.


\subsection{Monoidal categories and monoids}\label{SS:IntroMonoidal}

A category delivers a supply of monoids: any object $X$ of a category $\category$ gives 
rise to the monoid $\monoid=\End_{\category}(X)$ of its endomorphisms. It is further natural to 
consider \emph{monoidal categories}, where objects can be tensored 
$\otimes$ subject to suitable axioms, for the following reasons:

\begin{enumerate}[resume]
\item It would be preferable to have a family 
of monoids $\monoid_{n}$, say one for each $n\in\N$.
This is where monoidal categories enter.
A single object $X$ of a monoidal category $\category$ produces a family of monoids $\{\monoid_{n}=\End_{\category}(X^{\otimes n})|n\in\N\}$.

\item Commuting actions play a key role in cryptography, {\cf} the SU protocol recalled above. Such commuting actions 
exists naturally in the setting of categories and monoidal categories.
For any pair of objects $X$, $Y$ of a category $\category$, not necessarily monoidal, there is a commuting action of the monoids $\End_{\category}(X)^{op}$ (the opposite monoid) and $\End_{\category}(Y)$ on the set $\Hom_{\category}(X,Y)$. Thus, categories immediately produce a significant amount of commuting actions. 
Furthermore, monoidal categories provide an even richer supply of such actions: for any two objects $X$, $Y$ the actions of 
the monoids $\End_{\category}(X)$ and $\End_{\category}(Y)$ on $X\otimes Y$ commute. It is easy to convert these to commuting actions on sets, for instance, on the set $\Hom_{\category}(Z,X\otimes Y)$ for $Z\in\mathrm{Ob}(\category)$.

\item Monoidal categories are naturally two-dimensional structures. They 
often can be described via generating objects, generating morphisms and defining relations. The latter can be understood as relations on planar diagrams or networks, see {\eg} \cite{Se-survey-monoidal-diagrams}, \cite{TuVi-monoidal-tqft}. A natural problem is to construct examples of diagrammatically defined monoidal categories that may be useful for cryptographic purposes. We start tackling this for planar diagrammatics in \autoref{S:TL} 
and for diagrammatics involving permutation symmetries in \autoref{S:Brauer}.

\end{enumerate}

These are our reasons to study (diagram) monoids coming 
from monoidal categories and we elaborate on their 
potential usefulness in cryptography in the main body of the text.

There are then three additional facts regarding this project that we stress 
and that we think makes our discussion interesting:

\begin{enumerate}[resume]
\item The current literature on monoidal categories (see 
for example \cite{EtGeNiOs-tensor-categories}, \cite{TuVi-monoidal-tqft} and references therein) mostly studies $\K$-linear categories or variations of such. This means hom-spaces between the objects are $\K$-vector spaces for some field $\K$. Such categories are not immediately useful from the cryptographic or any classical computation viewpoint, since it usually takes a prohibitive amount of data to record an element of the hom-space between two objects (those hom-spaces tend to have exponentially big dimensions). One the other hand, protocols in $\K$-linear categories with homs between objects having moderate dimensions can be dealt with via linear decomposition attacks, see \cite{MyRo-linear-attack}. 

It makes sense to develop \emph{set-theoretic counterparts} of categories that appear in quantum algebra, quantum topology, mathematical physics, and TQFTs, and see whether related monoids have big representation gaps. We provide easy examples of such in the present paper and discuss their usefulness 
for cryptography, see parts of \autoref{S:TL} 
and \autoref{S:Brauer}.

\item It seems hard to build secure cryptographic protocols from noncommutative finite groups, due to finite simple groups having small representation gaps relative to their size. For example, among finite simple groups only the cyclic groups (of prime order) appear to be well-behaved 
for cryptographical purposes, {\cf} \autoref{E:RepGapSimpleGroups} 
and \autoref{E:RepGapSimpleGroups2}.

One of our points is that representation gaps and faithfulness \emph{tend to be bigger for suitable monoids} than for groups when controlling for size. The abstract theory of monoid 
representations should be useful for some general statements in this direction, see \autoref{S:Cells} for some first steps.

\item Finally, \emph{lower bounds on dimensions} of representations of monoids or \emph{growth rates} of such dimensions
are not yet extensively studied in the literature, even not for group (representation theorists seem to prefer precise numbers). 

Part of this project is also
to get good bounds and growth rates for simple and faithful $\monoid$-representations, 
for finite monoids $\monoid$, see 
\autoref{T:TLIsAGoodExample} for an example.
\end{enumerate}


\subsection{Cell theory and cryptography}\label{SS:IntroCells}

Our main tool to study monoid representations 
are \emph{Green's relations} {\aka} \emph{Green's theory of cells}. We explain 
the details in \autoref{S:RepGap}.

An example of how cell theory enters the paper is that
a monoid $\monoid$ can be truncated by considering a large cell submonoid $\monoid^{\geq\jcell}$, see \autoref{S:RepGap} for definition. 
Since simple $\monoid$-representations are ordered by cells, 
$\monoid^{\geq\jcell}$ will inherit precisely the simple $\monoid$-representations for large cells. The monoids of the form $\monoid^{\geq\jcell}$ sometimes 
have very few small representations. This truncation works actually 
in two ways, from above and from below, 
using Rees factors and cell truncations, and provides    
a good way to get rid of unwanted representations, {\cf} \autoref{SS:CellsSubmon}.

Moreover, in \autoref{S:Cells} we will discuss so-called $H$-cells,  
how they control the representation theory of the monoids and how large cells  
resist against linear attacks.

Another way cells help to determine whether a given monoid could resist 
linear attacks is that they give rise to 
what we call the \emph{semisimple representation gap}, which measures the normalized size of the cells.  
This numerical value is not as fine as the representation gap or the faithfulness, but easier to compute and agrees with the representation gap in the semisimple situation.

The representation gap, 
the semisimple representation gap and the faithfulness seem to be good 
first tests for determining whether a given monoid resists linear attacks.
Throughout the text we list a few additional properties, partially motivated 
by cell theory, that may be useful for cryptographical applications.


\subsection{Representation gap in some diagrammatic monoids}\label{SS:IntroDiaMonoids}

Let us take the opportunity to recall some diagrammatic 
monoids which we will discuss in this paper. All of these will be very familiar 
to the reader with background in quantum algebra, quantum topology and alike, but they also are prominent examples in monoid theory.

We will be very brief and details and references will follow in the main text. We also indicate 
whether these monoids might be useful for cryptography in the sense 
of having substantial (semisimple) representation gaps or only big faithful representations.

Most of the monoids which we will use can be obtained 
as hom-subsets of the \emph{set-theoretical partition category}. We will use 
matchings from $n$ bottom to $n$ top points of the following types (all of these are classical example, see {\eg} \cite{HaRa-partition-algebras} or \cite{HaJa-representations-diagram-algebras} for summaries):
\begin{enumerate}[label=$\bullet$]

\item The \emph{partition monoid} $\pamon$ of all diagrams of partitions of a $2n$-element set.

\item The \emph{rook-Brauer monoid} $\robrmon$ consisting of all diagrams with components of size $1,2$.

\item The \emph{Brauer monoid} $\brmon$ consisting of all diagrams with components of size $2$.

\item The \emph{rook monoid} $\romon$ consisting of all diagrams with components of size $1,2$, and all partitions have at most one component 
at the bottom and at most one at the top.

\item The \emph{symmetric group} $\sym$ consisting of all matchings with components of size $1$.

\item \emph{Planar} versions of these: $\ppamon$, $\probrmon=\momon$, $\pbrmon=\tlmon$, $\promon$ and $\psym\cong\onemon$ (the latter denotes the 
\emph{trivial monoid}). The planar rook-Brauer monoid is also called \emph{Motzkin monoid}, the planar Brauer monoid is also known as the \emph{Temperley--Lieb monoid}, and the planar symmetric group is trivial.

\end{enumerate}

\begin{Remark}\label{R:IntroNames}
The above diagram monoids appear in many different 
fields of mathematics. This makes them on the one hand very appealing, but on the other hand tends to cause confusion from time to time. For example, as we already indicated above, these diagram monoids have different names that vary with the field, {\eg} the Temperley--Lieb monoid is also known as the Jones monoid or the Kauffman monoid in monoid theory, but that name appears 
to be unheard-of in the representation theoretical literature on the algebra versions of these monoids.
\end{Remark}

\autoref{Eq:IntroDiaMonoids} summarizes our list, see also 
\cite[Section 2.3]{HaJa-representations-diagram-algebras}. 
In order to make components of size one 
visible we use loose dotted ends. 
We also indicate whether their nontrivial representations are reasonably big
(the ``Big reps''
column), meaning after 
appropriate cell truncation. Hereby ${}^{\ast}$ means that they have such representations but still come with an aftertaste (such as being semisimple in some cases), ${}_{c}$ means conjectural, and EX means excluded from the 
discussion due to triviality.
This is explained in more details in \autoref{C:TLConclusion} and 
\autoref{C:BrauerConclusion}.
\begin{gather}\label{Eq:IntroDiaMonoids}
\begin{tabular}{c|c|c||c|c|c}
\arrayrulecolor{tomato}
Symbol & Diagrams & Big reps 
& Symbol & Diagrams & Big reps
\\
\hline
\hline
$\ppamon$ & \begin{tikzpicture}[anchorbase]
\draw[usual] (0.5,0) to[out=90,in=180] (1.25,0.45) to[out=0,in=90] (2,0);
\draw[usual] (0.5,0) to[out=90,in=180] (1,0.35) to[out=0,in=90] (1.5,0);
\draw[usual] (0.5,1) to[out=270,in=180] (1,0.55) to[out=0,in=270] (1.5,1);
\draw[usual] (1.5,1) to[out=270,in=180] (2,0.55) to[out=0,in=270] (2.5,1);
\draw[usual] (0,0) to (0,1);
\draw[usual] (2.5,0) to (2.5,1);
\draw[usual,dot] (1,0) to (1,0.2);
\draw[usual,dot] (1,1) to (1,0.8);
\draw[usual,dot] (2,1) to (2,0.8);
\end{tikzpicture} & YES${}^{\ast}$
& $\pamon$ & \begin{tikzpicture}[anchorbase]
\draw[usual] (0.5,0) to[out=90,in=180] (1.25,0.45) to[out=0,in=90] (2,0);
\draw[usual] (0.5,0) to[out=90,in=180] (1,0.35) to[out=0,in=90] (1.5,0);
\draw[usual] (0,1) to[out=270,in=180] (0.75,0.55) to[out=0,in=270] (1.5,1);
\draw[usual] (1.5,1) to[out=270,in=180] (2,0.55) to[out=0,in=270] (2.5,1);
\draw[usual] (0,0) to (0.5,1);
\draw[usual] (1,0) to (1,1);
\draw[usual] (2.5,0) to (2.5,1);
\draw[usual,dot] (2,1) to (2,0.8);
\end{tikzpicture} & YES${}^{\ast}_{c}$
\\
\hline
$\momon$ & \begin{tikzpicture}[anchorbase]
\draw[usual] (0.5,0) to[out=90,in=180] (1.25,0.5) to[out=0,in=90] (2,0);
\draw[usual] (1,0) to[out=90,in=180] (1.25,0.25) to[out=0,in=90] (1.5,0);
\draw[usual] (2,1) to[out=270,in=180] (2.25,0.75) to[out=0,in=270] (2.5,1);
\draw[usual] (0,0) to (1,1);
\draw[usual,dot] (2.5,0) to (2.5,0.2);
\draw[usual,dot] (0,1) to (0,0.8);
\draw[usual,dot] (0.5,1) to (0.5,0.8);
\draw[usual,dot] (1.5,1) to (1.5,0.8);
\end{tikzpicture} & YES${}_{c}$
& $\robrmon$ & \begin{tikzpicture}[anchorbase]
\draw[usual] (1,0) to[out=90,in=180] (1.25,0.25) to[out=0,in=90] (1.5,0);
\draw[usual] (1,1) to[out=270,in=180] (1.75,0.55) to[out=0,in=270] (2.5,1);
\draw[usual] (0,0) to (0.5,1);
\draw[usual] (2.5,0) to (2,1);
\draw[usual,dot] (0.5,0) to (0.5,0.2);
\draw[usual,dot] (2,0) to (2,0.2);
\draw[usual,dot] (0,1) to (0,0.8);
\draw[usual,dot] (1.5,1) to (1.5,0.8);
\end{tikzpicture} & YES${}^{\ast}_{c}$
\\
\hline
$\tlmon$ & \begin{tikzpicture}[anchorbase]
\draw[usual] (0.5,0) to[out=90,in=180] (1.25,0.5) to[out=0,in=90] (2,0);
\draw[usual] (1,0) to[out=90,in=180] (1.25,0.25) to[out=0,in=90] (1.5,0);
\draw[usual] (0,1) to[out=270,in=180] (0.25,0.75) to[out=0,in=270] (0.5,1);
\draw[usual] (2,1) to[out=270,in=180] (2.25,0.75) to[out=0,in=270] (2.5,1);
\draw[usual] (0,0) to (1,1);
\draw[usual] (2.5,0) to (1.5,1);
\end{tikzpicture} & YES
& $\brmon$ & \begin{tikzpicture}[anchorbase]
\draw[usual] (0.5,0) to[out=90,in=180] (1.25,0.45) to[out=0,in=90] (2,0);
\draw[usual] (1,0) to[out=90,in=180] (1.25,0.25) to[out=0,in=90] (1.5,0);
\draw[usual] (0,1) to[out=270,in=180] (0.75,0.55) to[out=0,in=270] (1.5,1);
\draw[usual] (1,1) to[out=270,in=180] (1.75,0.55) to[out=0,in=270] (2.5,1);
\draw[usual] (0,0) to (0.5,1);
\draw[usual] (2.5,0) to (2,1);
\end{tikzpicture} & YES${}^{\ast}$
\\
\hline
$\promon$ & \begin{tikzpicture}[anchorbase]
\draw[usual] (0,0) to (0.5,1);
\draw[usual] (0.5,0) to (1,1);
\draw[usual] (2,0) to (1.5,1);
\draw[usual] (2.5,0) to (2.5,1);
\draw[usual,dot] (1,0) to (1,0.2);
\draw[usual,dot] (1.5,0) to (1.5,0.2);
\draw[usual,dot] (0,1) to (0,0.8);
\draw[usual,dot] (2,1) to (2,0.8);
\end{tikzpicture} & YES${}^{\ast}$
& $\romon$ & \begin{tikzpicture}[anchorbase]
\draw[usual] (0,0) to (1,1);
\draw[usual] (0.5,0) to (0,1);
\draw[usual] (2,0) to (2,1);
\draw[usual] (2.5,0) to (0.5,1);
\draw[usual,dot] (1,0) to (1,0.2);
\draw[usual,dot] (1.5,0) to (1.5,0.2);
\draw[usual,dot] (1.5,1) to (1.5,0.8);
\draw[usual,dot] (2.5,1) to (2.5,0.8);
\end{tikzpicture} & YES${}^{\ast}$
\\
\hline
$\psym$ & \begin{tikzpicture}[anchorbase]
\draw[usual] (0,0) to (0,1);
\draw[usual] (0.5,0) to (0.5,1);
\draw[usual] (1,0) to (1,1);
\draw[usual] (1.5,0) to (1.5,1);
\draw[usual] (2,0) to (2,1);
\draw[usual] (2.5,0) to (2.5,1);
\end{tikzpicture} & EX
& $\sym$ & \begin{tikzpicture}[anchorbase]
\draw[usual] (0,0) to (1,1);
\draw[usual] (0.5,0) to (0,1);
\draw[usual] (1,0) to (1.5,1);
\draw[usual] (1.5,0) to (2.5,1);
\draw[usual] (2,0) to (2,1);
\draw[usual] (2.5,0) to (0.5,1);
\end{tikzpicture} & NO
\\
\end{tabular}
.
\end{gather}
The left half of the table above contains \emph{planar} monoids, the right half 
\emph{symmetric} monoids.

We discuss all of these monoids and their 
representation gaps, respectively faithfulness, in \autoref{S:TL} (planar) and 
\autoref{S:Brauer} (symmetric).


\subsection{Further direction not discussed in this paper}\label{SS:IntroFurther}

Although truncated versions of the monoids mentioned in 
\autoref{SS:IntroDiaMonoids} have big representation gaps, big semisimple representation gaps 
and are of high faithfulness, they might not be suitable for cryptographic purposes due to their other properties. 

We list here several additional examples and ideas
which might be interesting to study 
from the perspective of cryptography. For all of these 
making the setup set-theoretical is the first crucial 
(and nontrivial) step:
\begin{enumerate}

\item Web categories in the sense of Kuperberg \cite{Ku-spiders-rank-2}.
These monoidal categories generalize 
the Temperley--Lieb category from the viewpoint of 
representation theory of Lie groups with Temperley--Lieb 
being the $\mathrm{SL}(2)$ case. 

A naive lower bound for the semisimple representation gap of the associated 
endomorphism algebras can be easily obtained. 
This bound is bigger than for the Temperley--Lieb monoid 
itself, so this might be a fruitful direction. 

Note that it is not clear how to make the appearing 
endomorphism algebras set-theoretical. For the Temperley--Lieb 
category what 
one effectively does to make its endomorphism algebras set-theoretical is 
to look at products of light ladders (in the sense 
of \cite{El-ladders-clasps}). The same might work 
for other web categories. Light ladder bases for these 
web categories were discussed 
for example in \cite{AnStTu-cellular-tilting}, \cite{El-ladders-clasps} 
or \cite{Bo-c2-tilting}.

Note that, if one can make these web categories set-theoretical, 
one would get new examples for monoid theory as well, which is interesting in its own right.

\item Soergel bimodules or categorified quantum groups 
in various flavors. 

Soergel bimodules \cite{So-hcbim} form monoidal categories
attached to a Coxeter system. These were diagrammatically reinterpreted in \cite{ElKh-diagrams-soergel} 
and \cite{ElWi-soergel-calculus}, see also
\cite{ElMaThWi-soergel} for a summary. For starters, one can look 
at the dihedral case \cite{El-two-color-soergel} and see whether its set-theoretic modifications can give interesting monoids. Looking at the analogs of light ladders, 
called light leaves in \cite{Li-light-leaves}, might be crucial. 
Let us note that some set-theoretical variations of Soergel diagrammatics
exist in the literature, see for example \cite[Section 4]{CaGoGoSi-algebraic-weaves}, but their 
usefulness in cryptography has  not been explored.

Categorified quantum groups originate in \cite{La-categorification-sl2}, \cite{KhLa-cat-quantum-sln-first}
and \cite{Ro-2-kac-moody}, see also \cite{KhLa-cat-quantum-sln-second}, \cite{KhLa-cat-quantum-sln-third}. These are also diagrammatic in nature 
and promising candidates, but may be harder to work with than Soergel bimodules.

As for web categories, set-theoretical versions of these 
would give novel examples in monoid theory.

\item Foams are suitably decorated 2-dimensional CW-complexes, defined 
abstractly or embedded in $\R^{3}$. They originate and 
most prominently appear 
in the study of link homologies, see for example 
\cite{Kh-sl3-link-homology}, \cite{EhStTu-blanchet-khovanov}, \cite{RoWa-foam-formula} or \cite{EhTuWe-functoriality-link-homologies}.
Using the universal construction from \cite{BlHaMaVo-tqft-kauffman-bracket},
they can easily modified, see {\eg} 
\cite{EhStTu-gl2foams} or 
\cite{KhKi-deform-foam-evaluation}. 

Similarly as in the previous points, if foams could be made 
set-theoretical, that would provide a big supply of potentially interesting monoids.

\item The representation gap and the faithfulness of $\monoid$ 
depend on the underlying field.
To get rid of the dependence of the field, it should be 
useful to consider integral representation of groups or monoids. 
This direction 
is widely open and not much appears to be known. However, their 
categorifications, called $2$-representations, 
have been studied a lot in the recent years.

Potential directions are:

\begin{enumerate}

\item $2$-representations of tensor and fusion categories, 
see {\eg} \cite[Section 7]{EtGeNiOs-tensor-categories} 
for a book chapter discussing these. 
Various diagrammatic fusion categories might be of 
interest to study here, see \cite{MoPeSn-categories-trivalent-vertex} 
for a compelling list of examples. These diagrammatic fusion categories also generalize $\tlmon$, 
so it is expected that the list given in \cite{MoPeSn-categories-trivalent-vertex} 
has suitable big ranks.

\item $2$-representations of fiat $2$-categories, see 
{\eg} \cite{Ma-classification-problems-2reps} for a slightly outdated summary.
For example, Soergel bimodules tend to have simple $2$-representations 
of very big rank, see \cite{MaMaMiTuZh-soergel-2reps} for a classification.
Other versions of $2$-representations of Soergel bimodules might also be useful, see {\eg} \cite{MaTu-soergel} or \cite{MaMaMiTu-trihedral}.

Another advantage of studying $2$-representations of fiat $2$-categories from the viewpoint of cryptography is that cell theory generalizes 
from monoids to these $2$-categories, see {\eg} \cite{MaMaMiTuZh-bireps}, which served as a partial motivation for \autoref{S:Cells}.

\end{enumerate}

\item Another approach is to use semirings for building cryptographic protocols, as proposed in \cite{GrSh-tropical-cryptographyI}, \cite{GrSh-tropical-cryptographyII}, \cite{RaSh-mobs}, 
see also \cite{Du-semirings-cryptography}, which contains a detailed review of the literature.

A linear attack on a semiring-based protocol would require the semiring to act on a vector space or a module, and it is not even clear how a semiring can act linearly on anything. There is the notion of a semimodule over a semiring, which is much closer to set theory compared to that of a module over a ring, and the theory of semimodules over semirings is computationally difficult, even for semimodules over the Boolean semiring $\mathbb{B}=\{0,1|1+1=1\}$, see for example \cite{CoCo-homological-char-one}.
A semiring can appear from a linear structure, as the Grothendieck semiring of an additive category. However, realizing even the Boolean semiring (or the tropical semiring) in this way appears rather nontrivial, due to impossibility of an isomorphism $\one\oplus\one\cong\one$ in a monoidal category, {\cf} \cite{KhTi-cat-zonehalf} which discusses ways to resolve such problems in similar situations.

\end{enumerate}
\medskip

\noindent\textbf{Acknowledgments.}
M.K. was partially supported by NSF grant DMS-1807425 and D.T. was supported by the Australian Research Council while working on the paper. 

D.T. would like to thank Robert Spencer for
freely sharing ideas and computations, patiently answering questions and helpful discussions, and twenty-one green bricks for moral support.
Special thanks to James East, Joel Gibson, Gus Lehrer, Alexander Moret{\'o}, and Geordie 
Williamson for very helpful comments, several clarifications and many useful reference suggestions.
Last but not least, we would like to thank the referee for their 
careful reading of the manuscript, wonderful suggestions, freely sharing ideas and very useful pointers to the literature. The referee's comments improved the mathematical quality and the exposition of the paper a lot, and, indeed, some of the proofs, examples and remarks below are due to the referee.


\section{Representation gaps and faithfulness}\label{S:RepGap}

For background we refer the reader to standard textbooks such as \cite{Be-rep-cohomology}, respectively \cite{St-rep-monoid}, for the basic theory of finite-dimensional representations of finite-dimensional algebras (such as monoid algebras), respectively, finite monoids.

\begin{Notation}\label{N:RepGapMonoid}
We let $\monoid$ denote a finite monoid.
If not stated otherwise, we work over an arbitrary field 
$\K$ and consider only finite-dimensional (left) $\monoid$-representation 
with ground field $\K$. The adjective 
\emph{small} and \emph{big} used for $\monoid$-representations will mean dimension-wise, where 
dimension is measured with respect to $\K$.
\end{Notation}


\subsection{Representation gaps}\label{SS:RepGapRepGap}

We start with a subtle difference between groups and monoids: 
the latter may have two types of ``trivial'' representations.

\begin{Definition}\label{D:RepGapTrivial}
Let $\group\subset\monoid$ be the subgroup of all invertible 
elements of $\monoid$, {\ie} $\group$ is the group of units. Then we define \emph{trivial representations}
\begin{gather*}
\oneb\colon\monoid\to\K,\quad s\mapsto
\begin{cases}
1&\text{if }s\in\group,
\\
0&\text{else},
\end{cases}
\quad
\onet\colon\monoid\to\K,\quad s\mapsto 1.
\end{gather*}
A $\monoid$-representation $\module$ is called \emph{trivial} 
if $\module\cong\oneb$ or $\module\cong\onet$.
\end{Definition}

The subscripts \emph{b} and \emph{t} are short for \emph{bottom} and \emph{top}, respectively. The top trivial representation $\onet$ is also what is called the trivial representation $\one$ of $S$, the unit object of the monoidal category of representations of $S$ with $\one\otimes\module\cong\module$ for any $\monoid$-representation $\module$.

\begin{Remark}\label{R:RepGapTrivial}
The notation is justified as follows. The $\monoid$-representation 
$\oneb$ is one of the simple $\monoid$-representations associated with the bottom $J$-cell $\jb=\group$, while 
the $\monoid$-representation $\onet$ is associated with the top $J$-cell $\jt$, {\cf} \autoref{L:CellsMaximal} below.  
\end{Remark}

\begin{Remark}\label{R:RepGapOrder}
With respect to \autoref{R:RepGapTrivial} and \autoref{S:Cells} below, we warn the reader familiar with monoid theory that the order we use for $J$-cells ({\aka} Green's $J$-classes) is opposite of the one often used in monoid theory. Thus, what we call 
bottom/top is usually the top/bottom in monoid theory. 
In contrast, our convention 
matches most of the cellular algebra literature.
\end{Remark}

\begin{Lemma}\label{L:RepGapTrivial}
Both, $\oneb$ and $\onet$ are simple $\monoid$-representations of 
dimension one. Moreover, $\oneb\cong\onet$ if and only if $\monoid$ is a group.
\end{Lemma}

\begin{proof}
Immediate from the definitions.
\end{proof}

\begin{Notation}\label{N:RepGapTrivial}
We write $\onebt$ short for either $\oneb$ or $\onet$. In particular, 
$\onebt^{\oplus m}$ means any of the $2^{m}$ possible direct sums 
of $\oneb$ and $\onet$ with $m$ symbols in total.
\end{Notation}

For cryptographic purposes it should be 
interesting to collect examples of naturally 
occurring finite monoids $\monoid$ such that 
any representation of sufficiently small dimension 
relative to $|\monoid|$, the size of $\monoid$, is 
suitably trivial. Note that 
all elements of $\msg$ act in the same way on any of the direct 
sums $\onebt^{\oplus m}$ and these 
representations cannot distinguish any two 
elements of $\msg$. Thus, suitably trivial
could mean being isomorphic 
to $\onebt^{\oplus m}$ which we take as the definition.
To state our definition let 
$\onemon$ be the trivial monoid with one element, and 
let $\monoid[S]_{0;1}$ be the monoid on the set $\{a_{0}\}\cup\{1\}$ with unit $1$ and multiplication $a_{0}\cdot a_{0}=a_{0}$ otherwise.

\begin{Definition}\label{D:RepGapmTrivial}
A pair $(\monoid,\K)$ of a monoid, with $\monoid\not\cong\onemon$ and $\monoid\not\cong\monoid[S]_{0;1}$, and a field $\K$ is called \emph{$m$-trivial} if $\monoid$-representations $\module$ with $\dimk(\module)\leq m$ satisfy $\module\cong\onebt^{\oplus\dimk(\module)}$. Moreover, by conventions, $\onemon$ and $\monoid[S]_{0;1}$ are $({-}1)$-trivial for all $\K$.

The maximal $m$ such that $(\monoid,\K)$ is $(m-1)$-trivial is 
called the \emph{representation gap} of $(\monoid,\K)$ and is denoted by $\gap{\monoid}$.
\end{Definition}

\begin{Remark}\label{R:RepGapmTrivialFunnyMonoid}
The two monoids $\onemon$ and $\monoid[S]_{0;1}$ are the only two 
monoids for which every representation is a direct sum of trivial
representations. Hence, their representation gap would be infinity if we would use the same definition for $m$-triviality as for other monoids. 
Since we define $\onemon$ and $\monoid[S]_{0;1}$ 
to be $({-}1)$-trivial we have 
$\gap{\onemon}=\gap{\monoid[S]_{0;1}}=0$.
\end{Remark}

Note that the $m$-triviality is a lower bound on  
the dimension of the smallest nontrivial 
simple $\monoid$-representation, assuming the absence of extensions between trivial representations $\onet$ and $\oneb$,  
see also \autoref{L:RepGapSimple} and \autoref{L:RepGapmTrivial} below.

\begin{Definition}\label{D:RepGapmTrivialMonoid}
A monoid $\monoid$ is called 
\emph{$m$-trivial} if $(\monoid,\K)$ is $m$-trivial for all $\K$. 

The maximal $m$ such that $\monoid$ is $(m-1)$-trivial is 
called the \emph{representation gap} of $\monoid$ and is denoted by $\gap[\ast]{\monoid}$.
\end{Definition}

\begin{Remark}\label{R:RepGapLiterature}
In group theory the representation gap and similar notions are well-known invariants studied by many people and with a number of applications, 
see \cite{BoGa-bounds-cayley-graphs-slfp} or \cite{Go-quasirandom-groups} for examples. However, the motivations in those papers 
are different 
from the ones in this paper and it would be interesting to make a connection.
\end{Remark}

\begin{Notation}\label{D:RepGapmTrivialNotation}
Below we will meet several notions similar to $\gap{\monoid}$ 
and $\gap[\ast]{\monoid}$. For all of them it makes sense to 
vary the field which we indicated using $\ast$.
Whenever the difference does not play a role 
we simply write $\gap[]{\monoid}$.
\end{Notation}

\begin{Remark}[\textbf{\emph{Main Task 1}}]\label{R:RepGapCryptoMain1}
For cryptographic applications it should be useful to have a  
supply of monoids $\{\monoid_{n}|n\in\N\}$ with exponentially big $\gap[]{\monoid_n}$ as $n\to\infty$. 
\end{Remark}

\begin{Example}\label{E:RepGapmTrivial}
A pair $(\monoid\not\cong\onemon,\K)$ or 
$(\monoid\not\cong\monoid_{0,1},\K)$ is $0$-trivial if and only if any there exists a one-dimensional $\monoid$-representation which is nontrivial. 
In particular, if $\monoid$ has a nontrivial one-dimension representation, then $\gap{\monoid}=1$.
\end{Example}

\begin{Lemma}\label{L:RepGapmTrivialForFixedm}
The pair $(\monoid,\K)$ is $m$-trivial if 
and only if $\monoid$-representations 
$\module$ with $\dimk(\module)=m$ satisfy $\module\cong\onebt^{\oplus m}$.
\end{Lemma}

\begin{proof}
By the unique decomposition property of finite-dimensional representations.
\end{proof}

\begin{Lemma}\label{L:RepGapSimple}
Assume that $\monoid$ has at least one nontrivial simple representations. We have 
\begin{gather*}
\gap{\monoid}\leq\min\{\dimk(\simple)|\text{$\simple$ is a nontrivial simple $\monoid$-representation}\}\leq|\monoid|-1.
\end{gather*}
Moreover, when $\K$ is algebraically closed, then $|\monoid|-1$ on the right can be replaced by $\sqrt{|\monoid|-1}$. In all cases, when $\monoid$ is not a group, then every appearance of $|\monoid|-1$ can be replaced by $|\monoid|-2$.
\end{Lemma} 

\begin{proof}
The first inequality follows directly from the definitions. 
To see the second inequality observe 
that simple $\monoid$-representation appear 
in the Jordan--H{\"o}lder filtration of $\K\monoid$, 
the monoid algebra, so their 
dimensions are bounded by $\dimk\big(\K\monoid\big)=|\monoid|$. 
Since the trivial representations $\onebt$ must appear as composition factors we actually get $|\monoid|-2$ or $|\monoid|-1$ 
as an upper bound, depending on whether $\oneb\not\cong\onet$ or not.
When $\K$ is algebraically closed we have the inequality $\sum_{\simple[]}\dimk(\simple[])^{2}\leq|\monoid|$ where the sum runs over all simple $\monoid$-representations. This implies the 
final claim after again taking into account that $\onebt$ must appear as composition factors.
\end{proof}

\begin{Remark}\label{R:RepGapNoNonTrivialModules}
Note that we assume that $\monoid$ has at least one nontrivial simple representations in \autoref{L:RepGapSimple}. 
This restriction is necessary. For example, 
let $\monoid[S]_{0,\dots,n-1;1}$ be the monoid on $\{a_{0},\dots,a_{n-1}\}\cup\{1\}$ 
with unit $1$ and multiplication $a_{i}a_{j}=a_{i}$ otherwise.
Then the only simple $\monoid[S]_{0,\dots,n-1;1}$-representations are
$\onebt$, as follows directly from \autoref{P:CellsSimples} below. Thus, the middle number in \autoref{L:RepGapSimple} is ambiguous.
\end{Remark}

\begin{Example}\label{E:RepGapSym}
Let $\sym=\Aut(\{1,\dots,n\})$ be the \emph{symmetric group} on $\{1,\dots,n\}$. 
For $\cchar\neq 2$ there is a $1$-dimensional nontrivial simple $\sym$-representation, called 
the sign representation. Hence, $\gap{\sym}=1$ unless $\cchar=2$, which implies $\gap[\ast]{\sym}=1$. Since $|\sym|=n!$, the ratio between the representation gap and the 
size of $\sym$ is thus very small. 
Even if one would argue that the sign representation 
is close to trivial, there is still the standard $\sym$-representation 
of dimension $n-1$. So $\gap{\monoid}\leq n-1$ by 
\autoref{L:RepGapSimple}, which is still small compared to $n!$.
\end{Example}

\begin{Example}\label{E:RepGapNoNonTrivialModules}
For the monoid in \autoref{R:RepGapNoNonTrivialModules} 
we have $\gap[\ast]{\monoid[S]_{0;1}}=0$ for $n=1$ and $\gap[\ast]{\monoid[S]_{0,\dots,n-1;1}}=2$ otherwise. This is not hard to verify, see also \autoref{E:RepGapNoNonTrivialModulesTwo} below.
\end{Example}


\subsection{Extensions and representation gaps}\label{SS:RepGapExtensions}

We now discuss extensions. These results are essentially in the literature, but we decided to keep the proofs for convenience of 
the reader. We elaborate on the literature in \autoref{R:RepGapExtensions} below.

We start with an example showing that there can be arbitrary complicated extensions, even with only trivial composition factors:

\begin{Example}\label{E:RepGapNoNonTrivialModulesTwo}
Back to \autoref{E:RepGapNoNonTrivialModules}. 
One can check that $\K\monoid[S]_{0,\dots,n-1;1}$
is a split basic algebra whose quiver $\Gamma$ is of the form
\begin{gather*}
n=1\colon
\Gamma=\bullet\phantom{\rightarrow}\bullet
,\quad
n=2\colon
\Gamma=\bullet\rightarrow\bullet
,\quad
n=3\colon
\Gamma=\bullet\rightrightarrows\bullet
,
\end{gather*}
and so on, {\ie} one has two vertices and $n-1$ edges for $\K\monoid[S]_{0,\dots,n-1;1}$.

Let us use the convention on path algebras where paths are composed from right to left. Then an isomorphism that realizes these descriptions 
sends $a_{0}$ to the initial vertex (on the left-hand side above), $1-a_{0}$ to the terminal vertex 
and, for $n\geq 2$, $a_{i}-a_{0}$ to the $i$th edge, counting {\eg} from top to bottom in the illustration, for $i\in\{1,\dots,n-1\}$.

By usual quiver representation theory it follows that 
$\K\monoid[S]_{0,\dots,n-1;1}$ is semisimple for $n=1$, 
has finite representation type for $n=2$, tame representation type 
for $n=3$ and is of wild representation type for $n\geq 4$.

However, as we have seen in \autoref{E:RepGapNoNonTrivialModules}, $\monoid[S]_{0,\dots,n-1;1}$ 
has only the trivial simple representations $\onebt$ and is $1$-trivial unless $n=1$. Thus, in general, 
$\monoid[S]_{0,\dots,n-1;1}$ has many nontrivial extensions of the form  $0\lra\onebt\lra\module\lra\onebt\lra 0$ with 
only trivial composition factors for $\module$.
\end{Example}

\begin{Lemma}\label{L:RepGapmTrivial}
A pair $(\monoid,\K)$ is $m$-trivial 
if and only if any nontrivial 
simple $\monoid$-representation has dimension at least 
$m+1$ and all extensions $0\lra\onebt\lra\module\lra\onebt\lra 0$ for $\dimk(\module)\leq m$ split.
\end{Lemma} 

\begin{proof} 
Being $m$-trivial clearly implies the second statement. The converse follows by induction on $m$ showing that any $\monoid$-representation $\module$ with $\dimk(\module)\leq m$ is a direct sum of $\onebt$.
\end{proof}

\begin{Remark}\label{R:RepGapmTrivial}
If $\monoid=\group$ is a group so that $\oneb\cong\onet$, then 
having no nontrivial extensions $0\lra\onebt\lra\module\lra\onebt\lra 0$ 
is equivalent to $H^{1}(\monoid,\K)\cong 0$ here $\monoid$ acts on $\K$ trivially: $s\mapsto 1$ for all $s\in\monoid$.
Moreover, for any monoid 
$\monoid$, recall that $H^{1}(\monoid,\K)$
consists of all homomorphisms from $\monoid$ to $(\K,+)$.
In particular, $H^{1}(\monoid,\K)\cong 0$
if and only if the only 
homomorphism from $\monoid$ to $(\K,+)$ is the trivial one. 
We will use this below, in particular, maps from $\monoid$ are always 
to $(\K,+)$.
\end{Remark}

We consider now the four possible cases of extensions of $\onebt$ by $\onebt$. 
Precisely, let $\module$ be an $\monoid$-representation. Suppose there is a short exact sequence
\begin{gather*}
0\lra\onebt\lra\module\lra\onebt\lra 0,\quad\text{meaning all four possibilities.}
\end{gather*}
Choosing a basis of $\module$ compatible with the corresponding filtration, 
the action of each $a\in\monoid$ in the basis will be given by an upper-triangular matrix, 
with either $0$ or $1$ in each diagonal entry (when the corresponding term is either $\oneb$ 
or $\onet$, respectively). The remaining $(1,2)$-entry is denoted by $f(a)$, 
so that the extension is described by a function $f\colon\monoid\to\K$. The condition 
$(ab)m=a(bm)$ for $m\in\module$ translates into 
four possible relations on $f$ depending on the 
types of the trivial representations involved:
\medskip 

\emph{Case (tt).} This case is the same as for groups, {\cf} \autoref{R:RepGapmTrivial}, that is:

\begin{Lemma}\label{L:RepGapTT}
We have $\HH^{1}(\monoid,\K)\cong 0$ if and only if 
$\monoid$ has only the trivial extension of the form
$0\lra\onet\lra\module\lra\onet\lra 0$.
\end{Lemma} 

\begin{proof}
Extensions of the form $0\lra\onet\lra\module\lra\onet\lra 0$, 
viewed as elements of $\Ext^{1}(\onet,\onet)$, are classified by functions $f\colon\monoid\to\K$ 
such that $f(ab)=f(a)+f(b)$ for $a,b\in\monoid$. Any such 
extension is trivial if and only if $\mathrm{H}^{1}(\monoid,\K)\cong 0$.
\end{proof}
\medskip 

\emph{Case (bt).} Recall that $\group\subset\monoid$ denotes the 
group of units of $\monoid$.

Consider the symmetric and transitive closure of the relation 
$ab\approx_{r}a$ for $a,b\in\msg$, and denote the closure by $\approx_{r}$ as well. We call $\monoid$ with a unique equivalence class in $\msg$ under 
$\approx_{r}$ a right-connected monoid.

\begin{Remark}\label{L:RepGapGroupsAreNice}
Note that groups $\monoid=\group$ are not right-connected 
since for groups we have $\msg=\emptyset$, and the empty set
has no equivalence classes under 
$\approx_{r}$.
\end{Remark}

We obtain a sufficient condition for the triviality of extensions:

\begin{Lemma}\label{L:RepGapTB}
If $\monoid$ is right-connected, then
$\monoid$ has only the trivial extension of the form
$0\lra\oneb\lra\module\lra\onet\lra 0$.
\end{Lemma} 

\begin{proof}
Extensions of the form $0\lra\oneb\lra\module\lra\onet\lra 0$,
viewed as elements of $\Ext^{1}(\onet,\oneb)$, are classified by functions $f\colon\monoid\to\K$ such that
\begin{gather}\label{Eq:RepGapEquationFab}
f(ab)= 
\begin{cases}
f(a) 
&\text{if }a\in\msg,
\\
f(a)+f(b) 
&\text{if }a\in\group, 
\end{cases}
\end{gather}
modulo the one-dimensional subspace of functions that are constant on $\msg$ and zero on $\group$. To see this, in a compatible basis $\{v_{1},v_{2}\}$ of $\module$ the action of $a\in\msg$ and $b\in\group$ is given by 
\begin{gather*}
a\mapsto 
\begin{pmatrix}
0 & f(a)\\
0 & 1
\end{pmatrix}
,\quad 
b\mapsto 
\begin{pmatrix}
1 & f(b)\\
0 & 1
\end{pmatrix},
\end{gather*}
leading to the above equations. Moreover, the basis 
$\{v_{1},v_{2}\}$ can be changed to $\{v_{1},v_{2}+\lambda v_{1}\}$ 
while preserving its compatibility with 
the sequence $0\lra\oneb\lra\module\lra\onet\lra 0$, 
explaining why one needs to mod out by functions that are constant on $\msg$ and zero on $\group$. 

If $f$ satisfies \autoref{Eq:RepGapEquationFab}, 
then the fact that $f(ab)=f(a)$ for 
$a\in\msg$ and $b\in\monoid$, together with right-connectedness 
implies that $f$ is constant on $\msg$.
Fix $b\in\msg$ (the set $\msg$ is nonempty by right-connectedness). Then, if $a\in\group$, we have
$ab\in\msg$ and so $f(b)=f(ab)=f(a)+f(b)$, whence $f(a)=0$. Thus, $f$
vanishes on $\group$. We deduce that 
$\Ext^{1}(\onet,\oneb)\cong 0$ by the previous paragraph.
\end{proof}
\medskip

\emph{Case (tb).} A monoid $\monoid$ is called \emph{left-connected} if the opposite monoid $\monoid^{op}$ is right-connected.

\begin{Lemma}\label{L:RepGapBT}
If $\monoid$ is left-connected, then 
$\monoid$ has only the trivial extension of the form
$0\lra\onet\lra\module\lra\oneb\lra 0$.
\end{Lemma} 

\begin{proof}
Dual to \autoref{L:RepGapTB}.
\end{proof} 
\medskip

\emph{Case (bb).} Finally, we call a monoid $\monoid$ \emph{null-connected} if any noninvertible element of $\monoid$ can be written as a product of two noninvertible elements. That is, for $a\in\msg$ we have $a=bc$ for some $b,c\in\msg$. Note that groups are null-connected.

\begin{Lemma}\label{L:RepGapBB}
If $\monoid$ is null-connected and $\HH^{1}(\group,\K)\cong 0$, then 
$\monoid$ has only the trivial extension of the form
$0\lra\oneb\lra\module\lra\oneb\lra 0$.
\end{Lemma} 

\begin{proof}
The extensions as in the statement, when viewed as elements of $\Ext^{1}(\oneb,\oneb)$, are classified by functions $f\colon\monoid\to\K$ such that
\begin{gather*}
f(ab)= 
\begin{cases}
0 
& \text{if }a,b\in\msg,
\\
f(a)+f(b) 
& \text{if }a,b\in\group, 
\\
f(a) 
& \text{if }a\in\msg,b\in\group, 
\\
f(b) 
& \text{if }a\in\group,b\in\msg.
\end{cases}
\end{gather*}
Similarly as before, one can see this by writing the action on $\module$ in a compatible basis as
\begin{gather*}
a\mapsto 
\begin{pmatrix}
0 & f(a)\\
0 & 0
\end{pmatrix}
,\quad 
b\mapsto 
\begin{pmatrix}
1 & f(b)\\
0 & 1
\end{pmatrix},
\end{gather*}
where $a\in\msg$ and $b\in\group$.
The rest of the argument is similar to \autoref{L:RepGapTB} and omitted.
\end{proof}

We say that a monoid $\monoid$ is \emph{well-connected} if it is either a group or right-connected, left-connected and null-connected.

\begin{Theorem}\label{T:RepGapH1Condition}
Assume $\monoid$ is well-connected 
and $\HH^{1}(\group,\K)\cong 0$. Then:
\begin{enumerate}

\item Any short exact sequence 
\begin{gather*}
0\lra\onebt\lra\module\lra\onebt\lra 0
\end{gather*}
splits.

\item We have
\begin{gather}\label{Eq:RepGapMainEquation}
\gap{\monoid}=\min\{\dimk(\simple[])|\text{$\simple[]\not\cong\onebt$ is a simple $\monoid$-representation}\}.
\end{gather} 

\end{enumerate}
In particular, for groups
$\monoid=\group$ it suffices to check whether 
$\HH^{1}(\group,\K)\cong 0$ to ensure that \autoref{Eq:RepGapMainEquation} hold.

Moreover, if $\monoid$ is semisimple over $\K$, then 
$\monoid$ is well-connected and $\HH^{1}(\group,\K)\cong 0$, so (a) 
and (b) hold.
\end{Theorem}

\begin{proof}
\emph{Well-connected and $\HH^{1}(\group,\K)\cong 0$ imply Claim (a).} This claim follows from \autoref{R:RepGapmTrivial},
and the statements in \autoref{L:RepGapTT}, \autoref{L:RepGapBT}, \autoref{L:RepGapTB} and \autoref{L:RepGapBB}.

\emph{Well-connected and $\HH^{1}(\group,\K)\cong 0$ imply Claim (b).} This follows from (a) and the definitions.

\emph{Groups.} Since $\oneb\cong\onet$, \autoref{L:RepGapTT} handles this case. It hence suffices to check $\HH^{1}(\group,\K)\cong 0$ for groups.
\medskip

We now assume that $\monoid$ is semisimple over $\K$.

\textit{Left and right-connectivity.} Assume that $\monoid$ 
is not a group. To see that $\monoid$ is right-connected note that the 
$\monoid$-representation $\onet$ is projective. Thus, there exists 
$e\in\monoid$ with $e^{2}=e$ and $\onet\cong\K\monoid e$.
Let $a,b\in\monoid$ with $a$ in the support of $e$ and $b\in\msg$. Then, since $be=e$, 
we get that there exists $c\in\monoid$ with $a=bc$ so $\monoid$ is right-connected. Finally, taking the opposite 
monoid preserves semisimplicity, so the same arguments as for right-connectivity imply left-connectivity.

\textit{Null-connectivity.} Recall that ideals 
in semisimple algebras are (unital) semisimple algebras. Hence, $\K(\msg)/(\msg)^{2}$ is semisimple, so it cannot be nilpotent. This implies that $(\msg)=(\msg)^{2}$, and thus, $\monoid$ is well-connected.

\emph{The cohomology vanishes.} The surjection $\K\monoid\twoheadrightarrow\K\group$ given by $a\mapsto a$ for $a\in\group$ and $a\mapsto 0$ for $a\in\msg$ implies that $\group$ is semisimple over $\K$ if $\monoid$ is. Thus, we get $\HH^{1}(\group,\K)\cong 0$.
\end{proof}

\begin{Remark}\label{R:RepGapUpperLowerExtensions}
Note that for upper bounds for $\gap[]{\monoid}$ it suffices to 
find some nontrivial simple $\monoid$-representation, but for lower bounds 
or the explicit value of $\gap[]{\monoid}$ we will calculate 
$\HH^{1}(\monoid,\K)$ and $\HH^{1}(\group,\K)$.
\end{Remark}

\begin{Remark}\label{R:RepGapExtensions}
The paper \cite{MaSt-quiver-monoids} computes 
certain quivers for monoid algebras with the computation of 
a generalization of $\Ext^{1}(\onet,\oneb)$ being a main point. 
The above lemmas are deducible from their computations, more precisely 
from \cite[Section 7]{MaSt-quiver-monoids}. In fact, \cite[Section 7]{MaSt-quiver-monoids} work in much greater 
generality and the setting with $\onet$ and $\oneb$ is a very special case.
\end{Remark}

\begin{Remark}\label{R:RepGapExtensionsTwo}
Using ideas in \cite{MaSt-quiver-monoids}, one can get a description of 
$\Ext^{1}(\onet,\oneb)$ as in the proof of \autoref{L:RepGapTB}.
That is, one can prove that $\Ext^{1}(\onet,\oneb)\cong\widetilde{\mathrm{H}}^{0}
(\Delta(P_{r}),\K)^{\group}$
(reduced cohomology) where $P_{r}$ is the poset of proper principal
right ideals of $\monoid$ and $\Delta(P_{r})$ is its order complex. There is, of course, the dual version for $\Ext^{1}(\oneb,\onet)$ using proper 
principal left ideals of $\monoid$.
Let us also mention that the special case of this result where 
$\group$ is trivial was explicitly proved in 
\cite{MaSaSt-combinatorial-topology} 
and a different proof was given in \cite{MaSaSt-cell-complexes} for left regular bands.

Similarly, following the ideas in \cite{MaSt-quiver-monoids}, 
one can show that $\Ext^{1}(\onet,\onet)\cong\mathrm{H}^{1}(\group,\K)\oplus\K^{|\group\slash A\backslash\group|}$ where $A=(\msg)/(\msg)^{2}$.
\end{Remark}

Recall for monoid theory 
that $a\in\monoid$ is called \emph{von Neumann regular}
if it can be written as $a=aba$ for some $b\in\monoid$, 
and $\monoid$ is \emph{von Neumann regular} if all of its elements are. Examples of von Neumann regular monoids are the diagram monoids 
in \autoref{Eq:IntroDiaMonoids}. 
As a final statement in this section we add:

\begin{Lemma}\label{R:RepGapExtensionsRegular}
If $\monoid$ is von Neumann regular, then $\monoid$ is null-connected.
\end{Lemma}

\begin{proof}
Any $a\in\msg$ satisfies $a=aba$ for some $b\in\monoid$. Since 
$ba\in\msg$ whenever $a\in\msg$, null-connectivity follows.
\end{proof}


\subsection{Examples}\label{SS:RepGapExamples}

The following is well-known. But since it is an important 
example for cryptography, see {\eg} 
\autoref{E:RepGapCyclic} below, we state and prove it here.

\begin{Proposition}\label{P:RepGapCyclicGroups}
Let $\cyclic\cong\Z/n\Z$ be the \emph{cyclic group} of order $n>1$.
\begin{enumerate}

\item We have $\gap[\Q]{\cyclic}=\min\{r-1|r\text{ prime},\ r\mid n\}$. (In particular, $\gap[\Q]{\cyclic}=n-1$ if $n$ is prime.)

\item Let $\F_{q}$ denote a finite field with $q=p^{k}$ elements, where $p$ is a prime.
\begin{enumerate}

\item For $\gcd(n,q-1)>1$ we have $\gap[\F_{q}]{\cyclic}=1$.

\item For $\gcd(n,q-1)=1$ and $p\mid n$ we have
$\gap[\F_{q}]{\cyclic}=2$.

\item For $\gcd(n,q-1)=1$ and $p\nmid n$ 
we have $\gap[\F_{q}]{\cyclic}=\min\big\{d\in\N|\gcd(n,q^{d}-1)\neq 1\big\}$.

\end{enumerate}
\item For any field $\K$ we have $\gap[\K]{\cyclic}=\min_{r}(\gap[\K]{\cyclic[r]})$, where the minimum is taken over all prime divisors $r$ of $n$.
\end{enumerate}
\end{Proposition}

\begin{proof}
\emph{Case (a).} First we have $\HH^{1}(\cyclic,\Q)\cong 0$, so 
by \autoref{T:RepGapH1Condition} it suffices 
to look at the dimensions of simple $\cyclic$-representations. 

To this end, recall that representations of $\cyclic$ are semisimple over $\Q$. The polynomial $X^{n}-1$ has no repeated roots over $\Q$ and factors
as $X^{n}-1=\prod_{d\mid n}\Phi_{d}$ for $\Phi_{d}$ the $d$th cyclotomic polynomial.
The Chinese reminder theorem then
gives $\Q[\cyclic]\cong\bigoplus_{d\mid n}\Q[X]/(\Phi_{d})$, and 
we see that there are simple $\cyclic$-representations for each $\Phi_{d}$ 
which are of the respective degrees $\deg\Phi_{d}=\varphi(d)$.
This implies $\gap[\Q]{\cyclic}=\min\{\varphi(d)|d\text{ divides }n\}$. However, 
since $a\mid b$ implies $\varphi(a)\mid\varphi(b)$ we get the claimed 
formula from this expression.
\medskip

\emph{Case (b).} There is a nontrivial one-dimensional $\cyclic$-representation
over $\F_{q}$ exactly when $\gcd(n,q-1)>1$, implying (i). In case (ii), there exist a nontrivial homomorphism $\cyclic\to\F_{q}$, where the latter is considered an abelian group under addition, giving a nontrivial selfextension of the trivial representation of $\cyclic$. 

In the remaining case (iii), when $\gcd(n,q-1)=1$ and $p\nmid n$, the trivial representation has no selfextensions and it is the unique (up to isomorphism) representation of dimension one over $\F_{p}$. The representation gap $\gap[\F_{p}]{\cyclic}$ is then the dimension $d\geq 2$ of the smallest nontrivial simple representation. Such a representation corresponds to a nontrivial homomorphism $\cyclic\to\mathrm{GL}(d,\F_{q})$. Since $\gcd(n,q-1)=1$ this homomorphism does not take $\cyclic$ to multiples of the identity matrix. So $d$ is the smallest number 
such that $\gcd\big(n,|\mathrm{GL}(d,\F_{q})|\big)\neq 1$. 
The order of $\mathrm{GL}(d,\F_{q})$, up to factors of $q-1$, which are coprime to $n$, is $(q^{d}-1)(q^{d}-q)\dots (q^{d}-q^{d-1})$. We see that the smallest $d$ with $\gcd\big(n,|\mathrm{GL}(d,\F_{q})|\big)\neq 1$ is the smallest $d$ such that $gcd(n,q^{d}-1)\neq 1$.
\medskip

\emph{Case (c).} This follows from (a) and (b).
\end{proof}

\begin{Example}\label{E:RepGapCyclic}
The groups $\cyclic$ lie at the heart of many standard cryptographic protocols, see {\eg} \cite[Section 1.4]{Ko-algebraic-cryptography}.
By \autoref{P:RepGapCyclicGroups} these groups have a quite big representation gap over $\Q$. However, the situation varies depending on the ground field, and over $\C$ the representation gap is small. In particular, 
for cryptographical purposes the point is that protocols are broken 
as soon as $\cyclic$ is identified explicitly. For $n+1=p$ with $p$ a large prime the classical protocols ``disguise'' 
$\cyclic$ since finding a generator of $(\Z/p\Z)^{\ast}$, meaning finding an explicit isomorphism of groups 
$(\Z/p\Z)^{\ast}\cong\cyclic$, is difficult.

Let $n$ be a prime number. 
Over a characteristic zero field $\K$ that contains a primitive root of unity $\xi$ of order $r$, all simple $\cyclic$-representations are one-dimensional, and $\gap[\K]{\cyclic}=1$. Instead, as argued in the proof 
of \autoref{P:RepGapCyclicGroups}, over the prime field $\Q$ there are two simple $\cyclic$-representations: the trivial $\one$ and an $(n-1)$-dimensional representation $\module$, the complement of the trivial in the regular representation. The representation $\module$ over a larger field that contains $\xi$ splits into the direct sum of one-dimensional $\cyclic$-representations, which are Galois conjugates of each other.
\end{Example}

Thus, for $n$ prime \autoref{P:RepGapCyclicGroups} 
and \autoref{E:RepGapCyclic} imply that $\cyclic$ has a substantial representation gap $n-1$ over $\Q$, close to its cardinality $n=|\cyclic|$.

\begin{Example}\label{E:RepGapSimpleGroups}
\autoref{P:RepGapCyclicGroups} discusses the cyclic groups $\cyclic$. These are simple if $n$ is a prime and the only commutative groups among the finite simple groups. 

Let us briefly discuss other finite simple groups:
\begin{enumerate}

\item The alternating groups $\alt\subset\sym$ of size 
$\tfrac{n!}{2}$ behave similarly to the symmetric groups, {\cf} \autoref{E:RepGapSym}.
They are a bit better in the sense that they do not have a sign representation. However, 
over $\Q$ the standard representation of $\sym$ restricts to a simple $\alt$-representations. Over other 
fields this representation might not be simple. But if its not, then it contains an even smaller nontrivial simple 
in its Jordan--H{\"o}lder filtration. Hence, $\gap[\ast]{\alt}\leq n-1$.

\item The biggest part of the periodic table of simple groups are the finite groups of Lie type.
(We consider the family of finite groups of Lie type in a very vague sense. 
In fact, the symmetric groups are secretly also part of this family, using the 
analogy that $\sym\leftrightsquigarrow\GL{\F_{1}}$.)
Most of these should have small representation gap over the defining field. 
To see this consider the group $\PSL{\F_{q}}$ for 
$q=p^{k}$ and $p$ a prime. This is a finite simple group 
(unless $n=2$ and $q\in\{2,3\}$) with $\tfrac{q^{n(n-1)/2}}{\gcd(n,q-1)}\prod_{i=2}^{n}(q^{i}-1)$ elements. 
(Thus, the number of elements grows exponentially in $n$.) However, $\PSL{\F_{q}}$ has a small nontrivial simple $\F_{q}$-representation 
of dimension $n^{2}-1$, namely $\big(\F_{q}^{n}\otimes(\F_{q}^{n})^{\ast}\big)/\F_{q}$.

\item Sporadic simple groups tend to have big representation gaps, see {\eg} \cite{CoCuNoPaWi-atlas}. However, 
they do not come in $\N$-families and are all only moderately big. 
So they are probably not of immediate use for cryptography. 

Let us discuss the monster group $M$ 
as an example. Its smallest nontrivial and faithful representation over $\C$ has dimension $196883$, see \cite{CoCuNoPaWi-atlas} under the entry $M=F_{1}$ therein (see also 
\cite[Chapter 12]{FrLeMe-vertex-monster} where this number$+1$ appears as the graded dimension of the moonshine representation), and the smallest nontrivial and faithful representation over any field has dimension $196882$, see \cite{LiPaWaWi-computer-monster}. With the minimal representation of a sufficiently big dimension, there is a potential chance for cryptographic protocols built from the monster. However, the monster still is sporadic and does not 
come in an infinite family. We are not aware of any literature on the subject.

\end{enumerate}
Thus, one could argue that noncommutative finite groups do not seem to be very useful for cryptography purposes 
by the above.
\end{Example}

\begin{Example}\label{E:RepGapPGroups}
Finite groups that often have a big representation gap 
are $p$-groups for a prime $p$. Under the name minimal 
character degree, there is a big literature on 
the representation gap of these 
groups, see for example \cite{Hu-char-degree-pgroups} or \cite{JaZaMo-char-degrees-pgroups}, often aiming for an upper bound 
and not a lower bound as we would need it. Having a large representation 
gap might make them useful in cryptography, see {\eg} 
\cite[Section 3]{Ro-algebraic-cryptography}.
\end{Example}


\subsection{Field size and representation gap}\label{SS:RepGapFieldSize}

In our definition of the representation gap we do not differentiate between a particular field used and our measure of complexity is the dimension of the smallest nontrivial representation over that field. More practically, we can keep track of the complexity of working over a specific field. 

For the finite field $\F_{q}$ a natural measure of complexity is $\log_{2}(\F_{q})=n\log_{2}(p)$, the log of the size of the field or some related complexity that measures the difficulty of manipulating elements of the field. Given an $\monoid$-representation $\module$ over $\F_{q}$, the \emph{complexity} of $\module$ 
over $\F_{q}$ can then be defined as 
\begin{gather*}
c(\module)=\dimk[\F_{q}](\module)c(\F_{q}),
\quad\text{where }
c(\F_{q})=\log_{2}|\F_{q}|.
\end{gather*}
Note that $c(\module)$ is preserved when viewing $\module$ as an $\monoid$-representation over any subfield of $\F_{q}$. 

\begin{Definition}\label{D:RepGapFieldSize}
Define the \emph{finite characteristic representation gap}
$\gap[f]{\monoid}$ of $\monoid$ as the minimum of $c(\module)$, over all nontrivial representations $\module$ over finite fields.
\end{Definition}

We can alternatively restrict to 
$\monoid$-representations $\module$ over finite 
extensions $\Q\subset\K$ and define
\begin{gather*}
c_{0}(\module)= 
\dimk(\module)[\K:\Q]= 
\dimk[\Q](\module).
\end{gather*}
Again, $c_{0}(\module)$ does not change if $\module$ is viewed as an $\monoid$-representation over a subfield $\LL\subset\K$. 

\begin{Definition}\label{D:RepGapFieldSizeTwo}
Define the \emph{characteristic zero representation gap} $\gap[0]{\monoid}$ of $\monoid$ as the minimum of $c_{0}(\module)$, over all nontrivial $\monoid$-representations over finite extensions of $\Q$.
\end{Definition}

The pair $\big(\gap[0]{\monoid},\gap[f]{\monoid}\big)$ is a measure of the representation complexity of $\monoid$ over both $\Q$ and finite fields.

\begin{Remark}\label{R:RepGapFieldSize}
Recall from above that the groups $\cyclic$ have large (exponential) representation gap over $\Q$. The more refined notion of representation gap, introduced in this section, might be a better measure of the complexity of $\monoid$ from the linear attacks viewpoint.
\end{Remark}


\subsection{Faithfulness}\label{SS:RepGapFaithful}

By a \emph{faithful} $\monoid$-representation we mean a representation on which any two elements of $\monoid$ act differently.

\begin{Remark}\label{R:RepGapFaithful}
Since there is no $\K$-linear structure involved, this notion of faithfulness is slightly different from that of a faithful representation of the monoid algebra $\K\monoid$.
\end{Remark}

Besides the notion of the representation gap, we introduce a related 
(weaker) notion: 

\begin{Definition}\label{D:RepGapFaithful}
Let $\faith{\monoid}$ be the number
\begin{gather*}
\faith{\monoid}=\min\{\dimk(\module)|\text{$\module$ is a faithful $\monoid$-representation}\}.
\end{gather*}
We call $\faith{\monoid}$ the \emph{faithfulness} of $(\monoid,\K)$. 
We also define $\faith[\ast]{\monoid}$ to be the minimum of 
$\faith{\monoid}$ over all fields.
\end{Definition}

In words, $\faith[]{\monoid}$ is the dimension of 
the smallest faithful $\monoid$-representation.

\begin{Remark}[\textbf{\emph{Main Task 2}}]\label{R:RepGapCryptoMain2}
Similarly as in \autoref{R:RepGapCryptoMain1},
for cryptographic applications it should be useful to have a  
supply of monoids with exponentially big $\faith{\monoid}$. 
\end{Remark}

\begin{Remark}\label{R:RepGapFaithfulLiterauture}
For finite groups $\faith[]{\monoid}$ is a well-known invariant 
studied since the early days of representation theory. It is
sometimes called representation dimension, and
has attracted recent attention, see \cite{Mo-faithful-finite-group} 
and the references therein, including \cite{CeKaRe-representation-dimension} or \cite{BaMaKaSa-faithful-p-groups}. Various versions 
of faithfulness have been studied in monoid theory 
as well, see for example \cite{MaSt-effective-dimension-semigroups} who call the faithfulness the effective dimension.
\end{Remark}

\begin{Remark}\label{R:RepGapmBurau}
Faithfulness is only one measure of the complexity of $\monoid$.
As one example of a small size representation that is not faithful in general but still gives rise to efficient attacks is the Burau representation of the braid group $\brgr$ on $n$ strands. (The braid group is not a finite monoid, but that does not play a role for our discussions involving it.) 
The Burau representation has dimension $n$, or $n-1$ for 
the reduced Burau representation, and in the proposed protocols $n$ is very small. Furthermore, the kernel of the Burau representation is also small, in an appropriate sense, and the action of an element of $\brgr$ on the representation carries full information about the element for the protocol's purposes. Many of these protocols admit efficient attacks, as documented in the literature.
\end{Remark}

\begin{Example}\label{E:RepGapFaithful1}
The symmetric group $\sym$ has its $n$-dimensional 
permutation representation, which is faithful. Hence,  
$\faith[\ast]{\sym}\leq n$. 

In fact, one can do better. 
If the characteristic of $\K$ does not divide $n$, then $\faith[\K]{\sym}=n-1$. The corresponding $\sym$-representation is the 
standard representation. Otherwise and if $n\geq 5$ one has $\faith[\K]{\sym}=n-2$, and hence, still assuming $n\geq 5$, 
we have $\faith[\ast]{\sym}=n-2$. This is a fact from the early days of representation theory, see {\eg} \cite[Section 9.3]{MaSt-effective-dimension-semigroups} for a modern formulation.
\end{Example}

\begin{Example}\label{E:RepGapFaithful1Cyclic}
We have $\faith[\C]{\cyclic[i,p]}=i+1$ for the cyclic 
monoid that we will meet in \autoref{E:CellsGenOneElement} below, see {\eg}
\cite[Section 10]{MaSt-effective-dimension-semigroups} 
where the author's list $\faith[\C]{\monoid}$ for various 
monoids, including the cyclic ones.
\end{Example}

\begin{Lemma}\label{L:RepGapFaithful}
Assume that $\monoid$ has at least one nontrivial simple representations. Then we have
\begin{gather*}
\gap{\monoid}\leq\faith{\monoid}\leq|\monoid|.
\end{gather*}
\end{Lemma}

\begin{proof}
Every $S$-representation has a Jordan--H{\"o}lder
filtration by simple representations, which therefore are 
of smaller (or equal) dimensions. The first claim then follows from \autoref{L:RepGapSimple}. The second inequality follows since 
every monoid admit a faithful representation on itself.
\end{proof} 

\begin{Remark}\label{R:RepGapNoNonTrivialModulesAgain}
The assumption in 
\autoref{L:RepGapFaithful} is necessary for the same reasons as in \autoref{R:RepGapNoNonTrivialModules}.
\end{Remark}

\begin{Example}\label{E:RepGapFaithful2}
Let $\brgr$ be the braid group on $n$ strands.
We already mentioned its Burau representation in \autoref{R:RepGapmBurau}, but this representation is not faithful in general. 
However, a faithful $\brgr$-representation over $\Q(q,t)$ is the 
Laurence--Krammer--Bigelow representation, see \cite{Bi-linear-artin}
and \cite{Kr-linear-artin}, which is 
of dimension $\tfrac{n(n-1)}{2}$. Thus, 
$\gap[\Q(q,t)]{\brgr}\leq\faith[\Q(q,t)]{\brgr}\leq\tfrac{n(n-1)}{2}$,
which creates obstacles of applications of $\brgr$ to cryptography, see also \cite{MyShUs-attack-braid}.
\end{Example}

The following is useful in examples:

\begin{Lemma}\label{L:RepGapFaithfulEmbedding}
Assume that there is a embedding of monoids
$\monoid\hookrightarrow\monoid[T]$.
\begin{gather*}
\faith[]{\monoid}\leq\faith[]{\monoid[T]}.
\end{gather*}
\end{Lemma}

\begin{proof}
This follows since a faithful $\monoid[T]$-representation 
restricts to a faithful $\monoid$-representation.
\end{proof}

We come back to \autoref{E:RepGapCyclic}, but now from the viewpoint of 
faithfulness.

\begin{Proposition}\label{P:RepGapCyclicGroups2}
Let us consider the setting of \autoref{P:RepGapCyclicGroups}.
\begin{enumerate}

\item We have $\faith[\Q]{\cyclic}=\sum_{i=1}^{k}(r_{i}^{d_{i}}-r_{i}^{d_{i}-1})$, where $n$ has the prime factor decomposition $n=\prod_{i=1}^{k}r_{i}^{d_{i}}$. (In particular, $\faith[\Q]{\cyclic}=n-1$ if $n$ is prime.)

\item Let $n$ be prime and $\cchar\nmid n$. Then $\faith[\K]{\cyclic}=\gap[\K]{\cyclic}$
for all the cases in \autoref{P:RepGapCyclicGroups}.

\end{enumerate}
\end{Proposition}

\begin{proof}
\emph{Case (a).} Recall that $\Q[\cyclic]\cong\bigoplus_{d\mid n}\Q[X]/(\Phi_{d})$, see the proof of \autoref{P:RepGapCyclicGroups}. 
The simple $\cyclic$-representations $\Q[X]/(\Phi_{d})$ can be identified with 
$\Q(\zeta_{d})$ for $\zeta_{d}$ a primitive $d$th root of unity. It is then 
easy to see that $\bigoplus_{i=1}^{k}\Q(\zeta_{r(i)})$ for $r(i)=r_{i}^{d_{i}}$ is a faithful $\cyclic$-representation. The dimensions of the summands are 
the degrees of the associated $\Phi_{d}$. Hence, these
summands are of dimensions 
$r_{i}^{d_{i}}-r_{i}^{d_{i}-1}$, which shows 
$\faith[\Q]{\cyclic}\leq\sum_{i=1}^{k}(r_{i}^{d_{i}}-r_{i}^{d_{i}-1})$.
The decomposition of $\Q[\cyclic]$ into $\Q(\zeta_{d})$ also implies that 
one can not find a smaller faithful $\cyclic$-representation since $\Phi_{d}$ with $d=kr(i)$ and $k$ coprime to $r_{i}$ has bigger degree than 
$\Phi_{r(i)}$.
\medskip

\emph{Case (b).} This follows since $\cyclic$ is a simple group when $n$ is a prime, and because the representation theory of $\cyclic$ is semisimple under the assumption $\cchar\nmid n$.
\end{proof}

The analog of \autoref{E:RepGapSimpleGroups} is:

\begin{Example}\label{E:RepGapSimpleGroups2}
For finite simple groups faithfulness
is not much different from \autoref{E:RepGapSimpleGroups}. That is, 
\autoref{P:RepGapCyclicGroups2} treats the cyclic groups and:
\begin{enumerate}

\item The alternating groups $\alt$ has a faithful representation 
of dimension $n$, which is the restriction of the permutation 
representation of $\sym$ to $\alt$, see also 
\autoref{L:RepGapFaithfulEmbedding}. Thus, $\faith[\ast]{\alt}\leq n$.

\item The $\GL{\F_{q}}$-representation $\F_{q}^{n}$ is faithful, giving an example of a group acting faithfully on a small representation. To pass to a simple group, one can take $\PSL{\F_{q}}$, which then acts faithfully on $\F_{q}^{n}\otimes(\F_{q}^{n})^{\ast}$. Hence, $\faithb[\F_{q}]{\PSL{\F_{q}}}\leq n^{2}$.

\item For sporadic groups
the same remarks as in \autoref{E:RepGapSimpleGroups} apply.
The smallest faithful representations for sporadic groups 
are listed in \cite{Ja-minimal-faithful}.

\end{enumerate}
\autoref{E:RepGapSimpleGroups} and this example motivate 
to study monoids that are not groups.
\end{Example}

\begin{Example}\label{E:RepGapPGroupsFaith}
Similarly as in \autoref{E:RepGapPGroups}, $p$-groups tend to have a large faithfulness and this is well-studied, see {\eg} 
\cite{Ja-faithful-pgroups} for some early results 
and \cite{Mo-faithful-finite-group} for a more recent treatment.
\end{Example}


\subsection{Ratios}\label{SS:RepGapRatios}

As argued earlier, for potential 
cryptographic purposes one wants to specialize to 
monoids with the representation gap of 
size comparable to $|\monoid|^{\epsilon}$, 
for some $\epsilon>0$, as opposed to monoids 
where representation gap is exponentially 
smaller than the size of $\monoid$. As a 
measure of complexity, we 
can define:

\begin{Definition}\label{D:RepGapRatio}
The \emph{gap-ratio} and the \emph{faithful-ratio} of $\monoid$ are
\begin{gather}\label{Eq:RepGapRatio}
\gratio{\monoid}=\frac{\gap{\monoid}}{\sqrt{|\monoid|}},
\quad
\fratio{\monoid}=\frac{\faith{\monoid}}{|\monoid|}.
\end{gather}
\end{Definition}

\begin{Remark}[\textbf{\emph{Additional Task 1}}]\label{R:RepGapCryptoAdditional}
For cryptographic applications it makes sense to search for naturally occurring 
families of monoids $\{\monoid_{n}|n\in\N\}$ with 
$\lim_{n\to\infty}\gratio{\monoid_{n}}$ or 
$\lim_{n\to\infty}\fratio{\monoid_{n}}$ that do not approach $0$ exponentially fast. 

Note that these are rather crude: They are motivated 
by the search for families of monoids $\{\monoid_{n}|n\in\N\}$ where representation gap grows exponentially while computations in the monoid grow polynomial, but oversimplify this problem.
\end{Remark} 

\begin{Remark}\label{R:RepGapCrypto2}
The square root in \autoref{Eq:RepGapRatio} comes from the 
observation that over an algebraically closed field a simple $\monoid$-representation has dimension at most $\sqrt{|\monoid|}$.
We stress that we have a slightly better bound of $\sqrt{|\monoid|-1}$ or $\sqrt{|\monoid|-2}$ in \autoref{Eq:RepGapRatio}, but the differences to $\sqrt{|\monoid|}$ do not 
play significant roles so we ignored these bounds in \autoref{Eq:RepGapRatio} for the sake of simplicity.
\end{Remark}

\begin{Example}\label{R:RepGapGapRatio}
For the symmetric group $\sym$, {\cf} \autoref{E:RepGapSym} and 
\autoref{E:RepGapFaithful1}, we have $\gratio[\ast]{\sym}=\big(\sqrt{n!}\big)^{-1}$ and $\fratio[\ast]{\sym}=\big((n-3)!(n-1)n\big)^{-1}$, again indicating 
that $\sym$ is not very useful for cryptography. 
The alternating group as in \autoref{E:RepGapSimpleGroups} and \autoref{E:RepGapSimpleGroups2} 
has $\gratio[\ast]{\alt}\leq\sqrt{2}(n-1)\big(\sqrt{n!}\big)^{-1}$ and $\fratio[\ast]{\alt}\leq2\big((n-1)!\big)^{-1}$, which are still tiny.
\end{Example}

\begin{Example}\label{R:RepGapGapRatioTwo}
For monoids it is not hard to find examples 
with $\fratio[\ast]{\monoid}=1$, see \cite[Proposition 28]{MaSt-effective-dimension-semigroups} for an explicit example.
Moreover, the main monoids under study in this paper 
have also large $\gratio[\ast]{\monoid}$, see {\eg} \autoref{T:TLIsAGoodExample} below.
\end{Example}


\section{Cell theory}\label{S:Cells}

An important tool to study representations of monoids 
are \emph{Green cells} or \emph{Green's relations}. 
In this section we explain how these help to calculate 
$\gap[]{\monoid}$ and $\faith[]{\monoid}$, and also 
give us another numerical measure which we will 
call \emph{semisimple representation gap}.

\begin{Remark}\label{R:CellsCells}
We will summarize the 
main constructions using the language of cells 
as in \cite{GrLe-cellular}, which is 
more common in representation theory. The classical description
using Green's relations from monoid theory can be found in 
many (older and newer) papers {\eg} \cite{Gr-structure-semigroups} or \cite{GaMaSt-irreps-semigroups}, and also in books such as 
\cite{ClPr-algebraic-semigroups}, \cite{ClPr-algebraic-semigroups-2} or \cite{St-rep-monoid}. The cell based discussion is not so easy to find in the 
literature, see however \cite{GuWi-almost-cellular}, \cite{TuVa-handlebody} or \cite{Tu-sandwich}.
\end{Remark}


\subsection{The basics}\label{SS:CellsBasics}

Recall that $\monoid$ denotes a finite monoid. (Cell theory also works for infinite monoids, but the theory is technically more involved. We will not discuss it here.)

We define preorders
on $\monoid$, called \emph{left, right and two-sided cell order}, by
\begin{align*}
(a\leq_{l}b)&\Leftrightarrow
\exists c:b=ca,
\\
(a\leq_{r}b)&\Leftrightarrow
\exists c:b=ac,
\\
(a\leq_{lr}b)&\Leftrightarrow
\exists c,d:b=cad.
\end{align*}
In words, $a$ is left lower than $b$ if $b$ can be obtained from $a$ 
by left multiplication, and similarly 
for right ans two-sided. 

\begin{Remark}\label{R:CellsOrder}
As in \autoref{R:RepGapOrder}, these orders are
in-line with the most common convention used in the theory of cellular algebras but the opposite of the one usually used in monoid theory.
\end{Remark}

We define equivalence relations, 
the \emph{left, right and two-sided equivalence}, by
\begin{align*}
(a\sim_{l}b)&\Leftrightarrow
(a\leq_{l}b\text{ and }b\leq_{l}a),
\\
(a\sim_{r}b)&\Leftrightarrow
(a\leq_{r}b\text{ and }b\leq_{r}a),
\\
(a\sim_{lr}b)&\Leftrightarrow
(a\leq_{lr}b\text{ and }b\leq_{lr}a).
\end{align*}
The respective equivalence classes are called 
left, right respectively two-sided \emph{cells}. 
We denote all these by $\lcell$, $\rcell$ and $\jcell$
and call two-sided cells \emph{$J$-cells}. 
Finally, an \emph{$H$-cell} $\hcell=\hcell(\lcell,\rcell)=\lcell\cap\rcell$ is an intersection of a 
left $\lcell$ and a right cell $\rcell$. 

The picture to keep in mind (stolen from \cite[Section 2]{TuVa-handlebody}) is
\begin{gather}\label{Eq:CellsIllustration}
\begin{tikzpicture}[baseline=(A.center),every node/.style=
{anchor=base,minimum width=1.4cm,minimum height=1cm}]
\matrix (A) [matrix of math nodes,ampersand replacement=\&] 
{
\hcell_{11} \& \hcell_{12} 
\& \hcell_{13} \& \hcell_{14}
\\
\hcell_{21} \& \hcell_{22} 
\& \hcell_{23} \& \hcell_{24}
\\
\hcell_{31} \& \hcell_{32} 
\& \hcell_{33} \& \hcell_{34}
\\
};
\draw[fill=blue,opacity=0.25] (A-3-1.north west) node[blue,left,xshift=0.15cm,yshift=-0.5cm,opacity=1] 
{$\rcell$} to (A-3-4.north east) 
to (A-3-4.south east) to (A-3-1.south west) to (A-3-1.north west);
\draw[fill=red,opacity=0.25] (A-1-3.north west) node[red,above,xshift=0.7cm,opacity=1] 
{$\lcell$} to (A-3-3.south west) 
to (A-3-3.south east) to (A-1-3.north east) to (A-1-3.north west);
\draw[very thick,black,dotted] (A-1-2.north west) to (A-3-2.south west);
\draw[very thick,black,dotted] (A-1-3.north west) to (A-3-3.south west);
\draw[very thick,black,dotted] (A-1-4.north west) to (A-3-4.south west);
\draw[very thick,black,dotted] (A-2-1.north west) to (A-2-4.north east);
\draw[very thick,black,dotted] (A-3-1.north west) to (A-3-4.north east);
\draw[very thick,black] (A-1-1.north west) node[black,above,xshift=-0.5cm] {$\jcell$} to 
(A-1-4.north east) to (A-3-4.south east) 
to (A-3-1.south west) to (A-1-1.north west);
\draw[very thick,black,->] ($(A-1-1.north west)+(-0.4,0.4)$) to (A-1-1.north west);
\draw[very thick,black,->] ($(A-3-4.south east)+(0.5,0.2)$) 
node[right]{$\hcell(\lcell,\rcell)=\hcell_{33}$} 
to[out=180,in=0] ($(A-3-3.south east)+(0,0.2)$);
\end{tikzpicture}
,
\end{gather}
where we use matrix notation for the twelve $H$-cells in $\jcell$. 
In this notation 
left cells are columns, right cells are rows, the $J$-cell is the whole block and 
$H$-cells are the small blocks.

We will also write $<_{l}$ or $\geq_{r}$ 
{\etc}, having the evident meanings.
Note that the three preorders also give rise to 
preorders on the set of cells, as well as between elements of $\monoid$ and cells.
For example, the notations $\lcell\geq_{l}a$ or $\lcell\leq_{l}\lcell^{\prime}$ 
make sense. In particular, for a fixed left cell $\lcell$ we can define
\begin{gather*}
\monoid_{\geq_{l}\lcell}=\{a\in\monoid|a\geq_{l}\lcell\},
\end{gather*}
as well as various versions which we will distinguish by the subscript.

\begin{Remark}\label{R:CellsWarning}
The cell orders need not to be total orders. In all of our 
examples the $\leq_{lr}$-order is a total order, but that is a coincidence.
\end{Remark}

\begin{Example}\label{E:CellsGroup}
If $\monoid$ is a group, then it has only one cell, the whole group, which 
is a left, right, $J$- and $H$-cell at the same time.
\end{Example}

\begin{Remark}\label{R:CellsNoGroups}
\autoref{E:CellsGroup} shows why the reader familiar with
the theory of groups might have never heard about cell theory:
for groups cell theory is trivial.
\end{Remark}

We write $\hcell(e)$ if $\hcell$ contains an idempotent $e\in S$. 
The $H$-cells of the form $\hcell(e)$ are called \emph{idempotent $H$-cells}, and the $J$-cells $\jcell(e)$ containing these 
$\hcell(e)\subset\jcell(e)$ are called \emph{idempotent $J$-cells}.

\begin{Remark}\label{R:CellsIdem}
In monoid theory idempotent $H$ and $J$-cells are called 
regular to avoid confusion with the property that 
{\eg} $\jcell\jcell=\jcell$ 
and because it is equivalent to each element
of the cell being von Neumann regular in the sense of the definition 
before \autoref{R:RepGapExtensionsRegular}.
However, for us the existence of an idempotent is crucial, so we use the above nomenclature.
\end{Remark}

$H$-cells are crucial as justified by:

\begin{Proposition}\label{P:CellsHCells}
For the monoid $\monoid$ we have:
\begin{enumerate}
\item Every $H$-cell is contained in some $J$-cell, and every $J$-cell is a disjoint union of $H$-cells.

\item $\hcell(e)$ is a group with identity $e$. In this case $\hcell(e)=\jcell(e)\cap(e\monoid e)$.
\end{enumerate}
\end{Proposition}

\begin{proof}
Part (a) is clear, while (b) is classical, see \cite[Theorem 7]{Gr-structure-semigroups}.
\end{proof}

\begin{Notation}\label{N:CellsHCells}
One case will play a special role, namely the case 
where $\hcell(e)$ is the trivial group. In this case 
we say $\hcell(e)$ is \emph{trivial} and write $\hcell(e)\cong\onemon$.
\end{Notation}

We have minimal and maximal $J$-cells in the $\leq_{lr}$-order.
In our illustrations the minimal cell will be at the bottom, so we call it the \emph{bottom cell} $\jb$, while the maximal cell will be at the top, so 
we call it the \emph{top cell} $\jt$.

\begin{Lemma}\label{L:CellsMaximal}
Every monoid has a unique bottom and top $J$-cell 
which are minimal respectively maximal in the $\leq_{lr}$-order. 
Both are idempotent $J$-cells.
\end{Lemma}

This is classical, {\eg} \cite[Chapter 6]{ClPr-algebraic-semigroups-2} discusses ordering relations on $J$-cells, but we will give a short proof for completeness.

\begin{proof}
The bottom $J$-cell is easy to find: Let $\group\subset\monoid$ 
be the group of units of $\monoid$, {\ie} 
the set of invertible elements of $\monoid$. Then $\group$ forms a left, a right and a $J$-cell at the same time, and is the smallest in all cell orders. 
To see this note that $1\leq_{l}a$ for all $a\in\monoid$ 
since we can choose $c=a$. 
But every invertible element $b\in\monoid$ satisfies $1=b^{-1}b$, which implies 
$b\leq_{l}1$, thus $b\sim_{l}1$. Similarly for $\sim_{r}$ and 
$\sim_{lr}$. The converse also holds, {\ie} every element in a minimal $J$-cell is invertible, so $\group$ is the 
unique bottom cell $\jb$. Moreover, the unit is an idempotent in $\jb$.

The top $J$-cell is not much harder to find:
If $\jcell$ and $\jcell^{\prime}$ are maximal $J$-cells, then 
$\jcell=\jcell\jcell^{\prime}=\jcell^{\prime}$ by maximality. 
Existence of a maximal $J$-cell follows from the finiteness of $\monoid$.
Furthermore, the $J$-cell $\jt$ contains an idempotent since $\jt\jt=\jt$ by maximality. 
This ensures the existence of an idempotent, see \cite[Proposition 1.23]{St-rep-monoid}.
\end{proof}

\begin{Example}\label{E:CellsTraMon}
The \emph{transformation monoid} $\tmon$ on the set $\{1,\dots,n\}$ 
is $\End(\{1,\dots,n\})$.
The cells of $\tmon[3]$, whose elements are written in one-line notation, with $(ijk)$ denoting the map $1\mapsto i,2\mapsto j,3\mapsto k$, are 
as follows. Using the illustration conventions as in \autoref{Eq:CellsIllustration} we have
\begin{gather*}
\xy
(0,0)*{\begin{gathered}
\begin{tabular}{C}
\arrayrulecolor{tomato}
\cellcolor{mydarkblue!25}(111) \\
\hline 
\cellcolor{mydarkblue!25}(222) \\
\hline 
\cellcolor{mydarkblue!25}(333)
\end{tabular}
\\[1pt]
\begin{tabular}{C|C|C}
\arrayrulecolor{tomato}
\cellcolor{mydarkblue!25}(122),(211) & \cellcolor{mydarkblue!25}(121),(212) & (221),(112)
\\
\hline
\cellcolor{mydarkblue!25}(133),(311) & (313),(131) & \cellcolor{mydarkblue!25}(113),(331)
\\
\hline
(233),(322) & \cellcolor{mydarkblue!25}(323),(232) & \cellcolor{mydarkblue!25}(223),(332)
\end{tabular}
\\[1pt]
\begin{tabular}{C}
\arrayrulecolor{tomato}
\cellcolor{mydarkblue!25}
\begin{gathered}
\phantom{.}
\\[-0.35cm]
(123),(213),(132)
\\[-1pt]
(231),(312),(321)
\end{gathered}
\end{tabular}
\end{gathered}};
(-46,13.8)*{\jt};
(-46,-2.5)*{\jm};
(-46,-16.5)*{\jb};
(46,13.8)*{\hcell(e)\cong\sym[1]};
(46,-2.5)*{\hcell(e)\cong\sym[2]};
(46,-16.5)*{\hcell(e)\cong\sym[3]};
\endxy
\quad
.
\end{gather*}
That is, $a\sim_{l}b$ if and only if $a(x)=a(y)\Leftrightarrow b(x)=b(y)$ (as functions), and $a\sim_{r}b$ if and only if they have the same image.
All idempotent $H$-cells are symmetric groups $\sym[k]$ of varying sizes. Note that not all $H$-cells contain idempotents: we have colored/shaded the $H$-cells containing idempotents.
\end{Example}

Let $|\lcell|$, $|\rcell|$ and $|\hcell|$ denote the 
sizes of fixed left, right and $H$-cells in a
$J$-cell $\jcell$ of size $|\jcell|$.

\begin{Lemma}\label{L:CellSizes}
Within one $J$-cell we have $|\lcell|=|\lcell^{\prime}|$, 
$|\rcell|=|\rcell^{\prime}|$ and $|\hcell|=|\hcell^{\prime}|$, 
and we have $|\lcell|\cdot|\rcell|/|\hcell|=|\jcell|$.
Moreover, $|\hcell|$ divides both, $|\lcell|$ and $|\rcell|$.
\end{Lemma}

\begin{proof}
The first three equalities follow from \cite[Theorem 1]{Gr-structure-semigroups}, 
the final two statements can then be shown from the previous three.
\end{proof}

Note that $|\lcell|,|\rcell|,|\jcell|,|\hcell|\in\N$, and \autoref{L:CellSizes} gives us additionally $|\lcell|/|\hcell|,|\rcell|/|\hcell|\in\N$.
These are important measures of the complexity of $\monoid$.

\begin{Example}\label{E:CellsTraMon2}
The middle $J$-cell in \autoref{E:CellsTraMon} 
has $|\hcell|=2$, $\jm=18=6\cdot 6/2=|\lcell|\cdot|\rcell|/|\hcell|$ and $|\lcell|/|\hcell|=|\rcell|/|\hcell|=3$.
\end{Example}

A \emph{left ideal} $I\subset S$ is a set such that $aI\subset I$. Right and two-sided ideals are defined similarly.
The following lemma explains the matrix notation:

\begin{Lemma}\label{L:CellsIdeal}
For fixed left cell $\lcell$ the set $\monoid_{\geq_{l}\lcell}$ is a left ideal in $\monoid$.
Similarly, $\monoid_{\geq_{r}\rcell}$ is a right and 
$\monoid_{\geq_{lr}\jcell}$ is a two-sided ideal. The same works 
when replacing $\geq$ by $>$.
\end{Lemma}

\begin{proof}
Directly from the definitions: 
given $b\in\monoid_{\geq_{l}\lcell}$, the element $ab$ 
is still left greater or equal to $l\in\lcell$ since $b=cl$ for some $c$.
\end{proof}

Let us state how cell theory helps to understand periods of 
elements, which in turn are of importance 
in cryptography. To this end, 
recall that the \emph{index} $i(a)\in\N$ for $a\in\monoid$ is the smallest number such that
$a^{i(a)}=a^{i(a)+d}$ for some $d\in\Z_{>0}$. The smallest possible $d$ 
is then in turn called the \emph{period} of $a$ and we denote it by $p(a)$.

\begin{Theorem}\label{T:CellsPeriod}
There exists an $H$-cell $\hcell(e)$ such that $\cyclic[{p(a)}]\cong\{a^{s}\mid s\geq i(a)\}\subset\hcell(e)$ is a subgroup. In particular, 
$p(a)\mid|\hcell(e)|$.
\end{Theorem}

\begin{proof} 
As a consequence of \cite[Theorem 7]{Gr-structure-semigroups}, 
the $H$-cells of the form $\hcell(e)$ 
are the maximal subgroups of $\monoid$, so no other 
subgroup will be contained in some $\hcell(e)$.
\end{proof}

\begin{Example}\label{E:CellsGenOneElement}
Given $i\in\N$, $p\in\Z_{\geq 1}$ form the finite cyclic monoid 
$\cyclic[i,p]=\langle a|a^{i+p}=a^{i}\rangle$ of cardinality $i+p$. 
The element $a$ has index $i(a)=i$ and period $p(a)=p$.
Moreover, the monoid $\cyclic[i,p]$ 
is commutative, so left, right and $J$-cells coincide. The elements $1,a,\dots,a^{i-1}$ each constitute a single $J$-cell, 
in total $i-1$ such $J$-cells.
All the remaining elements $\jt=\{a^i,a^{i+1},\dots,a^{i+p-1}\}$ constitute one $J$-cell (the top cell) which is a cyclic group of order $p$ under multiplication. The element $e=a^{pj}$ where $j$ is such that $i\leq pj<i+p$ is the idempotent for $\jt=\hcell(e)$ and the identity of that group. Out of the $i+1$ cells in $\cyclic[i,p]$ two cells are idempotent: $\jb=\{1\}$ 
and $\jt$.

To be completely explicit, let us consider $\cyclic[3,2]$, 
which is the monoid generated by one element $a$ 
of index $3$ and period $2$. Then $\cyclic[3,2]=\{1,a,a^{2},a^{3},a^{4}\}$	
and its cell structure is
\begin{gather*}
\xy
(0,0)*{\begin{gathered}
\begin{tabular}{C}
\arrayrulecolor{tomato}
\cellcolor{mydarkblue!25}\raisebox{-0.07cm}{$a^{3},a^{4}$}
\end{tabular}
\\[1pt]
\begin{tabular}{C}
\arrayrulecolor{tomato}
a^{2}
\end{tabular}
\\[1pt]
\begin{tabular}{C}
\arrayrulecolor{tomato}
a
\end{tabular}
\\[1pt]
\begin{tabular}{C}
\arrayrulecolor{tomato}
\cellcolor{mydarkblue!25}1
\end{tabular}
\end{gathered}};
(-20,10)*{\jt};
(-20,3.5)*{\jcell_{a^{2}}};
(-20,-3.5)*{\jcell_{a}};
(-20,-10)*{\jb};
(23,10)*{\hcell(e)\cong\cyclic[2]\cong\Z/2\Z};
(20,-10)*{\hcell(e)\cong\onemon};
\endxy
\quad
.
\end{gather*}
Note that $\monoid_{a}$ is commutative, so left, right, $J$- and $H$-cells 
agree.
\end{Example}

\begin{Remark}[\textbf{\emph{Additional Task 2}}]\label{R:CellsCrypto}
Using the DH protocol with protocol 
monoid $\monoid$ other than a group, it would be important 
to find elements $g\in\monoid$ of big period that has a large prime factor, see {\eg} the original DH key exchange \cite{Ko-algebraic-cryptography},  \cite[Section 1.2]{MyShUs-group-cryptography}.
So, by \autoref{T:CellsPeriod}, it would be preferable to have a monoid $\monoid$ 
with $H$-cells whose orders have large prime divisors since the period of $a\in\monoid$ divides the order 
of the idempotent $H$-cell of $\monoid$ that contains the top cell of $\cyclic[i,p]$.
\end{Remark}


\subsection{Classification of simple representations}\label{SS:CellsClassification}

Recall that we consider $\monoid$-representations defined over $\K$.

Cells can be considered $\monoid$-representations, called \emph{cell representations} or \emph{Sch{\"u}tzenberger representations}, up to higher order terms:

\begin{Lemma}\label{L:CellsMod}
Each left cell $\lcell$ of $\monoid$ gives rise to a left $\monoid$-representation $\lmod=\K\lcell$ by
\begin{gather*}
a\acts l\in\lmod=
\begin{cases}
al&\text{if }al\in\lcell,
\\
0&\text{else.}
\end{cases}
\end{gather*}
Similarly, right cells give right $\monoid$-representations $\rmod$ and $J$-cells give $\monoid$-birepresentations (often called $\monoid$-birepresentations). We have $\dimk(\lmod)=|\lcell|$ and $\dimk(\rmod)=|\rcell|$.
\end{Lemma}

\begin{proof}
Directly from the definitions.
\end{proof}

The annihilator $\mathrm{Ann}_{\monoid}(\module)=\{s\in\monoid|s\acts\module=0\}$ 
of an $\monoid$-representation $\module$ is a two-sided ideal of $\monoid$. 
An \emph{apex} of $\module$ is a $J$-cell $\jcell$ 
such that, firstly, $\jcell\cap\mathrm{Ann}_{\monoid}(\module)=\emptyset$, and 
secondly, all $J$-cells $\jcell^{\prime}$ 
with $\jcell^{\prime}\cap\mathrm{Ann}_{\monoid}(\module)=\emptyset$ satisfy 
$\jcell^{\prime}\leq_{lr}\jcell$.  
In other words, an apex is the $\leq_{lr}$-maximal $J$-cell not annihilating $\module$.
The following justifies the terminology of the \emph{apex of a simple $\monoid$-representation}:

\begin{Lemma}\label{L:CellsApex}
Every simple $\monoid$-representation has a unique apex.
\end{Lemma}

\begin{proof}
This is classical, see {\eg} \cite[Theorem 5]{GaMaSt-irreps-semigroups}.
\end{proof}

\begin{Example}\label{E:CellsSimples}
The apex of $\oneb$ is always $\jb=\group$.
On the other hand, the apex of $\onet$ 
is $\jt$ since every $s\in\monoid$ acts as $1$.
\end{Example}

Recall that the nonunital way to induce is 
$\ind(\module)=\K\monoid e\otimes_{\K e\monoid e}\module$
for some idempotent $e\in\monoid$, see {\eg} \cite[Section 4.1]{St-rep-monoid} (inducing from the submonoid $e\monoid e$ to $\monoid$, or rather using their 
monoid algebras).
It follows from \cite{Gr-structure-semigroups} that $\lmod$ 
is a free right $\hcell(e)$-representation, and this action commutes 
with the left $\monoid$-action. Thus, $\lmod$ is a $\monoid$-$\hcell(e)$-birepresentation. We can then define
an induction functor
\begin{gather*}
\ind_{\hcell(e)}^{\monoid}\module
=
\lmod\otimes_{\hcell(e)}\module
,
\end{gather*}
where $\module$ is a left $\module$-representation.

\begin{Example}\label{E:CellsRegular}
Let $\K[\hcell(e)]$ denote the regular 
$\hcell(e)$-representation, which as a $\K$-vector space 
is just $\K\hcell(e)$ and the $\hcell(e)$-action is the multiplication action.
We have $\ind_{\hcell(e)}^{\monoid}\K[\hcell(e)]\cong\lmod$ 
as left $\monoid$-representations.
\end{Example}

Recall also that the \emph{head} 
$\hd(\module)$ of an $\monoid$-representation $\module$ is
the maximal semisimple quotient of $\module$. It is well-defined, up to isomorphism, for any representation over a finite monoid and is isomorphic to the 
quotient $\module/\rad(\module)$. Here $\rad(\module)$ 
denotes the radical, which is the intersection of all maximal subrepresentations of $\module$.

We get the \emph{Clifford--Munn--Ponizovski\u{\i} theorem} 
or \emph{$H$-reduction}:

\begin{Proposition}\label{P:CellsSimples}
For a monoid $\monoid$:
\begin{gather*}
\{\text{simple $\monoid$-representations of apex $\jcell$}\}/\cong\;
\xleftrightarrow{1{:}1}
\{\text{simple $\hcell(e)$-representations}\}/\cong\;
,
\end{gather*}
where $\hcell(e)\subset\jcell$ is any 
arbitrarily chosen idempotent $H$-cell in 
an idempotent $J$-cell $\jcell$. Moreover, an explicit bijection (from right to left) is given by
\begin{gather*}
K\mapsto\simple[K]\cong\hd(\ind_{\hcell(e)}^{\monoid}K).
\end{gather*}
\end{Proposition}

\begin{proof}
The above is an easy reformulation of \cite[Theorem 7]{GaMaSt-irreps-semigroups} 
or \cite[Theorem 5.5]{St-rep-monoid}.
\end{proof}

Note that only idempotent $J$-cells contribute to the classification.
We usually omit to write {\eg} ``simples up to isomorphism'' in the rest of the paper.

\begin{Remark}\label{R:CellsSimples}
The $1{:}1$ correspondence 
in \autoref{P:CellsSimples} always exists regardless of $\K$. 
However, the classification still depends on $\K$ since 
the number of simple $\hcell(e)$-representation does.
\end{Remark}

\begin{Example}\label{E:CellsTraMonSimples}
Let $\cchar$ be such that $\cchar\nmid 3!=6$, {\eg} $\cchar=0$. The cell structure from \autoref{E:CellsTraMon} 
shows that $\tmon[3]$ has three simple $\tmon[3]$-representations 
of apex $\jb$, two of apex $\jm$ 
and one of apex $\jt$ since the associated $\hcell(e)$ 
are the symmetric groups $\sym[3]$, 
$\sym[2]$ and $\sym[1]$ (and the number of simple $\sym$-representations 
is given by the number of partitions of $n$). 

For $\cchar=3$ one gets only two simple $\tmon[3]$-representations of apex $\jb$ since $\sym[3]$ has only two simple representations in this 
characteristic; the rest remains the same as for $\cchar=0$. Similarly, for $\cchar=2$ both apexes $\jb$ and $\jm$ 
have one fewer associated simple $\tmon[3]$-representation than for $\cchar=0$, but $\jt$ still has the same count.
\end{Example}

We can thus define a partial order, also denoted by $\leq_{lr}$, on the set of simple $\monoid$-representations by saying that one simple is strictly smaller than another if its apex is strictly smaller. Note that simples of the same apex are incomparable.

\begin{Example}\label{E:CellsHTrivial}
Note that if $\hcell(e)$ is trivial, then \autoref{P:CellsSimples} implies 
that one can say that the simples 
are indexed by the poset of apexes. 
\end{Example}

\begin{Remark}\label{R:CellsHTrivial}
When working over $\C$ and when all $J$-cells are idempotent, it is shown in 
\cite[Theorem 2.1]{Pu-semigroups-hw-categories} that $\leq_{lr}$ makes the 
representation category of $\monoid$ 
into a highest weight category in the sense of \cite{ClPaSc-h-weight-qh}.
In fact, for the reader familiar with cellular algebras as in \cite{GrLe-cellular}, \cite{TuVa-handlebody} or \cite{Tu-sandwich} we point out that \cite[Theorem 2.1]{Pu-semigroups-hw-categories} shows that, if all $J$-cells are idempotent, then the monoid algebra $\C\monoid$ is a quasi-hereditary sandwich cellular algebra.

As a historical remark, 
the fact that the monoid algebra $\C\monoid$ of a regular 
monoid (a regular monoid satisfies any of the conditions in 
\autoref{L:CellsAdmissible} below) in characteristic zero
is a quasi-hereditary sandwich cellular algebra was first proven in 
\cite{Ni-semigroups-hw-categories} in the early 1970s. Of course 
the result was phrased in a different language since \cite{Ni-semigroups-hw-categories} appeared before quasi-hereditary or (sandwich) cellular algebras were defined.
\end{Remark}


\subsection{Cells and (semisimple) representation gaps}\label{SS:CellsDimSimples}

Note that \autoref{P:CellsSimples} makes it easy to classify simple $\monoid$-representations but does not give much information about their dimensions.

\begin{Theorem}\label{T:CellsDimensions}
The dimension of the simple $\monoid$-representation $\simple[K]$ 
associated to the simple $\hcell(e)$-representation
$K$ via \autoref{P:CellsSimples} can be bounded by
\begin{gather*}
\dimk(\simple[K])\leq
|\lcell|/|\hcell|\cdot\dimk(K).
\end{gather*}
\end{Theorem}

\begin{proof}
First, recall from \autoref{L:CellSizes} that all left and $H$-cells within one $J$-cell are of the same size, so for the bound we can and will omit writing $\lcell(e)$ and $\hcell(e)$.
Then this follows from the explicit bijection in \autoref{P:CellsSimples} and the fact that $\lmod$ is a free $\hcell(e)$-representation of rank $|\lcell|/|\hcell|$.
\end{proof}

Note that dimension of $\hd(\ind_{\hcell(e)}^{\monoid}K)$ depends on the field, in general, and can be hard to compute. 
The quantity $|\lcell|/|\hcell|\cdot\dimk(K)$ is often easy to compute in practice so we define:

\begin{Definition}\label{E:CellsRepGap}
We call $\ssdimk(\simple[K])=|\lcell|/|\hcell|\cdot\dimk(K)$ 
the \emph{semisimple dimension} of $\simple[K]$.
The minimal $m$ such that there is a nontrivial simple $\monoid$-representation
with $\ssdimk(\simple[K])=m$ is called the 
\emph{semisimple representation gap} $\ssgap{\monoid}$ of $\monoid$.

We also call $\ssratio{\monoid}=\tfrac{\ssgap{\monoid}}{\sqrt{|\monoid|}}$
the \emph{semisimple-gap-ratio}.
\end{Definition}

The square root in the definition of 
$\ssratio{\monoid}$	is used for the same reasons as in 
\autoref{R:RepGapCrypto2}.
With the same assumptions as in {\eg} \autoref{L:RepGapSimple} we have:

\begin{Theorem}\label{T:CellsBound}
Assume that $\monoid$ has at least one nontrivial simple representations. We have
\begin{gather*}
\gap{\monoid}\leq\min\{\dimk(\simple)|\text{$\simple\not\cong\onebt$ is a simple $\monoid$-representation}\}\leq\ssgap{\monoid}\leq|\monoid|.
\end{gather*}
\end{Theorem} 

\begin{proof}
Clear by definition and \autoref{L:RepGapSimple}.
\end{proof}

\begin{Remark}[\textbf{\emph{Additional Task 3}}]\label{R:CellCryptoMain}
As before, it is important 
for potential cryptographic applications to find monoids 
with $\ssgap{\monoid}$ exponentially big.
\end{Remark}

\begin{Example}\label{E:CellsTraDim}
In the setting of \autoref{E:CellsTraMon}
and \autoref{E:CellsTraMonSimples} 
(in particular, $\cchar\nmid 6$) we have the following.

The three simple $\tmon[3]$-representations of apex $\jb$ are the simple $\sym[3]$-representations inflated to $\tmon[3]$, so they 
are of dimensions $1$, $2$ and $1$ 
(one of these is $\oneb$). These are also their 
semisimple dimensions. 

The simple $\sym[3]$-representation 
of apex $\jt$ can be identified with $\onet$, 
so is of dimension one, which is also its semisimple 
dimension. 

The two simple $\sym[3]$-representations 
of apex $\jm$ are induced from the 
respective $\sym[2]$-representations, and are of semisimple dimension $3$. 
One can check that they are of dimensions $3$ respectively $2$.

In general, for the representation theory of $\tmon[n]$ see \cite[Section 4]{Pu-semigroups-hw-categories} or \cite[Section 5.3]{St-rep-monoid}.
\end{Example}

The name semisimple representation gap is justified by the following.

\begin{Proposition}\label{P:CellsSemisimple}
The following are equivalent.
\begin{enumerate}

\item The monoid $\monoid$ is semisimple over $\K$.

\item All $J$-cells are idempotent, all $\hcell(e)$ are semisimple over $\K$ and
$\dimk(\simple[K])=\ssdimk(\simple[K])$ for all simple 
$\monoid$-representations $\simple[K]$.

\end{enumerate}
\end{Proposition}

\begin{proof}
This follows from \cite[Theorem 5.19]{St-rep-monoid} 
and the paragraph below that theorem.
\end{proof}


\subsection{Cells and Gram matrices}\label{SS:CellsGram}

Recall the following construction of Gram matrices, also called sandwich matrices in monoid theory, see {\eg} 
\cite[Section 5.2]{ClPr-algebraic-semigroups} or \cite[Section 5.4]{St-rep-monoid}.
Fix an idempotent $H$-cell $\hcell(e)=\lcell\cap\rcell$ 
in some idempotent $J$-cell $\jcell$. Then 
$\lcell$ is a free right $\hcell(e)$-set and $\rcell$ is a
free left $\hcell(e)$-set, so we can let $\{l_{1},\dots,l_{R}\}$ 
and $\{r_{1},\dots,r_{L}\}$ complete sets of representatives for
$\lcell/\hcell(e)$ respectively for $\hcell(e)\text{\textbackslash}\rcell$.
Here $R$ is the number of right cells and $L$ is the number of left cells in $\jcell$.

The \emph{Gram matrix} $P^{\jcell}=(P^{\jcell}_{i,j})_{i,j}$ is the matrix with values in $\K\hcell(e)$ defined by
\begin{gather*}
P^{\jcell}_{i,j}=
\begin{cases}
r_{i}l_{j}&\text{if }r_{i}l_{j}\in\hcell(e),
\\
0&\text{else}.
\end{cases}
\end{gather*}
Note that $P^{\jcell}$ depends on choices, but one can show that 
its important properties do not depend on these choices, see the references above.

Gram matrices are in particularly useful for $\hcell(e)\cong\onemon$ and $L=R$ 
as justified by part (a) of the following (which the reader familiar with \cite{GrLe-cellular} might recognize):

\begin{Proposition}\label{P:CellsGram}
Fix an idempotent $J$-cell $\jcell$. All cells in the statement 
are within $\jcell$.

\begin{enumerate}

\item Assume
$\hcell(e)\subset\jcell$ satisfies $\hcell(e)\cong\onemon$. 
Assume further that $P^{\jcell}$ is square and symmetric.
Let $\simple[\jcell]$ 
denote the associated simple $\monoid$-representation, see \autoref{P:CellsSimples}. Then:
\begin{gather*}
\dimk(\simple[\jcell])=\rank(P^{\jcell}).
\end{gather*}

\item More generally, let $K$ is a simple $\hcell(e)$-representation  
and let $\simple[K]$ be the associated simple $\monoid$-representation.
Let $P^{\jcell}_{K}$ denote the
matrix one gets by applying $K$ to each entry of $P^{\jcell}$.
Then:
\begin{gather*}
\dimk(\simple[K])=\rank(P^{\jcell}_{K}).
\end{gather*}

\end{enumerate}
\end{Proposition}

\begin{proof}
\textit{(a).} Let $\radd$ denote the radical of the symmetric bilinear form associated to 
$P^{\jcell}$. We claim that $\radd$ is a $\monoid$-submodule of 
the corresponding cell representation $\lmod$. 
To see this note that $r_{i}l_{j}\notin\hcell(e)$ can only occur 
if they end up in $\jcell_{>_{lr}\jcell}$, and multiplying by elements 
from $\monoid$ preserves this property.

We further claim that any element in $\lmod\setminus\radd$ generates $\lmod$.
This can be proven as in \cite[Lemma 3.4]{EhTu-relcell}.

It follows that $\lmod/\radd$ is a simple $\monoid$-representation since 
any proper submodule of it must be contained in $\radd$. Since the apex 
of $\lmod/\radd$ is $\jcell$, by construction, it follows that
$\lmod/\radd\cong\simple[\jcell]$. The proof completes.

\textit{(b).} Adjusting the arguments in (a), see {\eg} 
\cite[Corollary 5.30]{St-rep-monoid} for details.
\end{proof}

\begin{Theorem}\label{T:CellsGramSub}
Let $\monoid[R]\subset\monoid$ be a submonoid. Under the assumptions in 
\autoref{P:CellsGram}, if $\jcell$ restricts to
an idempotent $J$-cell of $\monoid[R]$, then
\begin{gather*}
\dimk(\simple[\jcell]^{\monoid})\geq\dimk(\simple[\jcell]^{\monoid[R]}),
\end{gather*}
for the associated simple $\monoid[R]$ and $\monoid$-representations.
\end{Theorem}

\begin{proof}
Note that under the assumptions we have that the Gram matrix for $\monoid[R]$ is 
a submatrix of $P^{\jcell}$.
The rank of a matrix is always greater or equal to the 
rank of a submatrix, so the statement follows by \autoref{P:CellsGram}.
\end{proof}

We stress that it is not generally true that $\jcell$ restricts to
a(n idempotent) $J$-cell of $\monoid[R]$, so the assumption in 
\autoref{T:CellsGramSub} is necessary.


\subsection{Cells, Burnside--Brauer--Steinberg and faithfulness}\label{SS:CellsFaithful}

Let $cl_{\hcell(e)}$ denote the number of conjugacy 
classes of the group $\hcell(e)$.
Let $\{e_{1},\dots,e_{r}\}$ be a choice of one idempotent per 
idempotent $J$-cell, and define
\begin{gather*}
cl(\monoid)=cl_{\hcell(e_{1})}+\dots+cl_{\hcell(e_{r})}.
\end{gather*}

\begin{Lemma}\label{L:CellsFaithful}
The number $cl(\monoid)\in\N$ is independent of the choice of $\{e_{1},\dots,e_{r}\}$.
\end{Lemma}

\begin{proof}
This is a consequence of \cite[Section 7.1]{St-rep-monoid}.
\end{proof}

Hence, $cl(\monoid)$ is a constant depending on $\monoid$ only.
One can use $cl(\monoid)$ for the \emph{Burnside--Brauer theorem} 
(characteristic zero)
and the \emph{Steinberg theorem} (arbitrary characteristic):

\begin{Proposition}\label{P:CellsBurnsideBrauer}
If $\fmodule$ is a faithful $\monoid$-representation,
then every simple $\monoid$-representation 
appears as a composition factor of $\fmodule^{\otimes k}$ 
for some $0\leq k\leq cl(\monoid)-1$. Moreover, if $\fmodule$ is a faithful $\monoid$-representation whose	composition factors are one-dimensional, then the composition factors of $\fmodule^{\otimes k}$ are also one-dimensional.
\end{Proposition}

\begin{proof}
For characteristic zero see \cite{St-burnside-brauer} or 
\cite[Section 7.4]{St-rep-monoid} and the observation 
that the $r$ in that theorem satisfies $r\leq cl(\monoid)$ 
by the discussion in \cite[Section 7.1]{St-rep-monoid}. 
For the characteristic free version see 
\cite[Corollary 10.7]{St-rep-monoid}, using the same observation.
\end{proof}

\begin{Example}\label{E:CellsBurnsideBrauer}
The bound given in \autoref{P:CellsBurnsideBrauer} 
is often not optimal but cannot be improved uniformly. For example, 
for $\cyclic$ we have $cl(\cyclic)=n$. 
Assume $n$ is prime. 
Over $\C$ the $n$th primitive root of unity $\exp(\tfrac{2\pi i}{n})$
gives rise to a $1$-dimensional faithful $\cyclic$-representation, 
and only the $(n-1)$th power of it will contain 
the simple $\cyclic$-representation associated to $\exp(\tfrac{2\pi i(n-1)}{n})$.
\end{Example}

The Burnside--Brauer--Steinberg theorem 
\autoref{P:CellsBurnsideBrauer} gives a bound for the dimension of faithful 
$\monoid$-representations:

\begin{Theorem}\label{T:CellsFaithful}
Let $\cchar=0$, and let $\simple[max]$ be a simple $\monoid$-representation 
of the biggest dimension.
If $\fmodule$ is a faithful $\monoid$-representation, then
$\dimk(\fmodule)\geq\sqrt[{cl(\monoid)-1}]{\dimk(\simple[max])}$.
Hence,
\begin{gather*}
\faith{\monoid}
\geq
\sqrt[{cl(\monoid)-1}]{\dimk(\simple[max])}.
\end{gather*}
\end{Theorem}

\begin{proof}
This follows from \autoref{P:CellsBurnsideBrauer}.
\end{proof}

Note that one can use \autoref{T:CellsFaithful} often in combination with 
\autoref{L:RepGapFaithfulEmbedding}.

\begin{Remark}[\textbf{\emph{Additional Task 4}}]\label{R:CellCryptoFaithful}
Thus, by \autoref{T:CellsFaithful} 
it is preferable for cryptographical applications
to find a monoid $\monoid$ with $c(\monoid)$ being small.
\end{Remark}

\begin{Example}\label{E:CellsTraMonFaith}
Applying \autoref{T:CellsFaithful} for 
$\tmon[3]$ gives $\sqrt[6]{3}$ as a lower bound, which rounds to $2$.
The smallest faithful $\tmon[3]$-representation is $\K\{1,2,3\}$ 
(with the defining action), 
so of dimension three.
\end{Example}

With respect to extensions as discussed in \autoref{SS:RepGapRepGap}
we get:

\begin{Proposition}\label{P:CellsBurnsideBrauerExtensions}
There is 
a faithful $\monoid$-representation  containing only 
$\onebt$ as composition factors if and only if $\monoid$ 
has at most two idempotent $J$-cells and all idempotent $H$-cells 
are trivial, {\ie} $\hcell(e)\cong\onemon$.
\end{Proposition}

\begin{proof}
\textit{$\Rightarrow$.}	If $\fmodule$ is a faithful $\monoid$-representations only containing 
$\onebt$ as composition factors, then \autoref{P:CellsBurnsideBrauer} implies 
that there can be no simple $\monoid$-representations except $\onebt$.
Thus, the result follows by \autoref{P:CellsSimples}.
\medskip

\textit{$\Leftarrow$.} In this case \autoref{P:CellsSimples}
implies that $\onebt$ are the only simple $\monoid$-representations. 
\end{proof}

\begin{Example}\label{E:CellsBurnsideBrauerExtensions}
Let $\cchar=0$.
When $\monoid$ is a group \autoref{P:CellsBurnsideBrauerExtensions} implies that only the trivial 
group has faithful representations entirely made of trivial representations. 
(Note that this is clear because of a different 
reason: the assumption is $\cchar=0$ so the representation theory 
of groups is semisimple.)
\end{Example}

\begin{Example}\label{E:CellsBurnsideBrauerExtensions2}
It follows from the discussion in \autoref{Eq:TLCells} 
that the Temperley--Lieb monoid on three strands $\tlmon[3]$ 
is an example of a nontrivial monoid that has a faithful representation entirely made of $\onebt$. This works in arbitrary characteristic.
\end{Example}


\subsection{Cell submonoids and subquotients}\label{SS:CellsSubmon}

Recall that simple $\monoid$-representations arrange themselves 
according to the cells, see \autoref{P:CellsSimples}. Let us 
in this motivational paragraph for simplicity assume that
$\hcell(e)\cong\onemon$ for all idempotent $H$-cells and that all $J$-cells 
are idempotent. Then the dimensions of the simple $\monoid$-representations 
very often have the following form, which is roughly as expected from combinatorial numbers:
\begin{gather}\label{Eq:CellsDimTL}
\begin{gathered}
\scalebox{0.9}{$
\begin{tikzpicture}[anchorbase]
\node at (0,0) {\includegraphics[height=4.4cm]{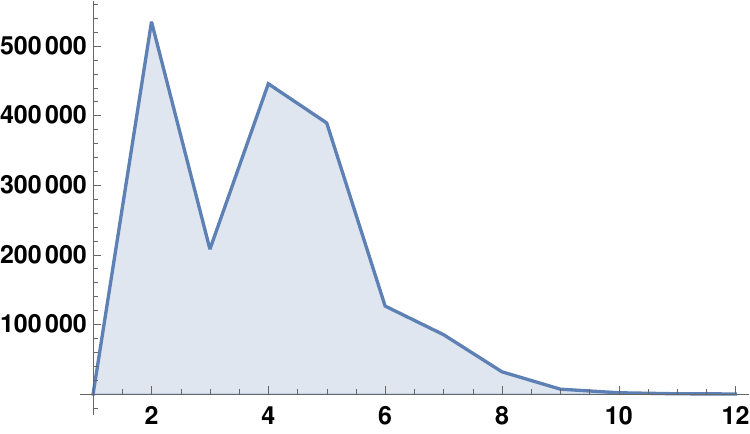}};
\node at (2.25,2) {$\tlmon[24]/\Q$};
\node at (2.25,1.5) {cells increase $\leftarrow$};
\draw[->] (-2.4,2.2) node[above,xshift=0.3cm]{y-axis: dim} to (-2.4,2.1) to (-2.8,2.1);
\draw[->] (2.75,-1.5) node[above,yshift=0.75cm]{x-axis:}node[above,yshift=0.3cm]{\# through}node[above,yshift=-0.15cm]{strands/2} to (2.75,-1.8);
\end{tikzpicture}
,
\begin{tikzpicture}[anchorbase]
\node at (-0.1,0) {\includegraphics[height=4.4cm]{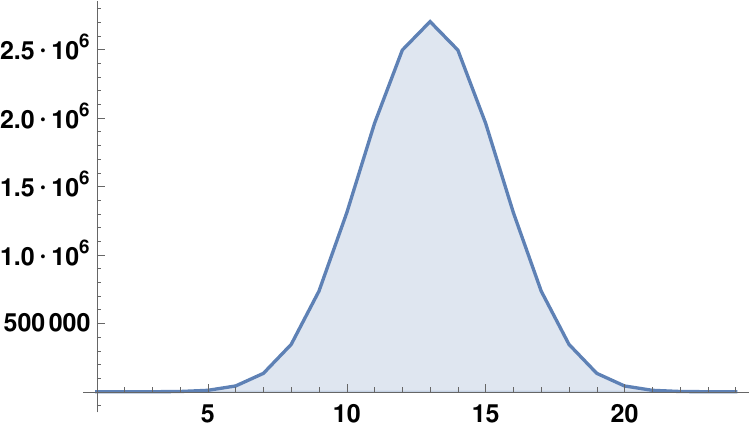}};
\node at (2.25,2) {$\promon[24]/\K$};
\node at (2.25,1.5) {cells increase $\leftarrow$};
\draw[->] (-2.4,2.2) node[above,xshift=0.3cm]{y-axis: dim} to (-2.4,2.1) to (-2.8,2.1);
\draw[->] (2.75,-1.5) node[above,yshift=0.75cm]{x-axis:}node[above,yshift=0.3cm]{\# through}node[above,yshift=0.0cm]{strands} to (2.75,-1.8);
\end{tikzpicture}
$}
,
\\
\scalebox{0.9}{$\text{TL dim}\colon\left(1,534888,208011,445741,389367,126292,85216,31878,6876,1726,252,22,1\right)$}
,
\\
\scalebox{0.7}{$\text{pRo dim}\colon\left(
\begin{gathered}
1,24,276,2024,10626,42504,134596,346104,735471,1307504,1961256,2496144,
\\
2704156,2496144,1961256,1307504,735471,346104,134596,42504,10626,2024,276,24,1
\end{gathered}
\right)$}
.
\end{gathered}
\end{gather}
These illustrations show the dimensions of the simple $\tlmon[24]$-representation (left) over $\Q$ (or any field of characteristic zero) and the simple $\promon[24]$-representations for general $\K$, respectively. See \autoref{S:TL} below 
for details. (Note the two trivial 
$\tlmon[24]$- respectively $\promon[24]$-representations 
of dimension one 
for the bottom and top cell.) Thus, it seems preferable to
cut-off the representations for small cells, and get rid 
of the fluctuations for very big cells.

The key to do the first is
are \emph{cell submonoids} as follows. 

\begin{Definition}\label{D:CellsSubmonoid}
For a $J$-cell 
$\jcell$ with $1\notin\jcell$ define the \emph{$\jcell$-submonoid}
\begin{gather*}
\monoid_{\geq\jcell}=
\monoid_{\geq_{lr}\jcell}\cup\{1\}.
\end{gather*}
\end{Definition}

In words, we artificially adjoint a unit $1$ 
(strictly speaking we should write $1^{\prime}$) to 
the two-sided ideal $\monoid_{\geq_{lr}\jcell}$ from \autoref{L:CellsIdeal}.

\begin{Lemma}\label{L:CellsSubmonoid}
For any $J$-cell 
$\jcell$ with $1\notin\jcell$, 
$\monoid_{\geq\jcell}$ is a submonoid of $\monoid$.
\end{Lemma}

\begin{proof}
By \autoref{L:CellsIdeal}.
\end{proof}

\begin{Remark}\label{R:CellsSemigroup}
There are minor, but not essential, differences between representations 
of monoids and semigroups. Adjoining a unit is for convenience 
only so that we do not need to leave the 
world of monoids.
\end{Remark}

Annihilating the bigger cells can be done using
the \emph{Rees factor} $\monoid/I$ of a monoid $\monoid$ by a 
two-sided ideal $I$.
The construction works as follows. As a set $\monoid/I=(\monoid\setminus I)\cup\{0\}$, where one artificially adjoints an element $0$. The 
multiplication is $s\bullet t=st$ if $s,t,st\in\monoid\setminus I$, and 
$s\bullet t=0$ otherwise.

\begin{Lemma}\label{L:CellsQuotient}
For any two-sided ideal, the Rees factor $\monoid/I$ is a well-defined monoid.
\end{Lemma}

\begin{proof}
An easy exercise, see also \cite[Exercise 1.6]{St-rep-monoid}.
\end{proof}

We can thus define \emph{cell subquotients}:

\begin{Definition}\label{D:CellsSubquotient}
For two $J$-cells $\jcell\leq_{lr}\kcell$ with $1\notin\jcell$ define 
the \emph{$\jcell$-$\kcell$-subquotient} as the Rees factor
\begin{gather*}
\monoid_{\jcell}^{\kcell}
=
\monoid_{\geq\jcell}/\monoid_{\geq\kcell}.
\end{gather*}
Here we additionally allow the following
extremal cases:
\begin{gather*}
\monoid_{\jcell}^{none}
=
\monoid_{\geq\jcell}
,\quad
\monoid_{none}^{\kcell}
=
\monoid/\monoid_{\geq\kcell}
,\quad
\monoid_{none}^{none}
=
\monoid.
\end{gather*}
We also call all of the above \emph{cell subquotients} for short.
\end{Definition}

By \autoref{L:CellsSubmonoid}
and \autoref{L:CellsQuotient}, $\monoid_{\jcell}^{\kcell}$ is 
a subquotient of $\monoid$.
Unless we are in one of the extreme cases,
$\monoid_{\jcell}^{\kcell}$ has $\jb=\{1\}$ and $\jt=\{0\}$.
Both are left, right, $J$- and $H$-cells at the same time.

\begin{Lemma}\label{L:CellsAdmissible}
The following conditions are equivalent:
\begin{enumerate}

\item For all left cells $\lcell$:
\begin{gather*}
\forall a,b\in\lcell\;\exists c\in\jcell\supset\lcell\text{ such that }a=cb.
\end{gather*}

\item For all right cells $\rcell$:
\begin{gather*}
\forall a,b\in\rcell\;\exists c\in\jcell\supset\rcell\text{ such that }a=bc.
\end{gather*}

\item For all $J$-cells $\jcell$:
\begin{gather*}
\forall a,b\in\jcell\;\exists c,d\in\jcell\text{ such that }a=cbd.
\end{gather*}

\item All $J$-cells are idempotent.

\item For all $a\in\monoid$ we have $a\in a\monoid a$.

\end{enumerate}
\end{Lemma}

\begin{proof}
Well-known, see {\eg} \cite[Theorem A.3.7]{RhSt-qtheory}.
\end{proof}

We say 
$\monoid$ is \emph{regular} 
(this is also sometimes called von Neumann regular) if any of the 
equivalent conditions in 
\autoref{L:CellsAdmissible} hold.

The regularity condition ensures that 
the cells are not affected when taking cell subquotients.

\begin{Lemma}\label{R:CellsJCellSubmonoid}
Let $\monoid$ be regular.
In the nonextremal cases 
the $J$-cells of $\monoid_{\jcell}^{\kcell}$ are given by
\begin{gather*}
\{\jb\}\cup\{\mcell|\text{$\mcell$ is a $J$-cell of $\monoid$ with }\jcell\leq_{lr}\mcell<_{lr}\kcell\}\cup\{\jt\}.
\end{gather*}
Similarly for left, and right cells, assuming the respective 
regularity condition, and $H$-cells.

An analog statement holds in the extremal cases.
\end{Lemma}

\begin{proof}
By the regularity assumption, the remaining elements 
of $\monoid_{\jcell}^{\kcell}$ arrange themselves into cells 
precisely as in $\monoid$.
\end{proof}

We require that $\monoid$ is regular for the remainder of this section.

Assume that we are in the nonextremal cases.
Then $\monoid_{\jcell}^{\kcell}$ has trivial representations $\oneb$ 
and $\onet$ associated to the apexes $\jb$ and $\jt$, 
and these 
are the only $\monoid_{\jcell}^{\kcell}$-representations of these apexes.
The other simple $\monoid_{\jcell}^{\kcell}$-representations 
and their dimensions are 
given by the following statement. Note hereby that any $\monoid_{\jcell}$-representation 
with apex $\mcell$ can be \emph{inflated} to a $\monoid$-representation by letting all elements in $\monoid_{<_{lr}\jcell}$ 
act by zero.

\begin{Proposition}\label{P:CellsSub}
Assume that we are in the nonextremal cases.
Let $\mcell\notin\{\jb,\jt\}$ be an apex of $\monoid_{\jcell}^{\kcell}$ 
which is also an apex of $\monoid$. Then we have:
\begin{gather*}
\{\text{simple $\monoid_{\jcell}^{\kcell}$-representations of apex $\mcell$}\}/\cong\;
\xleftrightarrow{1{:}1}
\{\text{simple $\monoid$-representations of apex $\mcell$}\}/\cong
.
\end{gather*}
Moreover, an explicit 
bijection (from left to right) is given by inflating simple $\monoid_{\jcell}^{\kcell}$-representations
to simple $\monoid$-representations.
The dimensions of the simples is preserved under this bijection.

An analog statement holds in the extremal cases.
\end{Proposition}

\begin{proof}
The first part follows from \autoref{P:CellsSimples}. 
For the final part note that inflation clearly does not change 
property of being simple nor the dimension.
\end{proof}

\begin{Theorem}\label{T:CellsSubmonoid}
For any two $J$-cells $\jcell\leq_{lr}\kcell$ we have
\begin{gather*}
\gap{\monoid_{\jcell}^{\kcell}}\geq\gap{\monoid}.
\end{gather*}
\end{Theorem}

\begin{proof}
By \autoref{P:CellsSub}.
\end{proof}

\begin{Remark}[\textbf{\emph{Additional Task 5}}]\label{R:CellsSubmonoid}
By \autoref{T:CellsSubmonoid}, a strategy is to 
find a monoid $\monoid$ with big representations 
for a slice of the cells. 
Then taking an appropriate cell subquotient 
the resulting monoid will have a suitable 
representation gap.
\end{Remark}

We will see examples of the task in \autoref{R:CellsSubmonoid} in the next two sections.


\section{Planar monoids}\label{S:TL}

We work over an arbitrary field $\K$.


\subsection{Temperley--Lieb categories and monoids}\label{SS:TLDef}

We now recall the \emph{Temperley--Lieb category} 
$\tlcat$. This category is a $\K$-linear monoidal category which depends 
on a parameter $\delta\in\K$. There are many references 
(the Temperley--Lieb calculus has been rediscovered many times, 
and there are too many papers to be cited here) for 
$\tlcat$ where more details can be found, see for example \cite{KaLi-TL-recoupling}.
The endomorphism spaces 
in the Temperley--Lieb category form $\K$-algebras, 
called \emph{Temperley--Lieb algebras}.
By appropriate reformulation we obtain set-theoretical 
versions of both of these.

\begin{Remark}\label{R:TLDelta}
It may be convenient to represent 
$\delta=-q-q^{-1}$ where $q$ is either 
in $\K$ or its quadratic extension. For our main application 
we need $\delta=1$, so $q$ in this case is a primitive third root of unity.
This is for example important when one wants to connect 
$\tlcat$ to the category of tilting 
representations for quantum $\mathrm{SL}_{2}$, 
see {\eg} \cite[Proposition 2.28]{TuWe-quiver-tilting} or 
\cite[Proposition 2.20]{SuTuWeZh-mixed-tilting} for a precise statement.
This perspective is sometimes useful, see for example 
\cite{An-simple-tl}, \cite{Sp-modular-tl} or \cite{TuWe-center}
for nontrivial results about the set-theoretical Temperley--Lieb algebras 
using tilting representations.
\end{Remark}

The Temperley--Lieb category $\tlcat$ has objects $n\in\Z_{\geq 0}$. The morphisms from $m$ to $n$ are $\K$-linear combinations of isotopy classes of diagrams of matchings of $m+n$ points in the strip $\R\times[0,1]$, with $m$ points at the bottom and $n$ points at the top line of the strip. These morphisms are known as \emph{crossingless matchings}. The relations on them are such that 
two diagrams represent the same morphism if and only if they 
represent the same crossingless matching.

Composition $\circ$ of crossingless matchings is given by \emph{vertical} gluing 
(and rescaling), using the convention to glue $a\colon m\to n$ on top of $b\colon k\to m$, which is denoted using the operator notation $a\circ b$.
This will give another crossingless matching, but with 
potentially internal circles. To get rid of this ambiguity, we remove 
such internal circles, say we have $k$ of these, and the resulting crossingless matching is multiplied by $\delta^{k}$. 

The monoidal structure $\otimes$ is given by $m\otimes n=m+n$ on objects and \emph{horizontal} juxtaposition on morphisms, extended bilinearly to $\K$-linear combinations.

\begin{Notation}\label{N:TLReading}
The following pictures summarize the main points from above, 
and also fixes the reading conventions that we will use for diagrammatics 
throughout.
\begin{gather*}
\begin{tikzpicture}[anchorbase]
\draw[usual] (0,0) to[out=90,in=180] (0.25,0.25) to[out=0,in=90] (0.5,0);
\draw[usual] (1,1) to[out=270,in=180] (1.25,0.75) to[out=0,in=270] (1.5,1);
\draw[usual] (1,0) to (0,1);
\draw[usual] (1.5,0) to (0.5,1);
\end{tikzpicture}
\;\circ\;
\begin{tikzpicture}[anchorbase]
\draw[usual] (1,0) to[out=90,in=180] (1.25,0.25) to[out=0,in=90] (1.5,0);
\draw[usual] (0,1) to[out=270,in=180] (0.25,0.75) to[out=0,in=270] (0.5,1);
\draw[usual] (0,0) to (1,1);
\draw[usual] (0.5,0) to (1.5,1);
\end{tikzpicture}
\;=\;
\begin{tikzpicture}[anchorbase]
\draw[usual] (1,0) to[out=90,in=180] (1.25,0.25) to[out=0,in=90] (1.5,0);
\draw[usual] (0,1) to[out=270,in=180] (0.25,0.75) to[out=0,in=270] (0.5,1);
\draw[usual] (0,0) to (1,1);
\draw[usual] (0.5,0) to (1.5,1);
\draw[usual] (0,1) to[out=90,in=180] (0.25,1.25) to[out=0,in=90] (0.5,1);
\draw[usual] (1,2) to[out=270,in=180] (1.25,1.75) to[out=0,in=270] (1.5,2);
\draw[usual] (1,1) to (0,2);
\draw[usual] (1.5,1) to (0.5,2);
\end{tikzpicture}
\;=\;\delta\cdot
\begin{tikzpicture}[anchorbase]
\draw[usual] (0,0) to (0,1);
\draw[usual] (0.5,0) to (0.5,1);
\draw[usual] (1,0) to[out=90,in=180] (1.25,0.25) to[out=0,in=90] (1.5,0);
\draw[usual] (1,1) to[out=270,in=180] (1.25,0.75) to[out=0,in=270] (1.5,1);
\end{tikzpicture}
\;,
\\[10pt]
\begin{tikzpicture}[anchorbase]
\draw[usual] (0,0) to[out=90,in=180] (0.75,0.5) to[out=0,in=90] (1.5,0);
\draw[usual] (0.5,0) to[out=90,in=180] (0.75,0.25) to[out=0,in=90] (1,0);
\draw[usual] (0,1) to[out=270,in=180] (0.25,0.75) to[out=0,in=270] (0.5,1);
\draw[usual] (1,1) to[out=270,in=180] (1.25,0.75) to[out=0,in=270] (1.5,1);
\end{tikzpicture}
\;\circ\;
\begin{tikzpicture}[anchorbase]
\draw[usual] (0,0) to[out=90,in=180] (0.75,0.45) to[out=0,in=90] (1.5,0);
\draw[usual] (0,1) to[out=270,in=180] (0.75,0.55) to[out=0,in=270] (1.5,1);
\draw[usual] (0.5,0) to[out=90,in=180] (0.75,0.25) to[out=0,in=90] (1,0);
\draw[usual] (0.5,1) to[out=270,in=180] (0.75,0.75) to[out=0,in=270] (1,1);
\end{tikzpicture}
\;=\;
\begin{tikzpicture}[anchorbase]
\draw[usual] (0,0) to[out=90,in=180] (0.75,0.45) to[out=0,in=90] (1.5,0);
\draw[usual] (0,1) to[out=270,in=180] (0.75,0.55) to[out=0,in=270] (1.5,1);
\draw[usual] (0.5,0) to[out=90,in=180] (0.75,0.25) to[out=0,in=90] (1,0);
\draw[usual] (0.5,1) to[out=270,in=180] (0.75,0.75) to[out=0,in=270] (1,1);
\draw[usual] (0,1) to[out=90,in=180] (0.75,1.5) to[out=0,in=90] (1.5,1);
\draw[usual] (0.5,1) to[out=90,in=180] (0.75,1.25) to[out=0,in=90] (1,1);
\draw[usual] (0,2) to[out=270,in=180] (0.25,1.75) to[out=0,in=270] (0.5,2);
\draw[usual] (1,2) to[out=270,in=180] (1.25,1.75) to[out=0,in=270] (1.5,2);
\end{tikzpicture}
\;=\;
\delta^{2}\cdot
\begin{tikzpicture}[anchorbase]
\draw[usual] (0,0) to[out=90,in=180] (0.75,0.5) to[out=0,in=90] (1.5,0);
\draw[usual] (0.5,0) to[out=90,in=180] (0.75,0.25) to[out=0,in=90] (1,0);
\draw[usual] (0,1) to[out=270,in=180] (0.25,0.75) to[out=0,in=270] (0.5,1);
\draw[usual] (1,1) to[out=270,in=180] (1.25,0.75) to[out=0,in=270] (1.5,1);
\end{tikzpicture}
\;,
\\[10pt]
\begin{tikzpicture}[anchorbase]
\draw[usual] (0,0) to[out=90,in=180] (0.75,0.45) to[out=0,in=90] (1.5,0);
\draw[usual] (0,1) to[out=270,in=180] (0.75,0.55) to[out=0,in=270] (1.5,1);
\draw[usual] (0.5,0) to[out=90,in=180] (0.75,0.25) to[out=0,in=90] (1,0);
\draw[usual] (0.5,1) to[out=270,in=180] (0.75,0.75) to[out=0,in=270] (1,1);
\end{tikzpicture}
\;\otimes\;
\begin{tikzpicture}[anchorbase]
\draw[usual] (0.5,0) to[out=90,in=180] (0.75,0.25) to[out=0,in=90] (1,0);
\draw[usual] (0,1) to[out=270,in=180] (0.25,0.75) to[out=0,in=270] (0.5,1);
\draw[usual] (1,1) to[out=270,in=180] (1.25,0.75) to[out=0,in=270] (1.5,1);
\end{tikzpicture}
\;=\;
\begin{tikzpicture}[anchorbase]
\draw[usual] (0,0) to[out=90,in=180] (0.75,0.45) to[out=0,in=90] (1.5,0);
\draw[usual] (0,1) to[out=270,in=180] (0.75,0.55) to[out=0,in=270] (1.5,1);
\draw[usual] (0.5,0) to[out=90,in=180] (0.75,0.25) to[out=0,in=90] (1,0);
\draw[usual] (0.5,1) to[out=270,in=180] (0.75,0.75) to[out=0,in=270] (1,1);
\draw[usual] (2.5,0) to[out=90,in=180] (2.75,0.25) to[out=0,in=90] (3,0);
\draw[usual] (2,1) to[out=270,in=180] (2.25,0.75) to[out=0,in=270] (2.5,1);
\draw[usual] (3,1) to[out=270,in=180] (3.25,0.75) to[out=0,in=270] (3.5,1);
\end{tikzpicture}
\;.
\end{gather*}
\end{Notation}

Let $C_{m}^{n}$ denote the set of crossingless matching with $m$ bottom
and $n$ top boundary points. Let $Ca(k)=\tfrac{1}{k+1}\binom{2k}{k}$ be the $k$th Catalan number. 
Note that the following lemma is 
independent of $\K$ and $\delta\in\K$.

\begin{Lemma}\label{L:TLBasis}
The set $C_{m}^{n}$ is a $\K$-linear basis of $\Hom_{\tlcat}(m,n)$.
Hence, the dimension of this space is either zero if $m\not\equiv n\bmod 2$, and otherwise given by $\dimk\big(\Hom_{\tlcat}(m,n)\big)=Ca(\tfrac{m+n}{2})$.
\end{Lemma}

\begin{proof}
This is well-known, see {\eg} \cite{RuTeWe-sl2} for the version with $\delta=-2$.
\end{proof}

\begin{Lemma}\label{L:TLInvolution}
The category $\tlcat$ has an antiinvolution ${\placeholder}^{\ast}$, {\ie} is a $\ast$-monoid, given by 
reflecting diagrams in a horizontal axis.
\end{Lemma}

\begin{proof}
Easy and omitted.
\end{proof}

The picture to keep in mind is
\begin{gather*}
\left(
\begin{tikzpicture}[anchorbase]
\draw[usual] (0.5,1) to[out=270,in=180] (0.75,0.75) to[out=0,in=270] (1,1);
\draw[usual] (-1,1) to[out=270,in=180] (-0.75,0.75) to[out=0,in=270] (-0.5,1);
\draw[usual] (0,0) to (0,1);
\draw[usual] (1.5,0) to (1.5,1);
\draw[usual] (2,0) to (2,1);
\draw[usual] (2.5,0) to[out=90,in=180] (2.75,0.25) to[out=0,in=90] (3,0);
\end{tikzpicture}
\right)^{\ast}
=
\begin{tikzpicture}[anchorbase,yscale=-1]
\draw[usual] (0.5,1) to[out=270,in=180] (0.75,0.75) to[out=0,in=270] (1,1);
\draw[usual] (-1,1) to[out=270,in=180] (-0.75,0.75) to[out=0,in=270] (-0.5,1);
\draw[usual] (0,0) to (0,1);
\draw[usual] (1.5,0) to (1.5,1);
\draw[usual] (2,0) to (2,1);
\draw[usual] (2.5,0) to[out=90,in=180] (2.75,0.25) to[out=0,in=90] (3,0);
\end{tikzpicture}
\;.
\end{gather*}

\begin{Remark}\label{R:TLInvolution}
It is easy to see (and we will use this silently) that 
${\placeholder}^{\ast}$ works for all the diagrammatic categories, 
algebras and monoids we use in this and the next section.
We call ${\placeholder}^{\ast}$ the \emph{diagrammatic antiinvolution}.
\end{Remark}

The Temperley--Lieb algebra on $n$-strands is then $\tlalg{\delta}=\End_{\tlcat}(n)$. This is the algebra of crossingless matchings with $n$ strands and only vertical composition.

\begin{Remark}\label{R:TLSchurWeyl}
The algebra $\tlalg{\delta}$ was introduced in the context of Schur--Weyl 
duality, see \cite{RuTeWe-sl2}.
Sometimes it is useful to use this perspective as {\eg} 
the reference \cite{An-simple-tl} does (using the 
connection to tilting representations, {\cf} \autoref{R:TLDelta}) which we will use below.
\end{Remark}

Now comes the main definition of this section.

\begin{Definition}\label{D:TLSet}
The \emph{set-theoretic Temperley--Lieb category} $\tlset$ 
is defined in almost the same way as $\tlcat$ above 
with two crucial differences:
\begin{enumerate}

\item The hom-spaces are $\Hom_{\tlset}(m,n)=C_{m}^{n}$, and,

\item the vertical composition $\circ$ is still given by vertical 
gluing, but all internal circles are just removed from the diagram, that is, without any factor.

\end{enumerate}
The \emph{Temperley--Lieb monoid} on $n$-strands is defined by $\tlmon=\End_{\tlset}(n)$. 
\end{Definition}

\begin{Remark}\label{R:TLMonVsLinear}
The Temperley--Lieb monoid appears in many works, 
way too many to be cited here, see however {\eg} 
\cite{HaRa-partition-algebras}, 
\cite{HaJa-representations-diagram-algebras} or
\cite{Si-topological-tl-actions}. In most papers coming 
from representation theory, quantum algebra and quantum topology
it is however studied as an algebra.
Note hereby that \autoref{D:TLSet} is not quite the same as $\tlcat[1]$ 
where the circle evaluates to $1$. The difference is that $\tlcat[1]$ 
is $\K$-linear, but $\tlset$ is not $\K$-linear. 
But the monoid algebra $\K[\tlmon]$ is isomorphic to 
$\tlalg{1}$.
Let us stress that the Temperley--Lieb monoid 
is also called the Jones monoid in monoid theory, or sometimes even the Kauffman monoid, see {\eg} \cite{LaFiGe-ideal-kauffman}.
\end{Remark}

The monoid $\tlmon$ has $Ca(n)$ elements. 
By \cite[Section 2.2]{KaLi-TL-recoupling} (or \cite{Ea-tl-presentation} 
for a new proof of the presentation), 
the monoid $\tlmon$ can be abstractly defined by the generators 
$\{u_{1},\dots, u_{n-1}\}$ and the defining relations 
\begin{gather}\label{Eq:TLPresentation}
u_{i}^{2}=u_{i}
,\quad
u_{i}u_{i\pm 1}u_{i}=u_{i}
,\quad
u_{i}u_{j}=u_{j}u_{i},\ |i-j|>1
.
\end{gather}

Denote by $\id_{k}$ the identity on $k\in\N$. The following determines the cell structure:

\begin{Lemma}\label{L:TLFactor}
For $a\in\Hom_{\tlset}(m,n)$ there is a unique factorization of the form $a=\gamma\circ\id_{k}\circ\beta$ for minimal $k$, and $\beta\in\Hom_{\tlset}(m,k)$ and $\gamma\in\Hom_{\tlset}(k,n)$.
\end{Lemma}

\begin{proof}
The following picture
\begin{gather}\label{Eq:TLFactor}
a=
\begin{tikzpicture}[anchorbase]
\draw[usual] (1,0) to[out=90,in=180] (1.25,0.25) to[out=0,in=90] (1.5,0);
\draw[usual] (0.5,1) to[out=270,in=180] (0.75,0.75) to[out=0,in=270] (1,1);
\draw[usual] (-1,1) to[out=270,in=180] (-0.75,0.75) to[out=0,in=270] (-0.5,1);
\draw[usual] (0,0) to (0,1);
\draw[usual] (0.5,0) to (1.5,1);
\end{tikzpicture}
\;=\;
\underbrace{\begin{tikzpicture}[anchorbase]
\draw[usual] (0.5,1) to[out=270,in=180] (0.75,0.75) to[out=0,in=270] (1,1);
\draw[usual] (-1,1) to[out=270,in=180] (-0.75,0.75) to[out=0,in=270] (-0.5,1);
\draw[usual] (0,0) to (0,1);
\draw[usual] (1.5,0) to (1.5,1);
\end{tikzpicture}}_{\gamma}
\;\circ\;
\underbrace{\begin{tikzpicture}[anchorbase]
\draw[usual] (0,0) to (0,1);
\draw[usual] (0.5,0) to (0.5,1);
\end{tikzpicture}}_{\id_{2}}
\;\circ\;
\underbrace{\begin{tikzpicture}[anchorbase]
\draw[usual] (1,0) to[out=90,in=180] (1.25,0.25) to[out=0,in=90] (1.5,0);
\draw[usual] (0,0) to (0,1);
\draw[usual] (0.5,0) to (0.5,1);
\end{tikzpicture}}_{\beta}
\;,
\end{gather}
generalizes without much work.
\end{proof}

We call $k$ the number of \emph{through strands} of $\alpha$, also 
known as the \emph{width}. Necessarily $k$ has the same parity as $m$ and $n$ and $k\leq m,n$. The diagram $a$ has $\tfrac{m-k}{2}$ caps and $\tfrac{n-k}{2}$ cups. The diagrams $\beta$ and $\gamma$ have no cups, respectively no caps, but the same number of caps, respectively cups, as $\alpha$. We call $\beta$ as in \autoref{L:TLFactor} the 
\emph{bottom half} and $\gamma$ the \emph{top half} of $a$.

Denote by $\nobottom{m}{n}\subset\Hom_{\tlset}(m,n)$ the set of diagrams 
without caps. An example for $m=2$ and $n=6$ is given by 
$\gamma$ in \autoref{Eq:TLFactor}. 
In other words,
$\nobottom{m}{n}$ consists of $m$ 
\emph{through strands} and $\tfrac{n-m}{2}$ cups.
Necessarily $m\leq n$ and $m+n$ is even. 
In the above factorization, in general, 
$\gamma,\beta^{\ast}\in\nobottom{m}{n}$. 
We may also write this factorization of $a$ as $a=a_{1}a_{2}^{\ast}$, 
$a_{1},a_{2}\in\nobottom{m}{n}$. 


\subsection{Cells of the Temperley--Lieb monoid}\label{SS:TLCell}

We now discuss the cell structure of $\tlmon$.

\begin{Remark}\label{R:TLCells}
The cell structure of the Temperley--Lieb monoid $\tlmon$ 
is very nice and easy to compute. It is well-known, see {\eg} \cite[Example 1.4]{GrLe-cellular},
and was rediscovered in many papers, see {\eg} \cite{RiSaAu-temperley-lieb}, 
or \cite{KhSa-cat-chebychev} or \cite{Sp-modular-tl}, 
although not always in the language of cells. 
The cell structure has also been 
rediscovered in monoid theory, see {\eg} \cite{LaFiGe-ideal-kauffman}.
In any case, the description of the cells 
is prototypical for diagram monoids and algebras
so we decided to repeat it here in that language.
\end{Remark}

The main pictures to keep in mind (which we will 
explain momentarily) are:
\begin{gather}\label{Eq:TLCells}
\begin{gathered}
\xy
(0,0)*{\begin{gathered}
\begin{tabular}{C|C}
\arrayrulecolor{tomato}
\cellcolor{mydarkblue!25}
\begin{tikzpicture}[anchorbase]
\draw[usual] (0,0) to[out=90,in=180] (0.25,0.2) to[out=0,in=90] (0.5,0);
\draw[usual] (0,0.5) to[out=270,in=180] (0.25,0.3) to[out=0,in=270] (0.5,0.5);
\draw[usual] (1,0) to (1,0.5);
\end{tikzpicture} &
\cellcolor{mydarkblue!25}
\begin{tikzpicture}[anchorbase]
\draw[usual] (0.5,0) to[out=90,in=180] (0.75,0.2) to[out=0,in=90] (1,0);
\draw[usual] (0,0.5) to[out=270,in=180] (0.25,0.3) to[out=0,in=270] (0.5,0.5);
\draw[usual] (0,0) to (1,0.5);
\end{tikzpicture}
\\
\hline
\cellcolor{mydarkblue!25}
\begin{tikzpicture}[anchorbase,xscale=-1]
\draw[usual] (0.5,0) to[out=90,in=180] (0.75,0.2) to[out=0,in=90] (1,0);
\draw[usual] (0,0.5) to[out=270,in=180] (0.25,0.3) to[out=0,in=270] (0.5,0.5);
\draw[usual] (0,0) to (1,0.5);
\end{tikzpicture} &
\cellcolor{mydarkblue!25}
\begin{tikzpicture}[anchorbase,xscale=-1]
\draw[usual] (0,0) to[out=90,in=180] (0.25,0.2) to[out=0,in=90] (0.5,0);
\draw[usual] (0,0.5) to[out=270,in=180] (0.25,0.3) to[out=0,in=270] (0.5,0.5);
\draw[usual] (1,0) to (1,0.5);
\end{tikzpicture}
\end{tabular}
\\[3pt]
\begin{tabular}{C}
\arrayrulecolor{tomato}
\cellcolor{mydarkblue!25}
\begin{tikzpicture}[anchorbase]
\draw[usual] (0,0) to (0,0.5);
\draw[usual] (0.5,0) to (0.5,0.5);
\draw[usual] (1,0) to (1,0.5);
\end{tikzpicture}
\end{tabular}
\end{gathered}};
(-45,4)*{\jcell_{1}};
(-45,-7)*{\jcell_{3}};
(45,4)*{\hcell(e)\cong\onemon};
(45,-7)*{\hcell(e)\cong\onemon};
(-51,0)*{\phantom{a}};
\endxy
\quad,
\\[0.27cm]
\hline
\\[-0.27cm]
\xy
(0,0)*{\begin{gathered}
\begin{tabular}{C|C}
\arrayrulecolor{tomato}
\cellcolor{mydarkblue!25}
\begin{tikzpicture}[anchorbase]
\draw[usual] (0,0) to[out=90,in=180] (0.25,0.2) to[out=0,in=90] (0.5,0);
\draw[usual] (0,0.5) to[out=270,in=180] (0.25,0.3) to[out=0,in=270] (0.5,0.5);
\draw[usual] (1,0) to[out=90,in=180] (1.25,0.2) to[out=0,in=90] (1.5,0);
\draw[usual] (1,0.5) to[out=270,in=180] (1.25,0.3) to[out=0,in=270] (1.5,0.5);
\end{tikzpicture} &
\cellcolor{mydarkblue!25}
\begin{tikzpicture}[anchorbase]
\draw[usual] (0,0) to[out=45,in=180] (0.75,0.20) to[out=0,in=135] (1.5,0);
\draw[usual] (0.5,0) to[out=90,in=180] (0.75,0.1) to[out=0,in=90] (1,0);
\draw[usual] (0,0.5) to[out=270,in=180] (0.25,0.3) to[out=0,in=270] (0.5,0.5);
\draw[usual] (1,0.5) to[out=270,in=180] (1.25,0.3) to[out=0,in=270] (1.5,0.5);
\end{tikzpicture}
\\
\hline
\cellcolor{mydarkblue!25}
\begin{tikzpicture}[anchorbase]
\draw[usual] (0,0) to[out=90,in=180] (0.25,0.2) to[out=0,in=90] (0.5,0);
\draw[usual] (1,0) to[out=90,in=180] (1.25,0.2) to[out=0,in=90] (1.5,0);
\draw[usual] (0,0.5) to[out=315,in=180] (0.75,0.3) to[out=0,in=225] (1.5,0.5);
\draw[usual] (0.5,0.5) to[out=270,in=180] (0.75,0.4) to[out=0,in=270] (1,0.5);
\end{tikzpicture} &
\cellcolor{mydarkblue!25}
\begin{tikzpicture}[anchorbase]
\draw[usual] (0,0) to[out=45,in=180] (0.75,0.20) to[out=0,in=135] (1.5,0);
\draw[usual] (0,0.5) to[out=315,in=180] (0.75,0.3) to[out=0,in=225] (1.5,0.5);
\draw[usual] (0.5,0) to[out=90,in=180] (0.75,0.1) to[out=0,in=90] (1,0);
\draw[usual] (0.5,0.5) to[out=270,in=180] (0.75,0.4) to[out=0,in=270] (1,0.5);
\end{tikzpicture}
\end{tabular}
\\[3pt]
\begin{tabular}{C|C|C}
\arrayrulecolor{tomato}
\cellcolor{mydarkblue!25}
\begin{tikzpicture}[anchorbase]
\draw[usual] (0,0) to[out=90,in=180] (0.25,0.2) to[out=0,in=90] (0.5,0);
\draw[usual] (0,0.5) to[out=270,in=180] (0.25,0.3) to[out=0,in=270] (0.5,0.5);
\draw[usual] (1,0) to (1,0.5);
\draw[usual] (1.5,0) to (1.5,0.5);
\end{tikzpicture} & 
\cellcolor{mydarkblue!25}
\begin{tikzpicture}[anchorbase]
\draw[usual] (0.5,0) to[out=90,in=180] (0.75,0.2) to[out=0,in=90] (1,0);
\draw[usual] (0,0.5) to[out=270,in=180] (0.25,0.3) to[out=0,in=270] (0.5,0.5);
\draw[usual] (0,0) to (1,0.5);
\draw[usual] (1.5,0) to (1.5,0.5);
\end{tikzpicture} &
\begin{tikzpicture}[anchorbase]
\draw[usual] (1,0) to[out=90,in=180] (1.25,0.2) to[out=0,in=90] (1.5,0);
\draw[usual] (0,0.5) to[out=270,in=180] (0.25,0.3) to[out=0,in=270] (0.5,0.5);
\draw[usual] (0,0) to (1,0.5);
\draw[usual] (0.5,0) to (1.5,0.5);
\end{tikzpicture}
\\
\hline
\cellcolor{mydarkblue!25}
\begin{tikzpicture}[anchorbase]
\draw[usual] (0,0) to[out=90,in=180] (0.25,0.2) to[out=0,in=90] (0.5,0);
\draw[usual] (0.5,0.5) to[out=270,in=180] (0.75,0.3) to[out=0,in=270] (1,0.5);
\draw[usual] (1,0) to (0,0.5);
\draw[usual] (1.5,0) to (1.5,0.5);
\end{tikzpicture} & 
\cellcolor{mydarkblue!25}
\begin{tikzpicture}[anchorbase]
\draw[usual] (0,0) to (0,0.5);
\draw[usual] (0.5,0) to[out=90,in=180] (0.75,0.2) to[out=0,in=90] (1,0);
\draw[usual] (0.5,0.5) to[out=270,in=180] (0.75,0.3) to[out=0,in=270] (1,0.5);
\draw[usual] (1.5,0) to (1.5,0.5);
\end{tikzpicture} &
\cellcolor{mydarkblue!25}
\begin{tikzpicture}[anchorbase]
\draw[usual] (0,0) to (0,0.5);
\draw[usual] (0.5,0) to (1.5,0.5);
\draw[usual] (1,0) to[out=90,in=180] (1.25,0.2) to[out=0,in=90] (1.5,0);
\draw[usual] (0.5,0.5) to[out=270,in=180] (0.75,0.3) to[out=0,in=270] (1,0.5);
\end{tikzpicture}
\\
\hline
\begin{tikzpicture}[anchorbase]
\draw[usual] (0,0) to[out=90,in=180] (0.25,0.2) to[out=0,in=90] (0.5,0);
\draw[usual] (1,0.5) to[out=270,in=180] (1.25,0.3) to[out=0,in=270] (1.5,0.5);
\draw[usual] (1,0) to (0,0.5);
\draw[usual] (1.5,0) to (0.5,0.5);
\end{tikzpicture} & 
\cellcolor{mydarkblue!25}
\begin{tikzpicture}[anchorbase]
\draw[usual] (0,0) to (0,0.5);
\draw[usual] (1.5,0) to (0.5,0.5);
\draw[usual] (0.5,0) to[out=90,in=180] (0.75,0.2) to[out=0,in=90] (1,0);
\draw[usual] (1,0.5) to[out=270,in=180] (1.25,0.3) to[out=0,in=270] (1.5,0.5);
\end{tikzpicture} &
\cellcolor{mydarkblue!25}
\begin{tikzpicture}[anchorbase]
\draw[usual] (0,0) to (0,0.5);
\draw[usual] (0.5,0) to (0.5,0.5);
\draw[usual] (1,0) to[out=90,in=180] (1.25,0.2) to[out=0,in=90] (1.5,0);
\draw[usual] (1,0.5) to[out=270,in=180] (1.25,0.3) to[out=0,in=270] (1.5,0.5);
\end{tikzpicture}
\end{tabular}
\\[3pt]
\begin{tabular}{C}
\arrayrulecolor{tomato}
\cellcolor{mydarkblue!25}
\begin{tikzpicture}[anchorbase]
\draw[usual] (0,0) to (0,0.5);
\draw[usual] (0.5,0) to (0.5,0.5);
\draw[usual] (1,0) to (1,0.5);
\draw[usual] (1.5,0) to (1.5,0.5);
\end{tikzpicture}
\end{tabular}
\end{gathered}};
(-45,13)*{\jcell_{0}};
(-45,-2.5)*{\jcell_{2}};
(-45,-16.5)*{\jcell_{4}};
(45,13)*{\hcell(e)\cong\onemon};
(45,-2.5)*{\hcell(e)\cong\onemon};
(45,-16.5)*{\hcell(e)\cong\onemon};
(-51,0)*{\phantom{a}};
\endxy
\quad.
\end{gathered}
\end{gather}
These are the cells of $\tlmon[3]$ and $\tlmon[4]$, which should be read as in 
\autoref{Eq:CellsIllustration}. We have also colored/shaded the idempotent $H$-cells. Note that $\jcell_{k}$ is the set of crossingless matchings with $k$ through strands, and $k$ and $n$ have the same parity. These diagrams have $c(k)=\tfrac{n-k}{2}$ caps respectively cups.

\begin{Proposition}\label{P:TLCells}
We have the following.
\begin{enumerate}

\item The left and right cells of $\tlmon$ 
are given by crossingless matchings where one fixes 
the bottom respectively top half of the diagram.
The $\leq_{l}$- and the $\leq_{r}$-order increases as the number of through strands 
decreases. Within $\jcell_{k}$ we have
\begin{gather*}
|\lcell|=|\rcell|=\tfrac{n-2c(k)+1}{n-c(k)+1}\binom{n}{c(k)}.
\end{gather*}

\item The $J$-cells $\jcell_{k}$ of $\tlmon$
are given by crossingless matchings with a fixed number of through strands $k$.
The $\leq_{lr}$-order is a total order and increases as the number of through strands 
decreases. For any $\lcell\subset\jcell_{k}$ we have
\begin{gather*}
|\jcell_{k}|=|\lcell|^{2}.
\end{gather*}

\item Each $J$-cell of $\tlmon$ is idempotent, and $\hcell(e)\cong\onemon$ for all idempotent $H$-cells. We have
\begin{gather*}
|\hcell|=1.
\end{gather*}

\end{enumerate}	

\end{Proposition}

\begin{proof}
\textit{(a)+(b).}
For (a) and (b) we recall that
the $\K$-linear version of this proposition can be found in 
{\eg} \cite[Example 1.4]{GrLe-cellular} 
or \cite[Section 2]{RiSaAu-temperley-lieb}. 
(Note that \cite[Section 2]{RiSaAu-temperley-lieb}
gives $|\lcell|=|\rcell|=\binom{n}{c(k)}-\binom{n}{c(k)-1}$, which 
we rewrite into the claimed expression 
via algebra autopilot.)
The arguments 
given in these papers do not depend on $\K$ nor on the parameter $\delta$ and go through in the set-theoretical case without change as well.
In monoid theory this appears again in many works, {\eg} in 
\cite{LaFiGe-ideal-kauffman}.
\medskip

\textit{(c).}
Observing that every crossingless matching that is 
symmetric under horizontal mirroring is an idempotent, this is then immediate from (a) and (b).
\end{proof}

\begin{Proposition}\label{P:TLSimples}
The set of apexes for simple 
$\tlmon$-representations can be indexed $1{:}1$ by the poset
$\Lambda=(\{n,n-2,\dots\},>)$ (ending on either $0$ or $1$, depending 
on the parity of $n$), and
there is precisely one simple 
$\tlmon$-representation of a fixed apex up to $\cong$.
\end{Proposition}

\begin{proof}
By \autoref{P:TLCells}, this is a direct application of \autoref{P:CellsSimples}.
\end{proof}

By \autoref{P:TLSimples} there is 
a poset
$\Lambda$ indexing the $J$-cells 
and the simple 
$\tlmon$-representations.
We can thus enumerate the 
$J$-cells by $\jcell_{k}$ for $k\in\Lambda$. 
We do the same for the simple 
$\tlmon$-representations and we write $\simple[k]$ for 
these. (Here we mean any choice 
of representatives of the isomorphism classes. Similarly 
below, and we stop stressing this.)

\begin{Lemma}\label{L:TLCellModules}
Within one $J$-cell, all left cell 
representations $\lmod$ and all right cell representations $\rmod$ are isomorphic. We write $\lmod[k]$ respectively $\rmod[k]$ for those in $\jcell_{k}$.

We have $\lmod[k]\cong\rmod[k]$ as $\K$-vector spaces
and $\dimk(\lmod[k])=\dimk(\rmod[k])=\tfrac{n-2c(k)+1}{n-c(k)+1}\binom{n}{c(k)}$.
\end{Lemma}

\begin{proof}
The diagrammatic antiinvolution ${\placeholder}^{\ast}$ 
is compatible with the cells structure and shows 
$\lmod[k]\cong\rmod[k]$. The dimension formula then follows 
from \autoref{P:TLCells} and \autoref{L:CellsMod}.
\end{proof}

\begin{Proposition}\label{P:TLSemisimpleDims}
The semisimple dimensions are $\ssdimk(\simple[k])=\tfrac{n-2c(k)+1}{n-c(k)+1}\binom{n}{c(k)}$.
\end{Proposition}

\begin{proof}
The equation follows immediately from \autoref{P:TLCells} and \autoref{P:TLSimples}.
\end{proof}

The numbers $\dimk(\simple[k])$ are as follows. 
These were computed in many papers, {\eg} in \cite{An-simple-tl} and 
\cite{Sp-modular-tl} which compute them for general $\K$ and 
$\delta\in\K$. (Strictly speaking \cite{An-simple-tl} needs 
$\delta=-q-q^{-1}$ because Andersen uses the connection to tilting representations 
as recalled in \autoref{R:TLDelta}.)
To state them 
we need some preliminary definitions.

\begin{Remark}\label{R:TLPAdicNotation}
The definitions below are fairly standard for 
Temperley--Lieb calculi over arbitrary fields, see {\eg} \cite{Sp-modular-tl}, \cite{Sp-valenced-tl} 
or \cite{SuTuWeZh-mixed-tilting}. 
The reader only interested in $\cchar=0$ (which is $\cchar=\infty$ below)
can ignore all definitions involving $p$-adic combinatorics. We elaborate on the $\cchar=0$ case in \autoref{E:TLDimsChar0} below.
\end{Remark}

Let $\cchar=p$, allowing $p=\infty$ 
which is the case $\cchar=0$. Let $\nu_{p}$ 
denote the \emph{$p$-adic valuation}. 
Let $\nu_{3,p}(x)=0$ if $x\not\equiv 0\bmod 3$, 
and $\nu_{3,p}(x)=\nu_{p}(\tfrac{x}{3})$ otherwise. Let further 
$x=[\dots,x_{1},x_{0}]$ denote the \emph{$(3,p)$-adic expansion} of $x$ 
given by
\begin{gather*}
[\dots,x_{1},x_{0}]=\sum_{i=1}^{\infty}3p^{i-1}x_{i}+x_{0}=x,
\quad
x_{i>0}\in\{0,\dots,p-1\},x_{0}\in\{0,1,2\}
.
\end{gather*}
The numbers $x_{j}$ are the digits of $x$, and 
most of these $x_{j}$ are zero.
Let now $x\tlord y$ if $[\dots,x_{1},x_{0}]$ is digit-wise
smaller or equal to $[\dots,y_{1},y_{0}]$. We also write $x\tlordsecond y$ if 
$x\tlord y$, $\nu_{3,p}(x)=\nu_{3,p}(y)$ and the $\nu_{3,p}(x)$th digit 
of $x$ and $y$ agree. Finally, set
\begin{gather}\label{Eq:TLDimENumber}
e_{n,k}=
\begin{cases}
1&\text{if }n\equiv k\bmod 2,\ \nu_{3,p}(k)=\nu_{3,p}(\tfrac{n+k}{2}),\ k\tlordsecond\tfrac{n+k}{2},
\\
-1&\text{if }n\equiv k\bmod 2,\ \nu_{3,p}(k)<\nu_{3,p}(\tfrac{n+k}{2}),
k\tlord\tfrac{n+k}{2}-1,
\\
0&\text{else.}
\end{cases}
\end{gather}

\begin{Example}\label{E:TLDimsChar0}
For $\cchar=0$ the above simplifies quite a bit. First, the only two relevant numbers $x_{1}\in\N,x_{0}\in\{0,1,2\}$ are given by $x=3x_{1}+x_{0}$, so 
$x_{0}$ is the reminder of $x$ upon division by $3$. 
The equation \autoref{Eq:TLDimENumber} simplifies to the following matrix whose entries are $e_{n,k}$:
\begin{gather*}
\scalebox{0.6}{\begin{tikzpicture}[anchorbase, ampersand replacement=\&]
\matrix (A) [matrix of nodes,every node/.style={anchor=base,text depth=0.5ex,text height=1ex,text width=2em}]
{
n\text{\textbackslash}k\& 0 \& 1 \& 2 \& 3 \& 4 \& 5 \& 6 \& 7 \& 8 \& 9 \& 10 \& 11 \& 12 \& 13 \& 14 \& 15 \& 16
\\
0 \& 1 \& \& \& \& \& \& \& \& \& \& \& \& \& \& \& \&
\\
1 \& \& 1 \& \& \& \& \& \& \& \& \& \& \& \& \& \& \&
\\
2 \& -1 \& \& 1 \& \& \& \& \& \& \& \& \& \& \& \& \& \&
\\
3 \& \& \& \& 1 \& \& \& \& \& \& \& \& \& \& \& \& \&
\\
4 \& \& \& \& \& 1 \& \& \& \& \& \& \& \& \& \& \& \&
\\
5 \& \& \& \& -1 \& \& 1 \& \& \& \& \& \& \& \& \& \& \&
\\
6 \& 1 \& \& -1 \& \& \& \& 1 \& \& \& \& \& \& \& \& \& \&	
\\
7 \& \& \& \& \& \& \& \& 1 \& \& \& \& \& \& \& \& \&
\\
8 \& -1 \& \& 1 \& \& \& \& -1 \& \& 1 \& \& \& \& \& \& \& \&
\\
9 \& \& \& \& 1 \& \& -1 \& \& \& \& 1 \& \& \& \& \& \& \&
\\
10 \& \& \& \& \& \& \& \& \& \& \& 1 \& \& \& \& \& \&
\\
11 \& \& \& \& -1 \& \& 1 \& \& \& \& -1 \& \& 1 \& \& \& \& \&
\\
12 \& 1 \& \& -1 \& \& \& \& 1 \& \& -1 \& \& \& \& 1 \& \& \& \&
\\
13 \& \& \& \& \& \& \& \& \& \& \& \& \& \& 1 \& \& \&
\\
14 \& -1 \& \& 1 \& \& \& \& -1 \& \& 1 \& \& \& \& -1 \& \& 1 \& \&	
\\
15 \& \& \& \& 1 \& \& -1 \& \& \& \& 1 \& \& -1 \& \& \& \& 1 \&
\\
16 \& \& \& \& \& \& \& \& \& \& \& \& \& \& \& \& \& 1
\\
};
\draw [very thick] ($(A-1-1.south west)+(0,0.05)$) to ($(A-1-18.south east)+(0,0.05)$);
\draw [very thick] ($(A-1-1.north east)+(-0.1,0)$) to ($(A-18-1.south east)+(-0.1,0)$);
\draw [very thick,->] (1.5,-0.15)node[right]{$e_{8,8}$} to (0.5,-0.15);
\draw [very thick,->] (5.7,-1.5)node[right]{$e_{15,11}$} to (3.5,-3.5);
\end{tikzpicture}}
\quad.
\end{gather*}
Here we have illustrated the case $n=16$. The pattern is that every 
third row has only one nonzero entry. Otherwise, the pattern $(-1,0,1)$ 
respectively $(1,0,-1)$ 
is shifted along rows with a distance of three zeros.
\end{Example}

We have the following alternating sum of $\dimk(\lmod)=\tfrac{n-2c(k)+1}{n-c(k)+1}\binom{n}{c(k)}$. 
(Recall that $c(k)$ denotes the number of caps respectively cups 
for diagrams in the $J$-cell $\jcell_{k}$.) 
That a dimension formula is of this form is expected 
from the cell structure, and the precise coefficients $e_{n,k}$ 
are the main point:

\begin{Proposition}\label{P:TLDims}
We have $\dimk(\simple[k])=
\sum_{r=0}^{c(k)}e_{n-2r+1,k+1}\left(\tfrac{n-2c(k)+1}{n-c(k)+1}\binom{n}{c(k)}\right)$. In particular, for $k\in\{0,1,n\}$ we have $\dimk(\simple[k])=1$.
\end{Proposition}

\begin{proof}
For the Temperley--Lieb algebra $\tlalg{1}$ these dimensions were computed in \cite[Corollary 9.3]{Sp-modular-tl}. 
These computations 
use the $\K$-linear cell structure of $\tlalg{1}$ given by it being 
a cellular algebra. These turn out to be the same calculations 
as for the cell structure of the Temperley--Lieb monoid $\tlmon$ 
discussed in \autoref{T:CellsDimensions} and the results 
in \cite[Corollary 9.3]{Sp-modular-tl} work 
thus for $\tlmon$ without change.
\end{proof}

\begin{Example}\label{E:TLDim}
It is easy to feed the above into a machine. Below 
we list the first few dimensions of the simple $\tlmon$-representations $\simple[k]$ 
for $\cchar=0$ (first table), and $\cchar=2$ (second table).
Here $0\leq n\leq 16$ is indexing the rows and $0\leq k\leq 16$ the columns.
\begin{gather*}
\scalebox{0.6}{\begin{tikzpicture}[anchorbase, ampersand replacement=\&]
\matrix (A) [matrix of nodes,every node/.style={anchor=base,text depth=0.5ex,text height=1ex,text width=2em}]
{
n\text{\textbackslash}k\& 0 \& 1 \& 2 \& 3 \& 4 \& 5 \& 6 \& 7 \& 8 \& 9 \& 10 \& 11 \& 12 \& 13 \& 14 \& 15 \& 16
\\
0 \& 1 \& \& \& \& \& \& \& \& \& \& \& \& \& \& \& \&
\\
1 \& \& 1 \& \& \& \& \& \& \& \& \& \& \& \& \& \& \&
\\
2 \& 1 \& \& 1 \& \& \& \& \& \& \& \& \& \& \& \& \& \&
\\
3 \& \& 1 \& \& 1 \& \& \& \& \& \& \& \& \& \& \& \& \&
\\
4 \& 1 \& \& 3 \& \& 1 \& \& \& \& \& \& \& \& \& \& \& \&
\\
5 \& \& 1 \& \& 4 \& \& 1 \& \& \& \& \& \& \& \& \& \& \&
\\
6 \& 1 \& \& 9 \& \& 4 \& \& 1 \& \& \& \& \& \& \& \& \& \&	
\\
7 \& \& 1 \& \& 13 \& \& 6 \& \& 1 \& \& \& \& \& \& \& \& \&
\\
8 \& 1 \& \& 28 \& \& 13 \& \& 7 \& \& 1 \& \& \& \& \& \& \& \&
\\
9 \& \& 1 \& \& 41 \& \& 27 \& \& 7 \& \& 1 \& \& \& \& \& \& \&
\\
10 \& 1 \& \& 90 \& \& 41 \& \& 34 \& \& 9 \& \& 1 \& \& \& \& \& \&
\\
11 \& \& 1 \& \& 131 \& \& 110 \& \& 34 \& \& 10 \& \& 1 \& \& \& \& \&
\\
12 \& 1 \& \& 297 \& \& 131 \& \& 144 \& \& 54 \& \& 10 \& \& 1 \& \& \& \&
\\
13 \& \& 1 \& \& 428 \& \& 429 \& \& 144 \& \& 64 \& \& 12 \& \& 1 \& \& \&
\\
14 \& 1 \& \& 1001 \& \& 428 \& \& 573 \& \& 273 \& \& 64 \& \& 13 \& \& 1 \& \&	
\\
15 \& \& 1 \& \& 1429 \& \& 1638 \& \& 573 \& \& 337 \& \& 90 \& \& 13 \& \& 1 \&
\\
16 \& 1 \& \& 3432 \& \& 1429 \& \& 2211 \& \& 1260 \& \& 337 \& \& 103 \& \& 15 \& \& 1
\\
};
\node at (8,3.5) {$\cchar=0$};
\draw [very thick] ($(A-1-1.south west)+(0,0.05)$) to ($(A-1-18.south east)+(0,0.05)$);
\draw [very thick] ($(A-1-1.north east)+(-0.1,0)$) to ($(A-18-1.south east)+(-0.1,0)$);
\draw [very thick,densely dotted] ($(A-1-2.north east)+(0.55,0)$) to ($(A-18-2.south east)+(0.55,0)$);
\draw [very thick,->] (1.5,-0.15)node[right]{$\dimk(\simple[8])$ for $\tlmon[8]$} to (0.5,-0.15);
\draw [very thick,->] (5.7,-1.5)node[right]{$\dimk(\simple[11])$ for $\tlmon[15]$} to (3.5,-3.5);
\draw [very thick,->] (1,3) to node[above,yshift=0.05cm]{Order of cells} (0.1,3);
\end{tikzpicture}}
\quad.
\end{gather*}
\begin{gather*}
\scalebox{0.6}{\begin{tikzpicture}[anchorbase, ampersand replacement=\&]
\matrix (A) [matrix of nodes,every node/.style={anchor=base,text depth=0.5ex,text height=1ex,text width=2em}]
{
n\text{\textbackslash}k\& 0 \& 1 \& 2 \& 3 \& 4 \& 5 \& 6 \& 7 \& 8 \& 9 \& 10 \& 11 \& 12 \& 13 \& 14 \& 15 \& 16
\\
0 \& 1 \& \& \& \& \& \& \& \& \& \& \& \& \& \& \& \&
\\
1 \& \& 1 \& \& \& \& \& \& \& \& \& \& \& \& \& \& \&
\\
2 \& 1 \& \& 1 \& \& \& \& \& \& \& \& \& \& \& \& \& \&
\\
3 \& \& 1 \& \& 1 \& \& \& \& \& \& \& \& \& \& \& \& \&
\\
4 \& 1 \& \& 3 \& \& 1 \& \& \& \& \& \& \& \& \& \& \& \&
\\
5 \& \& 1 \& \& 4 \& \& 1 \& \& \& \& \& \& \& \& \& \& \&
\\
6 \& 1 \& \& 9 \& \& 4 \& \& 1 \& \& \& \& \& \& \& \& \& \&	
\\
7 \& \& 1 \& \& 13 \& \& 6 \& \& 1 \& \& \& \& \& \& \& \& \&
\\
8 \& 1 \& \& 27 \& \& 13 \& \& 7 \& \& 1 \& \& \& \& \& \& \& \&
\\
9 \& \& 1 \& \& 40 \& \& 27 \& \& 7 \& \& 1 \& \& \& \& \& \& \&
\\
10 \& 1 \& \& 81 \& \& 40 \& \& 34 \& \& 9 \& \& 1 \& \& \& \& \& \&
\\
11 \& \& 1 \& \& 121 \& \& 110 \& \& 34 \& \& 10 \& \& 1 \& \& \& \& \&
\\
12 \& 1 \& \& 243 \& \& 121 \& \& 144 \& \& 54 \& \& 10 \& \& 1 \& \& \& \&
\\
13 \& \& 1 \& \& 364 \& \& 429 \& \& 144 \& \& 64 \& \& 12 \& \& 1 \& \& \&
\\
14 \& 1 \& \& 729 \& \& 364 \& \& 573 \& \& 272 \& \& 64 \& \& 13 \& \& 1 \& \&	
\\
15 \& \& 1 \& \& 1093 \& \& 1638 \& \& 573 \& \& 336 \& \& 90 \& \& 13 \& \& 1 \&
\\
16 \& 1 \& \& 2187 \& \& 1093 \& \& 2211 \& \& 1245 \& \& 336 \& \& 103 \& \& 15 \& \& 1
\\
};
\node at (8,3.5) {$\cchar=2$};
\draw [very thick] ($(A-1-1.south west)+(0,0.05)$) to ($(A-1-18.south east)+(0,0.05)$);
\draw [very thick] ($(A-1-1.north east)+(-0.1,0)$) to ($(A-18-1.south east)+(-0.1,0)$);
\draw [very thick,densely dotted] ($(A-1-2.north east)+(0.55,0)$) to ($(A-18-2.south east)+(0.55,0)$);
\draw [very thick,->] (1.5,-0.15)node[right]{$\dimk(\simple[8])$ for $\tlmon[8]$} to (0.5,-0.15);
\draw [very thick,->] (5.7,-1.5)node[right]{$\dimk(\simple[11])$ for $\tlmon[15]$} to (3.5,-3.5);
\draw [very thick,->] (1,3) to node[above,yshift=0.05cm]{Order of cells} (0.1,3);
\end{tikzpicture}}
\quad.
\end{gather*}
These tables also appear in \cite{An-simple-tl}. 
Note that the representation $\simple[0]$ 
for even $n$ and $\simple[1]$ for odd $n$, separated 
by a dotted line, are always of dimension one. 
This is a special coincidence of the involved combinatorics and was observed 
from a very different direction in 
\cite[Proposition 4.5]{SuTuWeZh-mixed-tilting}.
\end{Example}

Recall that $c(k)$ denotes the number of caps respectively cups in $\jcell_{k}$.
The following lower bound for the dimensions 
($k\notin\{0,1\}$ is covered in \autoref{P:TLDims}):

\begin{Proposition}\label{P:TLLowerBound}
Let $\cchar=0$.
For $k\notin\{0,1\}$ we have
\begin{gather*}
\dimk(\simple[k])\geq
\frac{1}{(n-c(k)+1)(n-c(k)+2)}\binom{n}{c(k)}
.
\end{gather*}
\end{Proposition}

\begin{proof}
See \cite[Propositions 9.4 and 9.5]{Sp-modular-tl}.
\end{proof}

\begin{Example}\label{E:TLLowerBound}
The dimensions of the simple $\tlmon[24]$-representations over $\Q$ and their lower bounds are given by the tuples
\begin{gather*}
\scalebox{0.95}{$\text{dim}\colon\left(1,534888,208011,445741,389367,126292,85216,31878,6876,1726,252,22,1\right)$},
\\
\scalebox{0.95}{$\text{lower bound}\colon\left(14858,11886,8171,4807,2403,1012,354,101,23,4,0.5,0.04,\tfrac{1}{650}\right)$}.
\end{gather*}
Here $\simple[0]$ correspond to the leftmost entry and then $k$ increases in steps of two from left to right.
Note that the lower bound does not work for $k=0$.
\end{Example}


\subsection{Truncating the Temperley--Lieb monoid}\label{SS:TLSubmonoid}

Recall that we need a regularity condition 
to ensure that taking cell subquotients works as expected, 
{\cf} \autoref{L:CellsAdmissible}.
We first establish:

\begin{Lemma}\label{L:TLAdmissible}
The monoid $\tlmon$ is regular.
\end{Lemma}

\begin{proof}
We check that $\tlmon$ satisfies condition (a) in \autoref{L:CellsAdmissible}. 
Take $a$ and $b$ 
with $k$ through strands, both in the same left cell within $\jcell_{k}$. 
Thus, $a$ and $b$ have the same bottom 
half $\beta_{a}=\beta_{b}$ but the top halves
$\gamma_{a}\neq\gamma_{b}$ can be different. We can now use 
$c=\gamma_{a}\circ(\gamma_{b})^{\ast}$ which implies that 
$a=cb$, as required. The picture is:
\begin{gather*}
\begin{tikzpicture}[anchorbase]
\draw[usual] (0,0) to[out=90,in=180] (0.25,0.2) to[out=0,in=90] (0.5,0);
\draw[usual] (0,0.5) to[out=270,in=180] (0.25,0.3) to[out=0,in=270] (0.5,0.5);
\draw[usual] (1,0) to (1,0.5);
\draw[usual] (1.5,0) to (1.5,0.5);
\draw[usual] (0,0.5) to[out=90,in=180] (0.25,0.7) to[out=0,in=90] (0.5,0.5);
\draw[usual] (1,1) to[out=270,in=180] (1.25,0.8) to[out=0,in=270] (1.5,1);
\draw[usual] (1,0.5) to (0,1);
\draw[usual] (1.5,0.5) to (0.5,1);
\end{tikzpicture}
=
\begin{tikzpicture}[anchorbase,yscale=-1]
\draw[usual] (1,0) to[out=90,in=180] (1.25,0.2) to[out=0,in=90] (1.5,0);
\draw[usual] (0,0.5) to[out=270,in=180] (0.25,0.3) to[out=0,in=270] (0.5,0.5);
\draw[usual] (0,0) to (1,0.5);
\draw[usual] (0.5,0) to (1.5,0.5);
\end{tikzpicture}
\;,
\end{gather*}
which is a calculation in $\jcell_{2}$ as in \autoref{Eq:TLCells}.

Alternatively, from \autoref{L:CellsAdmissible} we get that this lemma follows from Neumann regularity in monoid theory, and of course 
this is well-known for the Temperley--Lieb monoid, see {\eg} 
\cite{LaFiGe-ideal-kauffman}.
\end{proof}

Motivated by \autoref{E:TLtruncated} below we define:

\begin{Definition}\label{D:TLTruncated}
Define the \emph{$k$th truncated Temperley--Lieb monoid} by
\begin{gather*}
\tltru{k}=(\tlmon[n])_{\geq\jcell_{k}}.
\end{gather*}
\end{Definition}

This is the cell submonoid, see \autoref{SS:CellsSubmon}.
In words, $\tltru{k}$ consist of all crossingless matchings with fewer than $k$ through strands, together with an identity element.
Recall that, by \autoref{P:CellsSub} and the 
discussion in \autoref{SS:TLDef}, we know the 
simple $\tltru{k}$-representations and their dimensions.

\begin{Example}\label{E:TLtruncated}
Let us come back to \autoref{Eq:CellsDimTL}. Looking at 
the graphs of the dimensions 
and the semisimple dimensions of the simple $\tlmon[24]$-representations
\begin{gather}\label{Eq:TlDimTL}
\begin{gathered}
\scalebox{0.9}{$\begin{tikzpicture}[anchorbase]
\node at (0,0) {\includegraphics[height=4.4cm]{dimensions1.pdf}};
\node at (2.25,2) {$\tlmon[24]/\Q$};
\node at (2.25,1.5) {cells increase $\leftarrow$};
\draw[->] (-2.4,2.2) node[above,xshift=0.3cm]{y-axis: dim} to (-2.4,2.1) to (-2.8,2.1);
\draw[->] (2.75,-1.5) node[above,yshift=0.75cm]{x-axis:}node[above,yshift=0.3cm]{\# through}node[above,yshift=-0.15cm]{strands/2} to (2.75,-1.8);
\draw[very thick,densely dotted] (-0.7,-2) to (-0.7,1.5) node[above]{$k\approx 2\sqrt{24}$};
\end{tikzpicture}
,
\begin{tikzpicture}[anchorbase]
\node at (0,0) {\includegraphics[height=4.4cm]{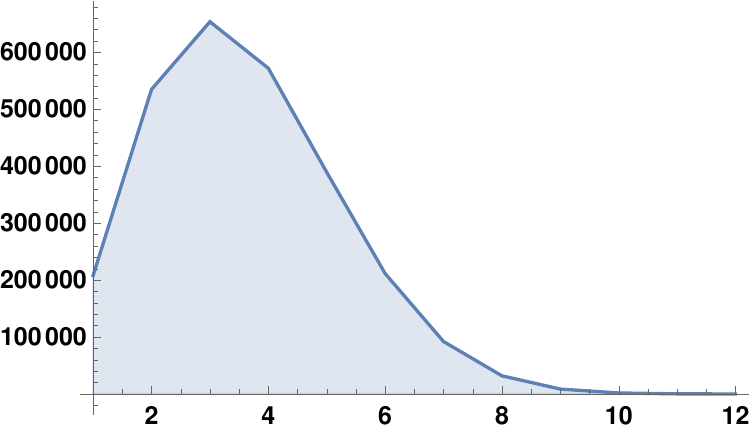}};
\node at (2.25,2) {$\tlmon[24]/\K$};
\node at (2.25,1.5) {cells increase $\leftarrow$};
\draw[->] (-2.4,2.2) node[above,xshift=0.3cm]{y-axis: ssdim} to (-2.4,2.1) to (-2.8,2.1);
\draw[->] (2.75,-1.5) node[above,yshift=0.75cm]{x-axis:}node[above,yshift=0.3cm]{\# through}node[above,yshift=-0.15cm]{strands/2} to (2.75,-1.8);
\draw[very thick,densely dotted] (-0.7,-2) to (-0.7,1.5) node[above,xshift=0.25cm]{$k\approx 2\sqrt{24}$};
\end{tikzpicture}$}
,
\\
\scalebox{0.9}{$\text{dim}\colon\left(1,534888,208011,445741,389367|126292,85216,31878,6876,1726,252,22,1\right)$},
\\
\scalebox{0.9}{$\text{ssdim}\colon\left(208012,534888,653752,572033,389367| 211508,92092,31878,8602,1748,252,23,1\right)$},
\end{gathered}
\end{gather}
it seems preferable to cut these graphs roughly 
at $k\approx2\sqrt{24}$ or at even lower values, as 
illustrated above. The submonoid $\tltru[24]{k}$ for this specific 
value of $k$ now does not have too small representations anymore and is still rich 
enough as a monoid. 
Note that the one-dimensional simple $\tltru[24]{k}$-representation for $\jcell_{0}$ is $\oneb$, so we do not need to get rid of it.
\end{Example}

Our main statement about the Temperley--Lieb case is 
a bound for the representation gap of $\tltru{k}$, but 
before we can prove it we need to discuss extensions.


\subsection{Trivial extensions in Temperley--Lieb monoids}\label{SS:TLTrivExtTL} 

Our next goal is to show that $\tlmon$ and $\tltru{k}$ have no extensions between $\onebt$ (under some minor restrictions on $n$ and $k$). Let $\xmon$ be either $\tlmon$ or $\tltru{k}$, and recall the notions of 
left-connected, right-connected, null-connected and well-connected
from \autoref{SS:RepGapRepGap}.

\begin{Lemma}\label{L:TLNullRounded}
The monoid $\xmon$ is null-connected.
\end{Lemma}

\begin{proof}
Note first that for each of these monoids the group $\group$ of invertible elements is trivial. For $a\in\xmon\setminus\group$ the decomposition $a=\gamma\circ\id_{m}\circ\beta=a_{1}a_{2}^{\ast}$ from \autoref{L:TLFactor} 
then implies that these monoids are null-connected since $a_{2}^{\ast}=a_{2}^{\ast}a_{2}a_{2}^{\ast}$ and 
$a=aa_{2}a_{2}^{\ast}$ is then a product of $a$ and $a_{2}^{\ast}a_{2}\in\xmon\setminus\group$. 
\end{proof}

Before we can prove the main statement of this section we need 
some terminology.

\begin{Remark}\label{R:TLRounded}
The reader might recognize the definitions 
below from the theory of Temperley--Lieb cells 
(or the many other occasions where this theory has appeared in disguise).
That is no coincidence as the notions of being left or right-connected 
are closely related to left and right cells.
\end{Remark}

Recall that a diagram $a\in\nobottom{m}{n}$ 
consists of $m$ through strands and $\frac{n-m}{2}$ cups. 
The through strands connect top and bottom endpoints in $a$, while cups connects top endpoints in pairs. 

For the following notion we naively compose diagrams, meaning 
that we do not remove internal circles.
We say that $a,b\in\nobottom{m}{n}$ are 
in a \emph{vertical position} if the 
diagram $b^{\ast}a$ is isotopic to $\id_{m}$, the identity 
diagram on $m$ strands. 
Elements $a,b\in\nobottom{m}{n}$ are said to be in a \emph{weakly vertical position} if $b^{\ast}a$ is isotopic to $\id_{m}$ together with potential internal circles.

\begin{Example}\label{E:TLVertical}
Consider $a,b\in\nobottom{2}{6}$ given by
\begin{gather*}
a=
\begin{tikzpicture}[anchorbase]
\draw[usual] (0.5,1) to[out=270,in=180] (0.75,0.75) to[out=0,in=270] (1,1);
\draw[usual] (-1,1) to[out=270,in=180] (-0.75,0.75) to[out=0,in=270] (-0.5,1);
\draw[usual] (0,0.5) to (0,1);
\draw[usual] (1.5,0.5) to (1.5,1);
\end{tikzpicture}
\;,\quad
b=
\begin{tikzpicture}[anchorbase,xscale=-1]
\draw[usual] (0.5,1) to[out=270,in=180] (0.75,0.75) to[out=0,in=270] (1,1);
\draw[usual] (-1,1) to[out=270,in=180] (-0.75,0.75) to[out=0,in=270] (-0.5,1);
\draw[usual] (0,0.5) to (0,1);
\draw[usual] (1.5,0.5) to (1.5,1);
\end{tikzpicture}
\;,\quad
b^{\ast}a
=
\begin{tikzpicture}[anchorbase]
\draw[usual] (0.5,1) to[out=270,in=180] (0.75,0.75) to[out=0,in=270] (1,1);
\draw[usual] (-1,1) to[out=270,in=180] (-0.75,0.75) to[out=0,in=270] (-0.5,1);
\draw[usual] (0,0.5) to (0,1);
\draw[usual] (1.5,0.5) to (1.5,1);
\draw[usual] (0,1) to[out=90,in=0] (-0.25,1.25) to[out=180,in=90] (-0.5,1);
\draw[usual] (1.5,1) to[out=90,in=0] (1.25,1.25) to[out=180,in=90] (1,1);
\draw[usual] (0.5,1) to (0.5,1.5);
\draw[usual] (-1,1) to (-1,1.5);
\end{tikzpicture}
\;.
\end{gather*}
Then $a$ and $b$ are in vertical position, as illustrated above. But neither $a$ and $a$ nor $b$ and $b$ are. The latter are only in weakly vertical position.
\end{Example}

Denote by $\vertical{n}{m}\subset\nobottom{m}{n}\times\nobottom{m}{n}$ the set of pairs of diagrams in a vertical position, and write $(a,b)\in\vertical{n}{m}$. This relation on diagrams is symmetric. 

Denote by $\wvertical{m}{n}\subset\nobottom{m}{n}\times\nobottom{m}{n}$ the set of pairs of diagrams in weakly vertical position, and write $(a,b)\in\wvertical{m}{n}$. Note that $(a,a)\in\wvertical{m}{n}$ for any $a\in\nobottom{m}{n}$.
Again, this relation on diagrams is symmetric.

If $(a_{2},b_{1})\in\wvertical{m}{n}$, then $a_{1}a_{2}^{\ast}b_{1}b_{2}^{\ast}=a_{1}(a_{2}^{\ast}b_{1})b_{2}^{\ast}=a_{1}b_{2}^{\ast}$. That is, inserting $a_{2}^{\ast}b_{1}$ in the middle of $a_{1}b_{2}^{\ast}$ does not change the latter.  

\begin{Definition}\label{D:TLGraph1}
Let $\tlgraph{m}{n}$ denote the unoriented graph
with vertex set $\nobottom{m}{n}$ and edges between $a$ and $b$ for all $(a,b)\in\vertical{m}{n}$. 
\end{Definition}

Note that $\tlgraph{m}{n}$ is 
nonempty if and only if $n\geq m$ and $n+m$ is even.

\begin{Lemma}\label{L:TLGraphConn}
The graph $\tlgraph{m}{n}$ is connected if $m>0$. 
\end{Lemma}

\begin{proof}
\textit{Case $m=1$.}
In this case the lemma can be proved by induction on $n$, 
by showing that any diagram $a\in\nobottom{1}{n}$ is 
connected by a path in $\tlgraph{1}{n}$ to the
diagram 
\begin{gather*}
\begin{tikzpicture}[anchorbase,xscale=-1]
\draw[usual] (0,0.5) to[out=270,in=180] (0.25,0.25) to[out=0,in=270] (0.5,0.5);
\node at (0.75,0.35) {$\dots$};
\draw[usual] (1,0.5) to[out=270,in=180] (1.25,0.25) to[out=0,in=270] (1.5,0.5);
\draw[usual] (2,0) to (2,0.5);
\end{tikzpicture}
\end{gather*}
with the through strand on the 
far left and $\tfrac{n-1}{2}$ unnested cups.
\medskip

\textit{General case.} Consider a diagram $a\in\nobottom{m}{n}$. 
Each through strand $c$ of $a$ may be surrounded by a cluster of cups on either side. The first case allows to bring each such cluster together with $c$ to a standard form as above (through strands followed by a sequence of unnested cups) via paths in suitable graphs $\tlgraph{1}{k}$, utilizing only one through strand $c$. Doing this transformation with each through strands in $a$ and moving all through strands all the way to the left transforms $a$ to a standard form of $m$ parallel vertical strands on the left followed by unnested $\frac{n-m}{2}$ cups. This shows that $\tlgraph{m}{n}$ is connected. 
\end{proof}

A cup is called \emph{outer} if it is not separated from the bottom of the diagram by any cup. 
A pair $(a,b)\in\nobottom{m}{n}\times\nobottom{m}{n}$ is called a \emph{flip pair} if $b$ is obtained from $a$ by converting an outer cup $c$ into a pair of through strands while simultaneously closing up a pair $p$ of adjacent through strands in $a$ into a cup. Note that $c$ must not be located between the two strands in $p$, and that the flip pair relation is symmetric.

\begin{Example}\label{E:TLFlippedPair}
\leavevmode	

\begin{enumerate}

\item In the element of $\nobottom{3}{10}$
\begin{gather*}
\begin{tikzpicture}[anchorbase]
\draw[usual] (0,0) to (0,0.5)node[above,yshift=-0.05cm]{1};
\draw[usual] (0.5,0.5)node[above,yshift=-0.05cm]{2} to[out=270,in=180] (1.25,0.1) to[out=0,in=270] (2,0.5)node[above,yshift=-0.05cm]{5};
\draw[usual] (1,0.5)node[above,yshift=-0.05cm]{3} to[out=270,in=180] (1.25,0.3) to[out=0,in=270] (1.5,0.5)node[above,yshift=-0.05cm]{4};
\draw[usual] (2.5,0) to (2.5,0.5)node[above,yshift=-0.05cm]{6};
\draw[usual] (3,0) to (3,0.5)node[above,yshift=-0.05cm]{7};
\draw[usual] (3.5,0.5)node[above,yshift=-0.05cm]{8} to[out=270,in=180] (3.75,0.3) to[out=0,in=270] (4,0.5)node[above,yshift=-0.05cm]{9};
\draw[usual] (4.5,0) to (4.5,0.5)node[above,yshift=-0.05cm]{10};
\end{tikzpicture}
\;,
\end{gather*}
the cups $(2,5)$ and $(8,9)$ are outer. The cups $(2,5)$ and $(3,4)$ are nested, with $(2,5)$ an outer nested cup.

\item The following is a flip pair:
\begin{gather*}
\begin{tikzpicture}[anchorbase]
\draw[usual] (0,0) to (0,0.5);
\draw[usual] (0.5,0.5) to[out=270,in=180] (0.75,0.3) to[out=0,in=270] (1,0.5);
\draw[usual] (1.5,0) to (1.5,0.5);
\draw[usual] (2,0) to (2,0.5);
\draw[usual] (2.5,0.5) to[out=270,in=180] (2.75,0.3) to[out=0,in=270] (3,0.5);
\draw[usual] (3.5,0) to (3.5,0.5);
\draw[ultra thick,spinach,dotted] (0.75,0) to (0.75,0.3);
\draw[ultra thick,spinach,dotted] (2,0.2) to (3.5,0.2);
\end{tikzpicture}
\ \xleftrightarrow{\text{flip pair}} \ 
\begin{tikzpicture}[anchorbase]
\draw[usual] (0,0) to (0,0.5);
\draw[usual] (0.5,0) to (0.5,0.5);
\draw[usual] (1,0) to (1,0.5);
\draw[usual] (1.5,0) to (1.5,0.5);
\draw[usual] (2,0.5) to[out=270,in=180] (2.75,0.1) to[out=0,in=270] (3.5,0.5);
\draw[usual] (2.5,0.5) to[out=270,in=180] (2.75,0.3) to[out=0,in=270] (3,0.5);
\end{tikzpicture}
\;.
\end{gather*}
We have indicated where we apply operations on the left-hand diagram.

\end{enumerate}
The reader might want to think of a flip pair as two diagrams 
related by opposite saddles moves.
\end{Example}

\begin{Definition}\label{D:TLGraph2}
Let $\tlgraphtwo{m}{n}$ denote the unoriented graph
with vertex set $\nobottom{m}{n}$ and edges between $a$ and $b$ for all flipped pairs $(a,b)$. 
\end{Definition}

Note that the graph $\tlgraphtwo{2}{n}$ is not connected for even $n\ge 4$, 
since a diagram with one through strand on the far left and one on the far right is not in any flip pair.

\begin{Lemma}\label{L:TLGraphConn2}
The graph $\tlgraphtwo{3}{n}$ is connected.
\end{Lemma}

\begin{proof}
We view the empty graph (with no vertices) as connected, which is the case when $n$ is even. For odd $n=3+2k$ the proof is by induction on $k$. 

\textit{Case $k=1$.} In this case the graph has the form 
\begin{gather*}
\begin{tikzpicture}[anchorbase]
\draw[usual] (0,0) to (0,0.5);
\draw[usual] (0.5,0.5) to[out=270,in=180] (0.75,0.3) to[out=0,in=270] (1,0.5);
\draw[usual] (1.5,0) to (1.5,0.5);
\draw[usual] (2,0) to (2,0.5);
\end{tikzpicture}
\ \xleftrightarrow{\text{flip}} \ 
\begin{tikzpicture}[anchorbase]
\draw[usual] (0,0) to (0,0.5);
\draw[usual] (0.5,0) to (0.5,0.5);
\draw[usual] (1,0) to (1,0.5);
\draw[usual] (1.5,0.5) to[out=270,in=180] (1.75,0.3) to[out=0,in=270] (2,0.5);
\end{tikzpicture}
\ \xleftrightarrow{\text{flip}} \ 
\begin{tikzpicture}[anchorbase,xscale=-1]
\draw[usual] (0,0) to (0,0.5);
\draw[usual] (0.5,0) to (0.5,0.5);
\draw[usual] (1,0) to (1,0.5);
\draw[usual] (1.5,0.5) to[out=270,in=180] (1.75,0.3) to[out=0,in=270] (2,0.5);
\end{tikzpicture}
\ \xleftrightarrow{\text{flip}} \ 
\begin{tikzpicture}[anchorbase,xscale=-1]
\draw[usual] (0,0) to (0,0.5);
\draw[usual] (0.5,0.5) to[out=270,in=180] (0.75,0.3) to[out=0,in=270] (1,0.5);
\draw[usual] (1.5,0) to (1.5,0.5);
\draw[usual] (2,0) to (2,0.5);
\end{tikzpicture}
\;, 
\end{gather*}
thus is connected. 
\medskip

\textit{Case $k>1$.} Represent $a\in\nobottom{3}{3+2k}$ as a composition $a=cb$ of a diagram $c\in\nobottom{1+2k}{3+2k}$ with a single cup and $b\in\nobottom{3}{1+2k}$, {\eg}: 
\begin{gather*}
a=
\begin{tikzpicture}[anchorbase]
\draw[usual] (0,0) to (0,0.5);
\draw[usual] (0.5,0.5) to[out=270,in=180] (1.25,0.1) to[out=0,in=270] (2,0.5);
\draw[usual] (1,0.5) to[out=270,in=180] (1.25,0.3) to[out=0,in=270] (1.5,0.5);
\draw[usual] (2.5,0) to (2.5,0.5);
\draw[usual] (3,0.5) to[out=270,in=180] (3.25,0.3) to[out=0,in=270] (3.5,0.5);
\draw[usual] (4,0) to (4,0.5);
\end{tikzpicture}
\;=\; 
\begin{tikzpicture}[anchorbase]
\draw[usual] (0,0) to (0,0.5);
\draw[usual] (0.5,0.5) to[out=270,in=180] (1.25,0.1) to[out=0,in=270] (2,0.5);
\draw[usual] (2.5,0) to (2.5,0.5);
\draw[usual] (3,0.5) to[out=270,in=180] (3.25,0.3) to[out=0,in=270] (3.5,0.5);
\draw[usual] (4,0) to (4,0.5);
\draw[ultra thick,spinach,densely dotted] (-0.5,0.5) to (4.5,0.5)node[above]{$c$}node[below]{$b$};
\draw[usual] (0,0.5) to (0,1);
\draw[usual] (0.5,0.5) to (0.5,1);
\draw[usual] (1,1) to[out=270,in=180] (1.25,0.8) to[out=0,in=270] (1.5,1);
\draw[usual] (2,0.5) to (2,1);
\draw[usual] (2.5,0.5) to (2.5,1);
\draw[usual] (3,0.5) to (3,1);
\draw[usual] (3.5,0.5) to (3.5,1);
\draw[usual] (4,0.5) to (4,1);
\end{tikzpicture}
\;.
\end{gather*}
By induction on $k$, the diagram $b$ is connected in the graph $\tlgraphtwo{3}{2k+1}$ to the diagram $b_{k-1}$, called standard, of three through strands on the far left and by $k-1$ unnested cups on the right. For example
\begin{gather}\label{Eq:TLProofConn1}
b=
\begin{tikzpicture}[anchorbase,xscale=-1]
\draw[usual] (0,0) to (0,0.5);
\draw[usual] (0.5,0.5) to[out=270,in=180] (0.75,0.3) to[out=0,in=270] (1,0.5);
\draw[usual] (1.5,0) to (1.5,0.5);
\draw[usual] (2,0.5) to[out=270,in=180] (2.25,0.3) to[out=0,in=270] (2.5,0.5);
\draw[usual] (3,0) to (3,0.5);
\end{tikzpicture}
\xleftrightarrow{\text{flip}}
b_{k-1}=
\begin{tikzpicture}[anchorbase,xscale=1]
\draw[usual] (0,0) to (0,0.5);
\draw[usual] (0.5,0) to (0.5,0.5);
\draw[usual] (1,0) to (1,0.5);
\draw[usual] (1.5,0.5) to[out=270,in=180] (1.75,0.3) to[out=0,in=270] (2,0.5);
\draw[usual] (2.5,0.5) to[out=270,in=180] (2.75,0.3) to[out=0,in=270] (3,0.5);
\end{tikzpicture}
\Rightarrow
cb_{k-1}=
\begin{tikzpicture}[anchorbase,xscale=1]
\draw[usual] (0,0) to (0,0.5);
\draw[usual] (0.5,0) to (0.5,0.5);
\draw[usual] (1,0.5) to[out=270,in=180] (1.25,0.3) to[out=0,in=270] (1.5,0.5);
\draw[usual] (2,0) to (2,0.5);
\draw[usual] (2.5,0.5) to[out=270,in=180] (2.75,0.3) to[out=0,in=270] (3,0.5);
\draw[usual] (3.5,0.5) to[out=270,in=180] (3.75,0.3) to[out=0,in=270] (4,0.5);
\end{tikzpicture}
\;.
\end{gather} 
Consequently, in $\tlgraphtwo{3}{2k+3}$ the diagrams $a$ and $cb_{k-1}$ are connected. If the extra cup in $cb_{k-1}$ coming from $c$ is not nested inside the rightmost cup of $b_{k-1}$, as shown in \autoref{Eq:TLProofConn2}, then $cb_{k-1}$ can be represented as a diagram in $\nobottom{3}{2k+1}$ union a cup on the far right, $cb_{k-1}=d\otimes\cup$, with $d$ in very specific form, see \autoref{Eq:TLProofConn1}. By induction, $d$ is connected to the standard diagram $b_{k-1}$ in $\tlgraphtwo{3}{2k+1}$. The union of the latter with a cup on the far right gives the standard diagram in $\tlgraphtwo{3}{2k+3}$, implying that $a$ is connected to the standard diagram $b_{k}\in\tlgraphtwo{3}{2k+3}$. 
\begin{gather}\label{Eq:TLProofConn2}
\begin{aligned}
cb_{k-1}=&\ \ 
\begin{tikzpicture}[anchorbase,xscale=1]
\draw[usual] (0,0) to (0,0.5);
\draw[usual] (0.5,0) to (0.5,0.5);
\draw[usual] (1,0) to (1,0.5);
\draw[usual] (1.5,0.5) to[out=270,in=180] (1.75,0.3) to[out=0,in=270] (2,0.5);
\draw[usual] (2.5,0.5) to[out=270,in=180] (3.25,0.1) to[out=0,in=270] (4,0.5);
\draw[usual] (3,0.5) to[out=270,in=180] (3.25,0.3) to[out=0,in=270] (3.5,0.5);
\end{tikzpicture}
\ 
\xleftrightarrow{\text{flip}}\    
\begin{tikzpicture}[anchorbase,xscale=1]
\draw[usual] (0,0) to (0,0.5);
\draw[usual] (0.5,0.5) to[out=270,in=180] (0.75,0.3) to[out=0,in=270] (1,0.5);
\draw[usual] (1.5,0) to (1.5,0.5);
\draw[usual] (2,0.5) to[out=270,in=180] (2.25,0.3) to[out=0,in=270] (2.5,0.5);
\draw[usual] (3,0) to (3,0.5);
\draw[usual] (3.5,0.5) to[out=270,in=180] (3.75,0.3) to[out=0,in=270] (4,0.5);
\end{tikzpicture}
\\
\xleftrightarrow{\text{flip}}\ \ &
\begin{tikzpicture}[anchorbase,xscale=1]
\draw[usual] (0,0.5) to[out=270,in=180] (0.25,0.3) to[out=0,in=270] (0.5,0.5);
\draw[usual] (1,0.5) to[out=270,in=180] (1.25,0.3) to[out=0,in=270] (1.5,0.5);
\draw[usual] (2,0) to (2,0.5);
\draw[usual] (2.5,0) to (2.5,0.5);
\draw[usual] (3,0) to (3,0.5);
\draw[usual] (3.5,0.5) to[out=270,in=180] (3.75,0.3) to[out=0,in=270] (4,0.5);
\end{tikzpicture}
\ \xleftrightarrow{\text{flip}} \ 
\begin{tikzpicture}[anchorbase,xscale=1]
\draw[usual] (0,0) to (0,0.5);
\draw[usual] (0.5,0) to (0.5,0.5);
\draw[usual] (1,0) to (1,0.5);
\draw[usual] (1.5,0.5) to[out=270,in=180] (1.75,0.3) to[out=0,in=270] (2,0.5);
\draw[usual] (2.5,0.5) to[out=270,in=180] (2.75,0.3) to[out=0,in=270] (3,0.5);
\draw[usual] (3.5,0.5) to[out=270,in=180] (3.75,0.3) to[out=0,in=270] (4,0.5);
\end{tikzpicture}
\;.
\end{aligned}
\end{gather} 
The remaining case is when the cup from $c$ is nested inside the rightmost cup of $b_{k-1}$, see \autoref{Eq:TLProofConn2}. Then a series of transformations along paths in the graph $\tlgraphtwo{3}{2k+3}$, possible by induction, show that $a$ is in the same connected component as the diagram $b_{k}$, concluding the induction step.
\end{proof}

\begin{Lemma}\label{L:TLGraphConnM}
The graph $\tlgraphtwo{m}{n}$ is connected for any $m\geq 3$. 
\end{Lemma}

\begin{proof} 
The proof is by induction on $m$. Case $m=3$ has already been established. Denote by $b_{n,m}$ the diagram with $m$ through strands on the far left followed by $\tfrac{n-m}{2}$ unnested cups to the right. If the leftmost strand of $a$ is a through strand, then the diagram $a$ can be written as a union $a=\vert\otimes a^{\prime}$ of a through strand and a diagram $a^{\prime}$ in $\nobottom{m-1}{n-1}$. By induction, $a^{\prime}$ is connected to the standard diagram $b_{n-1,m-1}$ implying that $a$ is connected to $b_{n,m}$. 

If the leftmost strand of $a$ is a cup, consider the three leftmost through strands of $a$ and form the subdiagram $a_{1}$ that consists of these strands and all cups to the left and in between of these through strands. We can write $a=a_{1}\otimes a_{2}$, with $a_{2}$ the complement $a_{2}$ of $a_{1}$ in $a$. By \autoref{L:TLGraphConn2}, $a_{1}$ is connected to some diagram $b_{r}$ 
with $r$ through strands to the left. Hence, $a=a_{1}\otimes a_{2}$ is connected to $b_{r}\otimes a_{2}$. In the latter diagram the leftmost strand is through, and the previous case allows to use the induction step.
\end{proof}

Recall the relation $\approx_{l}$ given by the closure of the relation 
$ba\approx_{l}a$, where $a,b\in\msg$.

\begin{Lemma}\label{L:TLLeftRoundedFirstStep}
\leavevmode
\begin{enumerate}

\item Suppose $(a,b)\in\nobottom{m}{n}$ is a flip pair, and $m\leq k$. Then $aa^{\ast}\approx_{l}bb^{\ast}$ in $\tltru{k}$. Moreover, 
for $a,b\in\nobottom{m}{n}$, $m\leq k$ we have $aa^{\ast}\approx_{l}bb^{\ast}$ in $\tltru{k}$.

\item For $a,b,c\in\nobottom{k}{n}$ we have $ac^{\ast}\approx_{l}bc^{\ast}$ in $\tltru{k}$, for $k\geq 3$. 

\end{enumerate}
\end{Lemma}

\begin{proof}
\textit{(a).} Suppose the flip is described via an outer cup $c$ and a pair $p$ of adjacent through strands in $a$, as in the definition of a flip. Let $d\in\nobottom{m-2}{n}$ is obtained from $a$ by closing up the pair $p$ into a strand. It is straightforward to check that 
$dd^{\ast}aa^{\ast}=dd^{\ast}=dd^{\ast}bb^{\ast}$,
which implies that $aa^{\ast}\approx_{l}bb^{\ast}$.
The second claim follows then from the first. 
\medskip 

\textit{(b).} Since $\tlgraphtwo{k}{n}$ is connected, we can choose a path $a=a_{1},a_{2},\dots,a_{r}=b$ in it, 
with each $(a_i,a_{i+1})$ an edge. Then $a_{i+1}^{\ast}a_i=\id_{k}$, and 
$a_{i+1}c=a_{i+1}a_{i+1}^{\ast}a_{i}c\approx_{l}a_{i}c$, and 
$ac=a_{0}c\approx_{l}a_{r}c=bc$.
\end{proof}

We are ready to prove that the Temperley--Lieb monoids 
are left-connected. 

\begin{Lemma}\label{L:TLLeftRounded}
We have the following.
\begin{enumerate}

\item The monoid $\tlmon$ is 
left-connected if $n\geq 5$.

\item The monoid $\tltru{k}$ is 
left-connected if $n\geq 5$ and $k\geq 3$. 

\end{enumerate}
\end{Lemma}

\begin{proof} 
Recall the generator-relation presentation from \autoref{Eq:TLPresentation}.	
\medskip

\textit{(a).} It is easy to see that $\tlmon[3]$ has two equivalence classes $\{u_{1},u_{2}u_{1}\}$ and $\{u_{2},u_{1}u_{2}\}$ under $\approx_{l}$, which are the top left cells in \autoref{Eq:TLCells}.
The monoid $\tlmon[4]$ also has two $\approx_{l}$ equivalence classes, represented by $u_{1}$ and $u_{2}$. In general,
since the $u_{i}$ generate $\tlmon$, any $\approx_{l}$ equivalence class is represented by some $u_{i}$. 
For $n>4$, each $u_{i}$ is in the same equivalence class as either $u_{1}$ or $u_{n-1}$. For instance, if $i>2$, $u_{i}$ and $u_{1}$ commute and $u_{i}\approx_{l}u_{1}u_{i}=u_{i}u_{1}\approx_{l}u_{1}$. Finally, $u_{1}\approx_{l}u_{n-1}u_{1}=u_{1}u_{n-1}\approx_{l}u_{n-1}$.
\medskip 

\textit{(b).} We need to show that there is a unique equivalence class under $\approx_{l}$ in $\tltru{k}\setminus\{1\}$. 
First, any element $u$ in the latter set is equivalent under $\approx_{l}$ to an element of width $k$. To see
this, write a minimal length presentation $u=u_{i_{r}}u_{i_{r-1}}\dots u_{i_{1}}$ of $u$ as a product of generators. The element $u$ has width $m\leq k$. Pick the smallest $p$ such that the suffix $v=u_{i_{p}}u_{i_{p-1}}\dots u_{i_{1}}$ of the presentation has width $k$ (this is possible since multiplication of an element by a generator $u_{i}$ either preserves the width or reduces it by one). Then $u=v^{\prime}v$ where $v^{\prime}$ is the product of the remaining terms. Note that $v=v v^{\ast}v$ and $u=v^{\prime}v=(v^{\prime}vv^{\ast})v$. Widths $\omega(v^{\prime}vv^{\ast})\leq k$, $\omega(v)=k$, so that both of these elements are in $\tltru{k}\setminus\{1\}$, and $u=v^{\prime}vv^{\ast}v\approx_{l}v$. We see that $u$ is equivalent to an element of width $k$. 
\medskip

Consequently, it is enough to show that 
$a\approx_{l}b$ for any two $a,b$ of width $k$. Factorize $a=a_{1}a_{2}^{\ast}, b=b_{1}b_{2}^{\ast}$ with $a_{1},a_{2},b_{1},b_{2}\in\nobottom{k}{n}$. From \autoref{L:TLLeftRounded} we have $a_{1}a_{2}^{\ast}\approx_{l} a_{2} a_{2}^{\ast}$ and $b_{1}b_{2}^{\ast}\approx_{l}b_{2}b_{2}^{\ast}$. From the same lemma, $a_{2}a_{2}^{\ast}\approx_{l}b_{2} b_{2}^{\ast}$, so that $a\approx_{l}b$.
\end{proof}

Note that the statement of part (b) of the lemma essentially contains part (a) by taking $k=n$. We have included both parts for clarity. 

\begin{Lemma}\label{L:TLWellRounded}
The monoid $\tlmon$ is well-connected if $n\geq 5$, and 
the monoid $\tltru{k}$ is well-connected if $n\geq 5$ and $k\geq 3$.
\end{Lemma}

\begin{proof} 
This is just the combination of the previous lemmas. 
Note hereby that the diagrammatic antiinvolution ${\placeholder}^{\ast}$ 
implies that the monoids $\tlmon$ and $\tltru{k}$ 
are left-connected if and only if they are right-connected.
\end{proof}

Let $\xmon$ be either $\tlmon$ or $\tltru{k}$ for $k\geq 3$.

\begin{Lemma}\label{L:TLHTrivial}
We have
$\HH^{1}(\xmon,\K)\cong 0$ for all $n\in\N$. 
\end{Lemma}

\begin{proof}
\textit{Case $\xmon=\tlmon$.}
A homomorphism $f\colon\tlmon\to\K$ takes each idempotent $e\in\tlmon$ to $0$. 
From the classical generators-relation presentation of $\tlmon$, see 
\autoref{Eq:TLPresentation}, 
it is clear that every nonidentity element is a product of idempotents, 
so we get $\HH^{1}(\tlmon,\K)\cong 0$.
\medskip

\textit{Case $\xmon=\tltru{k}$.} We now need a different argument. Suppose give a homomorphism $f\colon\tlmon\to\K$. Consider all diagrams of width $k$ in $\xmon$. They have the form $ab^{\ast}$, $a,b\in\nobottom{k}{n}$. Necessarily $f(aa^{\ast})=0$. If $(b,c)$ is an edge in $\tlgraph{k}{n}$ and $d\in\nobottom{k}{n}$, then $ab^{\ast}cd^{\ast}=ad^{\ast}$ and there is a relation
\begin{gather*}
f(ad^{\ast})=f(ab^{\ast})+f(cd^{\ast}). 
\end{gather*}
Choosing a path from $a$ to $d$ in $\tlgraph{k}{n}$ allows to write $f(ad^{\ast})$ as a sum over $f(bc^{\ast})$ where $(b,c)$ is an edge in $\tlgraph{k}{n}$. The relation $bc^{\ast}bb^{\ast}=bb^{\ast}$ implies
\begin{gather*}
f(bb^{\ast})=f(bc^{\ast})+f(bb^{\ast}), 
\end{gather*} 
so that $f(bc^{\ast})=0$ for an edge $(b,c)$. Consequently, $f(ad^{\ast})=0$ for $a,d$ as above, and $f(x)=0$ for any $x$ of width $k$ in $\xmon$. The elements $y$ of $\xmon$ of smaller width are products of elements of width $k$, showing that $f(y)=0$ as well. Thus, $f$ is identically $0$ on $\xmon$.
\medskip

\textit{General approach.}
If $\monoid$ is any finite monoid with trivial $H$-cells, then $\HH^{1}(\xmon,\K)\cong 0$, and even $\HH^{1}(\xmon,A)\cong 0$ 
for any abelian group $A$. 
To see this note that \autoref{T:CellsPeriod} 
and $\hcell(e)\cong\onemon$ imply that $\exists M\in\N$ 
with $x^{M}=x^{M+1}$ for all $x\in\monoid$. Therefore each element 
has trivial image under any $f\colon\monoid\to A$ since $A$ is a group.
\end{proof} 

\begin{Proposition}\label{P:TLSplits}
Let $\module$ be 
an $\xmon$-representation. Assume that $n\geq 5$ and in the truncated case $k\geq 3$. Then any short exact sequence 
\begin{gather*}
0\lra\onebt\lra\module\lra\onebt\lra 0
\end{gather*}
splits.
\end{Proposition} 

\begin{proof}
Note that the group of units $\group$ of $\xmon$ is trivial, so we get $\HH^{1}(\group,\K)\cong 0$.
Combine this with \autoref{L:TLWellRounded} and \autoref{T:RepGapH1Condition}.
\end{proof}


\subsection{Representation gap and faithfulness of the Temperley--Lieb monoid}\label{SS:TLMain}

We are ready to state and prove the main statements 
about the Temperley--Lieb monoid.

Let $\xmon$ be either $\tlmon$ or $\tltru{k}$ for $k\geq 3$.

\begin{Theorem}\label{P:TLRepGap}
Let $n\geq 5$, and
let $m(l)$ be the dimension of the simple 
$\xmon$-representation $\simple[l]$ as in \autoref{P:TLDims}. Then:
\begin{gather*}
\gap{\tlmon}=\min\big\{m(l)|l\notin\{0,1,n\}\big\},
\\
\gap{\tltru{k}}=\min\big\{m(l)|l\notin\{0,1,k+1,k+2,\dots,n\}\big\}.
\end{gather*}
\end{Theorem}

\begin{proof}
By \autoref{T:RepGapH1Condition} and \autoref{P:TLSplits}.
\end{proof}

Recall that $k$ denotes the number of through strands, and 
crossingless matchings with $k$ through strands have $\tfrac{n-k}{2}$ caps and cups. In particular, 
$\tltru{k}$ for $0\leq k\leq 2\sqrt{n}$ has crossingless matchings 
with at most $2\sqrt{n}$ through strands and at least $\tfrac{n-2\sqrt{n}}{2}$ 
(this number is bigger than $\sqrt{n}$ for $n>16$)
caps and cups. Also recall the Bachmann--Landau notation $f\in\Theta(g)$, 
meaning that $f$ is bounded both above and below by $g$ asymptotically.

\begin{Theorem}\label{T:TLIsAGoodExample}
Let $n\geq 5$ and fix $0\leq k\leq 2\sqrt{n}$. 
Let $\cchar=0$, and let $\LL$ be an arbitrary field.
We have the following lower bounds:
\begin{align*}
\gap[\K]{\tltru{k}}&\geq\frac{4}{(n+2\sqrt{n}+2)(n+2\sqrt{n}+4)}\binom{n}{\tfrac{n}{2}-\sqrt{n}}\in\Theta\big(2^{n}n^{-5/2}\big)
,\\
\ssgap[\LL]{\tltru{k}}&\geq\frac{2}{2n}\binom{n}{\floor{\tfrac{n}{2}}}
\in\Theta\big(2^{n}n^{-3/2}\big)
,\\
\faith[\K]{\tltru{k}}&\geq\frac{6}{n+4}\binom{n}{\tfrac{n^{\prime}}{2}-1}
\in\Theta\big(2^{n}n^{-3/2}\big)
.
\end{align*}
where in the final bound $n^{\prime}=n$, if $n$ is even, 
and $n^{\prime}=n-1$, if $n$ is odd.
\end{Theorem}

\begin{proof}
\emph{Representation gap.} We will make use of \autoref{P:TLSplits}. By 
\autoref{T:RepGapH1Condition}, this statement 
ensures that we only need to compute 
dimension bounds for simple $\tltru{k}$-representations.	
The first bound then follows from \autoref{P:TLLowerBound}.
The formula $\tfrac{1}{(n-c(k)+1)(n-c(k)+2)}\binom{n}{c(k)}$ 
has its minimum for $k=\floor{2\sqrt{n}}$. Plotting this $k$ into the formula 
and a bit of algebra autopilot gives the claimed lower bound. The asymptotic 
formula then follows by using that $\tfrac{4}{(n+2\sqrt{n}+2)(n+2\sqrt{n}+4)}$ 
is in $\Theta(\tfrac{1}{n^{2}})$, and using Stirling's approximation for $n!$ 
to get that the binomial is in $\Theta(2^{n}n^{-1/2})$.
\medskip

\emph{Semisimple representation gap.}
The second bound can be seen as follows. 
We need to minimize the formula in \autoref{P:TLSemisimpleDims}
for $0\leq k\leq 2\sqrt{n}$. Observe that the 
function $\tfrac{n-2c(k)+1}{n-c(k)+1}\binom{n}{c(k)}$ in $k$ has precisely 
one peak between $k=0$ and $k=\floor{2\sqrt{n}}$, and is monotone increasing 
respectively decreasing otherwise. So we only need to compare the 
two values for $k=0,1$ and $k=\floor{2\sqrt{n}}$, and it is then easy to see 
that the $k=0,1$ value is smaller. Since $c(0)=\tfrac{n-0}{2}$ 
and $c(1)=\floor{\tfrac{n}{2}}$, the result follows. 
The asymptotic formula follows also from Stirling's approximation for $n!$.
\medskip

\emph{Faithfulness.}
For the final bound we use 
\autoref{L:RepGapFaithfulEmbedding}. 
This lemma says that it suffices to find a lower bound 
for $\tltru{2}$: If $n$ is even, then $\tltru{2}\hookrightarrow\tltru{k}$. 
If $n$ is odd, then we can still use $\tltru{2}$ after adding another strand.
Note that a faithful 
$\tltru{2}$-representation can not be a nontrivial extension of $\onebt$ 
by \autoref{P:TLSplits} and also not a direct sum of $\onebt$. Hence, 
any faithful $\tltru{2}$-representation must contain
$\simple[2]$. The combinatorics 
from \autoref{E:TLDimsChar0} implies that $\dimk(\simple[2])=\ssdimk(\simple[2])$, 
so the claimed formula follows from \autoref{P:TLSemisimpleDims}. 
The asymptotic formula can be verified in the same way as for (a) and (b).

Let us note that alternatively to direct computations for the asymptotic formulas, the reader can also input the above 
bounds into a computer algebra system such as Mathematica and ask the computer to algebraically 
manipulate the symbols.
\end{proof}

\begin{Example}\label{E:TLMainTheorem}
The lower bounds in \autoref{T:TLIsAGoodExample} are far from being optimal. But they still grow very fast. Here are their plots:
\begin{gather*}
\scalebox{0.9}{$\begin{tikzpicture}[anchorbase]
\node at (0,0) {\includegraphics[height=4.4cm]{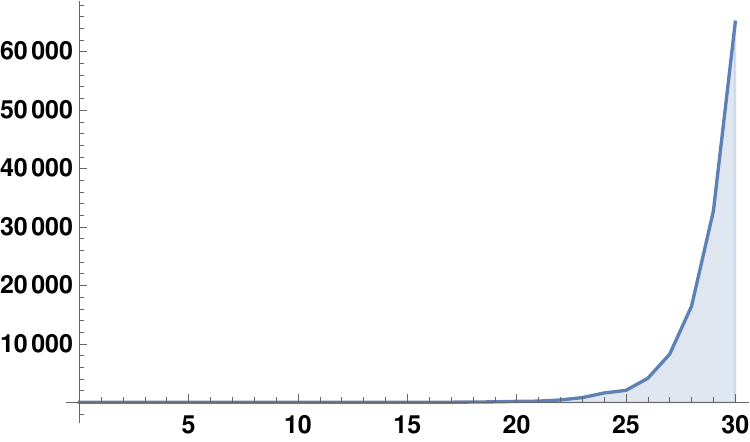}};
\node at (1.5,2) {$\tlmon[n]/\Q$};
\node at (1.5,1.5) {$n$ increases $\rightarrow$};
\draw[->] (-2.4,2.2) node[above,xshift=0.3cm]{y-axis: lower bound gap} to (-2.4,2.1) to (-2.8,2.1);
\draw[->] (1.5,-1.5)node[above]{x-axis: n} to (1.5,-1.8);
\end{tikzpicture}$}
,
\scalebox{0.9}{$\begin{tikzpicture}[anchorbase]
\node at (0,0) {\includegraphics[height=4.4cm]{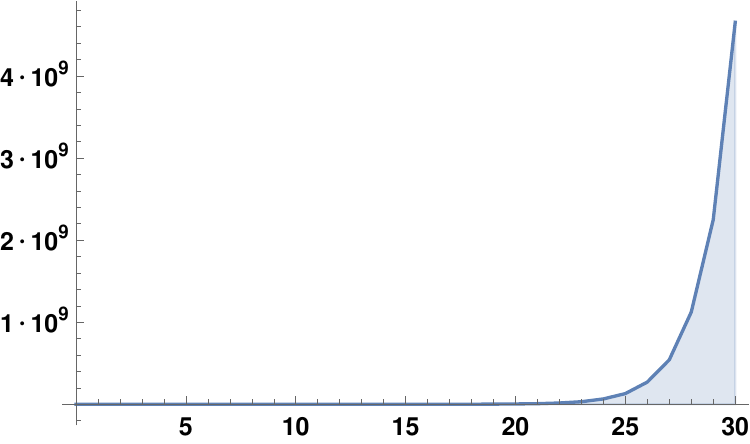}};
\node at (1.5,2) {$\tlmon[n]/\K$};
\node at (1.5,1.5) {$n$ increases $\rightarrow$};
\draw[->] (-2.4,2.2) node[above,xshift=0.3cm]{y-axis: lower bound ssgap} to (-2.4,2.1) to (-2.8,2.1);
\draw[->] (1.5,-1.5)node[above]{x-axis: n} to (1.5,-1.8);
\end{tikzpicture}$}
,
\\
\scalebox{0.9}{$\begin{tikzpicture}[anchorbase]
\node at (0,0) {\includegraphics[height=4.4cm]{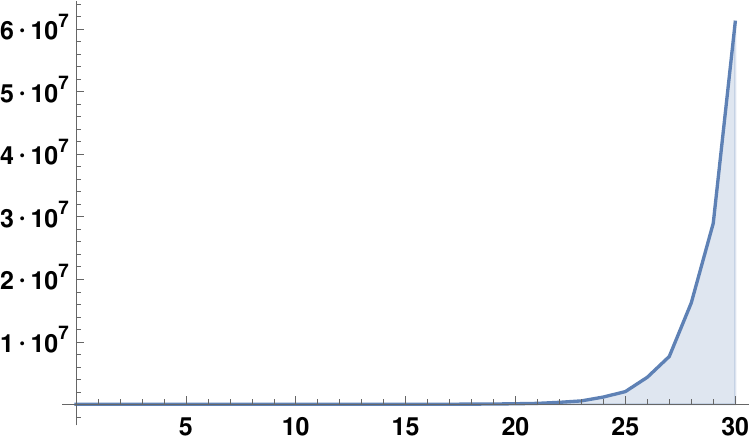}};
\node at (1.5,2) {$\tlmon[n]/\Q$};
\node at (1.5,1.5) {$n$ increases $\rightarrow$};
\draw[->] (-2.4,2.2) node[above,xshift=0.3cm]{y-axis: lower bound faith} to (-2.4,2.1) to (-2.8,2.1);
\draw[->] (1.5,-1.5)node[above]{x-axis: n} to (1.5,-1.8);
\end{tikzpicture}$}
.
\end{gather*}
In these plots $n$ increases form $0$ to $30$ when going left to right.
\end{Example}

Note that the bound $0\leq k\leq 2\sqrt{n}$ in \autoref{T:TLIsAGoodExample} means that 
the monoid $\tltru{k}$ has few through strands. This has the advantage that the 
dimensions of simple $\tltru{k}$-representations peak, but it also means that the 
information loss during multiplication is big. Alternatively one might want to 
keep $k$ close to $n$, so we also state:

\begin{Theorem}\label{T:TLIsAGoodExampleTwo}
Let $n\geq 8$ and fix $2\sqrt{n}\leq k\leq n-\sqrt{n}$. 
Let $\cchar=0$, and let $\LL$ be an arbitrary field.
We have the following lower bounds:
\begin{align*}
\gap[\K]{\tltru{k}}&\geq\frac{1}{(n-\tfrac{\sqrt{n}}{2}+1)(n-\tfrac{\sqrt{n}}{2}+2)}\binom{n}{\tfrac{\sqrt{n}}{2}}
\in\Theta\big(n^{\sqrt{n}/4}n^{-9/4}(2e)^{\sqrt{n}/2}\big)
,\\
\ssgap[\LL]{\tltru{k}}&\geq\frac{n-\sqrt{n}+1}{n-\sqrt{n}/2+1}\binom{n}{\tfrac{\sqrt{n}}{2}}
\in\Theta\big(n^{\sqrt{n}/4}n^{-3/4}(2e)^{\sqrt{n}/2}\big)
,\\
\faith[\K]{\tltru{k}}&\geq\frac{6}{n+4}\binom{n}{\tfrac{n^{\prime}}{2}-1}
\in\Theta\big(2^{n}n^{-3/2}\big)
.
\end{align*}
where in the final bound $n^{\prime}=n$, if $n$ is even, 
and $n^{\prime}=n-1$, if $n$ is odd.
\end{Theorem}

\begin{proof}
Similar to the proof of \autoref{T:TLIsAGoodExample} and omitted. 
(The assumption $n\geq 8$ ensures that $k\geq 3$, so we can use \autoref{L:TLWellRounded}.)
\end{proof}

\begin{Example}\label{E:TLRatio}
The various ratios from \autoref{SS:RepGapRatios} are easy 
to compute using \autoref{P:TLCells}, which gives $|\tltru[16]{8}|=\sum_{l\in\Lambda,l\geq k}|\jcell_{l}|$, 
and either of the theorems 
above. Explicitly, for $n=16$, $k=2\sqrt{16}=8$ and using 
$\K=\Q$, which is the setting from \autoref{T:TLIsAGoodExample}, we get
$\gratio[\Q]{\tltru[16]{8}}\geq 1.686\cdot 10^{-3}$. For comparison, 
the symmetric group $\sym[16]$ has $\gratio[\Q]{\sym[16]}\approx 2.186\cdot 10^{-7}$.
\end{Example}


\subsection{Other planar monoids}\label{SS:TLOther}

Let us now discuss cells, simples and bounds 
for the other planar monoids from \autoref{Eq:IntroDiaMonoids} 
in ascending order (of complexity). The constructions 
and statements are very similar to the Temperley--Lieb case, 
so we will be brief. The reader can find 
more details about the basics about the diagram 
monoids, and also references, in 
{\eg} \cite{HaJa-representations-diagram-algebras}.

\begin{Remark}\label{R:TLOther}
As we will see, the common theoretical feature 
of planar monoids is that their $H$-cells are all of size one.
As for the Temperley--Lieb monoid, the diagrammatically 
firm reader can deduce the cell structure of the planar
diagram monoids in 
this section themselves. There are also many references in the literature 
and the cells of these diagram monoids have been rediscovered many times. 
For example, \cite{DoEaGr-motzkin-pbrauer} computes the cells of $\promon$, 
and \cite{BeHa-motzkin} computes the cells of $\momon$.
\end{Remark}

We leave the case of the \emph{planar symmetric group} to the reader and start with
the \emph{planar rook monoid} $\promon$. 
This monoid was rediscovered several times, see {\eg} 
\cite{KhSa-cat-polyring}, and the reader might 
know it under a different name.
The construction of $\promon$ is almost the 
same as for $\tlmon$, but instead of caps and cups we have end and start dots, and all internal components are removed whenever they appear during composition.
The monoid $\promon$ has $\binom{2n}{n}$ elements and 
a typical cell is of the form
\begin{gather}\label{Eq:TLCellRook}
\xy
(0,0)*{\begin{gathered}
\begin{tabular}{C|C|C}
\arrayrulecolor{tomato}
\cellcolor{mydarkblue!25}
\begin{tikzpicture}[anchorbase]
\draw[usual,dot] (0,0) to (0,0.15);
\draw[usual,dot] (0.5,0) to (0.5,0.15);
\draw[usual] (1,0) to (1,0.5);
\draw[usual,dot] (0,0.5) to (0,0.35);
\draw[usual,dot] (0.5,0.5) to (0.5,0.35);
\end{tikzpicture} & 
\begin{tikzpicture}[anchorbase,xscale=-1]
\draw[usual,dot] (0,0) to (0,0.15);
\draw[usual] (0.5,0) to (0,0.5);
\draw[usual,dot] (1,0) to (1,0.15);
\draw[usual,dot] (0.5,0.5) to (0.5,0.35);
\draw[usual,dot] (1,0.5) to (1,0.35);
\end{tikzpicture} &
\begin{tikzpicture}[anchorbase]
\draw[usual] (0,0) to (1,0.5);
\draw[usual,dot] (0.5,0) to (0.5,0.15);
\draw[usual,dot] (1,0) to (1,0.15);
\draw[usual,dot] (0.5,0.5) to (0.5,0.35);
\draw[usual,dot] (0,0.5) to (0,0.35);
\end{tikzpicture}
\\
\hline
\begin{tikzpicture}[anchorbase]
\draw[usual,dot] (0,0) to (0,0.15);
\draw[usual,dot] (0.5,0) to (0.5,0.15);
\draw[usual] (1,0) to (0.5,0.5);
\draw[usual,dot] (0,0.5) to (0,0.35);
\draw[usual,dot] (1,0.5) to (1,0.35);
\end{tikzpicture} & 
\cellcolor{mydarkblue!25}
\begin{tikzpicture}[anchorbase]
\draw[usual,dot] (0,0) to (0,0.15);
\draw[usual] (0.5,0) to (0.5,0.5);
\draw[usual,dot] (1,0) to (1,0.15);
\draw[usual,dot] (0,0.5) to (0,0.35);
\draw[usual,dot] (1,0.5) to (1,0.35);
\end{tikzpicture} &
\begin{tikzpicture}[anchorbase]
\draw[usual] (0,0) to (0.5,0.5);
\draw[usual,dot] (0.5,0) to (0.5,0.15);
\draw[usual,dot] (1,0) to (1,0.15);
\draw[usual,dot] (0,0.5) to (0,0.35);
\draw[usual,dot] (1,0.5) to (1,0.35);
\end{tikzpicture}
\\
\hline
\begin{tikzpicture}[anchorbase]
\draw[usual,dot] (0,0) to (0,0.15);
\draw[usual,dot] (0.5,0) to (0.5,0.15);
\draw[usual] (1,0) to (0,0.5);
\draw[usual,dot] (1,0.5) to (1,0.35);
\draw[usual,dot] (0.5,0.5) to (0.5,0.35);
\end{tikzpicture} & 
\begin{tikzpicture}[anchorbase,xscale=-1]
\draw[usual,dot] (0,0) to (0,0.15);
\draw[usual] (0.5,0) to (1,0.5);
\draw[usual,dot] (1,0) to (1,0.15);
\draw[usual,dot] (0,0.5) to (0,0.35);
\draw[usual,dot] (0.5,0.5) to (0.5,0.35);
\end{tikzpicture} &
\cellcolor{mydarkblue!25}
\begin{tikzpicture}[anchorbase]
\draw[usual] (0,0) to (0,0.5);
\draw[usual,dot] (0.5,0) to (0.5,0.15);
\draw[usual,dot] (1,0) to (1,0.15);
\draw[usual,dot] (0.5,0.5) to (0.5,0.35);
\draw[usual,dot] (1,0.5) to (1,0.35);
\end{tikzpicture}
\end{tabular}
\end{gathered}};
(-30,0)*{\jcell_{1}};
(35,0)*{\hcell(e)\cong\onemon};
\endxy
\quad.
\end{gather}
This illustrates $\jcell_{1}$ of $\promon[3]$, which has one through strand.

The monoid containing both, $\tlmon$ and $\promon$, as submonoids 
is the \emph{Motzkin monoid} $\momon$. The definition of this monoid 
works {\muta} as for $\tlmon$ and $\promon$, now with caps and cups as well as start and end dots, and all internal components are removed whenever they appear during composition.
The Motzkin monoid has $\sum_{k=0}^{n}\tfrac{1}{k+1}\binom{2n}{2k}\binom{2k}{k}$ elements. The $J$-cells $\jcell_{i}$ are still given 
by through strands $k$, and a prototypical example is
\begin{gather*}
\xy
(0,0)*{\begin{gathered}
\begin{tabular}{C|C|C|C|C}
\arrayrulecolor{tomato}
\cellcolor{mydarkblue!25}
\begin{tikzpicture}[anchorbase]
\draw[usual] (0,0) to[out=90,in=180] (0.25,0.2) to[out=0,in=90] (0.5,0);
\draw[usual] (0,0.5) to[out=270,in=180] (0.25,0.3) to[out=0,in=270] (0.5,0.5);
\draw[usual] (1,0) to (1,0.5);
\end{tikzpicture} & 
\cellcolor{mydarkblue!25}
\begin{tikzpicture}[anchorbase]
\draw[usual] (0.5,0) to[out=90,in=180] (0.75,0.2) to[out=0,in=90] (1,0);
\draw[usual] (0,0.5) to[out=270,in=180] (0.25,0.3) to[out=0,in=270] (0.5,0.5);
\draw[usual] (0,0) to (1,0.5);
\end{tikzpicture} &
\cellcolor{mydarkblue!25}
\begin{tikzpicture}[anchorbase]
\draw[usual,dot] (0,0) to (0,0.15);
\draw[usual,dot] (0.5,0) to (0.5,0.15);
\draw[usual] (0,0.5) to[out=270,in=180] (0.25,0.3) to[out=0,in=270] (0.5,0.5);
\draw[usual] (1,0) to (1,0.5);
\end{tikzpicture} &
\begin{tikzpicture}[anchorbase]
\draw[usual,dot] (0,0) to (0,0.15);
\draw[usual,dot] (1,0) to (1,0.15);
\draw[usual] (0,0.5) to[out=270,in=180] (0.25,0.3) to[out=0,in=270] (0.5,0.5);
\draw[usual] (0.5,0) to (1,0.5);
\end{tikzpicture} &
\begin{tikzpicture}[anchorbase]
\draw[usual,dot] (1,0) to (1,0.15);
\draw[usual,dot] (0.5,0) to (0.5,0.15);
\draw[usual] (0,0.5) to[out=270,in=180] (0.25,0.3) to[out=0,in=270] (0.5,0.5);
\draw[usual] (0,0) to (1,0.5);
\end{tikzpicture}
\\[3pt]
\cellcolor{mydarkblue!25}
\begin{tikzpicture}[anchorbase]
\draw[usual] (0,0) to[out=90,in=180] (0.25,0.2) to[out=0,in=90] (0.5,0);
\draw[usual] (0.5,0.5) to[out=270,in=180] (0.75,0.3) to[out=0,in=270] (1,0.5);
\draw[usual] (1,0) to (0,0.5);
\end{tikzpicture} &
\cellcolor{mydarkblue!25}
\begin{tikzpicture}[anchorbase]
\draw[usual] (0.5,0) to[out=90,in=180] (0.75,0.2) to[out=0,in=90] (1,0);
\draw[usual] (0.5,0.5) to[out=270,in=180] (0.75,0.3) to[out=0,in=270] (1,0.5);
\draw[usual] (0,0) to (0,0.5);
\end{tikzpicture} &
\begin{tikzpicture}[anchorbase]
\draw[usual,dot] (0,0) to (0,0.15);
\draw[usual,dot] (0.5,0) to (0.5,0.15);
\draw[usual] (0.5,0.5) to[out=270,in=180] (0.75,0.3) to[out=0,in=270] (1,0.5);
\draw[usual] (1,0) to (0,0.5);
\end{tikzpicture} &
\begin{tikzpicture}[anchorbase]
\draw[usual,dot] (0,0) to (0,0.15);
\draw[usual,dot] (1,0) to (1,0.15);
\draw[usual] (0.5,0.5) to[out=270,in=180] (0.75,0.3) to[out=0,in=270] (1,0.5);
\draw[usual] (0.5,0) to (0,0.5);
\end{tikzpicture} &
\cellcolor{mydarkblue!25}
\begin{tikzpicture}[anchorbase]
\draw[usual,dot] (1,0) to (1,0.15);
\draw[usual,dot] (0.5,0) to (0.5,0.15);
\draw[usual] (0.5,0.5) to[out=270,in=180] (0.75,0.3) to[out=0,in=270] (1,0.5);
\draw[usual] (0,0) to (0,0.5);
\end{tikzpicture}
\\[3pt]
\cellcolor{mydarkblue!25}
\begin{tikzpicture}[anchorbase]
\draw[usual] (0,0) to[out=90,in=180] (0.25,0.2) to[out=0,in=90] (0.5,0);
\draw[usual,dot] (0,0.5) to (0,0.35);
\draw[usual,dot] (0.5,0.5) to (0.5,0.35);
\draw[usual] (1,0) to (1,0.5);
\end{tikzpicture} &
\begin{tikzpicture}[anchorbase]
\draw[usual] (0.5,0) to[out=90,in=180] (0.75,0.2) to[out=0,in=90] (1,0);
\draw[usual,dot] (0,0.5) to (0,0.35);
\draw[usual,dot] (0.5,0.5) to (0.5,0.35);
\draw[usual] (0,0) to (1,0.5);
\end{tikzpicture} &
\cellcolor{mydarkblue!25}
\begin{tikzpicture}[anchorbase]
\draw[usual,dot] (0,0) to (0,0.15);
\draw[usual,dot] (0.5,0) to (0.5,0.15);
\draw[usual,dot] (0,0.5) to (0,0.35);
\draw[usual,dot] (0.5,0.5) to (0.5,0.35);
\draw[usual] (1,0) to (1,0.5);
\end{tikzpicture} &
\begin{tikzpicture}[anchorbase]
\draw[usual,dot] (0,0) to (0,0.15);
\draw[usual,dot] (1,0) to (1,0.15);
\draw[usual,dot] (0,0.5) to (0,0.35);
\draw[usual,dot] (0.5,0.5) to (0.5,0.35);
\draw[usual] (0.5,0) to (1,0.5);
\end{tikzpicture} &
\begin{tikzpicture}[anchorbase]
\draw[usual,dot] (1,0) to (1,0.15);
\draw[usual,dot] (0.5,0) to (0.5,0.15);
\draw[usual,dot] (0,0.5) to (0,0.35);
\draw[usual,dot] (0.5,0.5) to (0.5,0.35);
\draw[usual] (0,0) to (1,0.5);
\end{tikzpicture}
\\[3pt]
\begin{tikzpicture}[anchorbase]
\draw[usual] (0,0) to[out=90,in=180] (0.25,0.2) to[out=0,in=90] (0.5,0);
\draw[usual,dot] (0,0.5) to (0,0.35);
\draw[usual,dot] (1,0.5) to (1,0.35);
\draw[usual] (1,0) to (0.5,0.5);
\end{tikzpicture} &
\begin{tikzpicture}[anchorbase]
\draw[usual] (0.5,0) to[out=90,in=180] (0.75,0.2) to[out=0,in=90] (1,0);
\draw[usual,dot] (0,0.5) to (0,0.35);
\draw[usual,dot] (1,0.5) to (1,0.35);
\draw[usual] (0,0) to (0.5,0.5);
\end{tikzpicture} &
\begin{tikzpicture}[anchorbase]
\draw[usual,dot] (0,0) to (0,0.15);
\draw[usual,dot] (0.5,0) to (0.5,0.15);
\draw[usual,dot] (0,0.5) to (0,0.35);
\draw[usual,dot] (1,0.5) to (1,0.35);
\draw[usual] (1,0) to (0.5,0.5);
\end{tikzpicture} &
\cellcolor{mydarkblue!25}
\begin{tikzpicture}[anchorbase]
\draw[usual,dot] (0,0) to (0,0.15);
\draw[usual,dot] (1,0) to (1,0.15);
\draw[usual,dot] (0,0.5) to (0,0.35);
\draw[usual,dot] (1,0.5) to (1,0.35);
\draw[usual] (0.5,0) to (0.5,0.5);
\end{tikzpicture} &
\begin{tikzpicture}[anchorbase]
\draw[usual,dot] (1,0) to (1,0.15);
\draw[usual,dot] (0.5,0) to (0.5,0.15);
\draw[usual,dot] (0,0.5) to (0,0.35);
\draw[usual,dot] (1,0.5) to (1,0.35);
\draw[usual] (0,0) to (0.5,0.5);
\end{tikzpicture}
\\[3pt]
\begin{tikzpicture}[anchorbase]
\draw[usual] (0,0) to[out=90,in=180] (0.25,0.2) to[out=0,in=90] (0.5,0);
\draw[usual,dot] (0.5,0.5) to (0.5,0.35);
\draw[usual,dot] (1,0.5) to (1,0.35);
\draw[usual] (1,0) to (0,0.5);
\end{tikzpicture} &
\cellcolor{mydarkblue!25}
\begin{tikzpicture}[anchorbase]
\draw[usual] (0.5,0) to[out=90,in=180] (0.75,0.2) to[out=0,in=90] (1,0);
\draw[usual,dot] (0.5,0.5) to (0.5,0.35);
\draw[usual,dot] (1,0.5) to (1,0.35);
\draw[usual] (0,0) to (0,0.5);
\end{tikzpicture} &
\begin{tikzpicture}[anchorbase]
\draw[usual,dot] (0,0) to (0,0.15);
\draw[usual,dot] (0.5,0) to (0.5,0.15);
\draw[usual,dot] (0.5,0.5) to (0.5,0.35);
\draw[usual,dot] (1,0.5) to (1,0.35);
\draw[usual] (1,0) to (0,0.5);
\end{tikzpicture} &
\begin{tikzpicture}[anchorbase]
\draw[usual,dot] (0,0) to (0,0.15);
\draw[usual,dot] (1,0) to (1,0.15);
\draw[usual,dot] (0.5,0.5) to (0.5,0.35);
\draw[usual,dot] (1,0.5) to (1,0.35);
\draw[usual] (0.5,0) to (0,0.5);
\end{tikzpicture} &
\cellcolor{mydarkblue!25}
\begin{tikzpicture}[anchorbase]
\draw[usual,dot] (1,0) to (1,0.15);
\draw[usual,dot] (0.5,0) to (0.5,0.15);
\draw[usual,dot] (0.5,0.5) to (0.5,0.35);
\draw[usual,dot] (1,0.5) to (1,0.35);
\draw[usual] (0,0) to (0,0.5);
\end{tikzpicture}
\end{tabular}
\end{gathered}};
(-45,0)*{\jcell_{1}};
(50,0)*{\hcell(e)\cong\onemon};
\endxy
\quad.
\end{gather*}
This illustrates $\jcell_{1}$ and $\momon[3]$.

Finally, the \emph{planar partition monoid} $\ppamon$ 
has all of the above mentioned planar monoids as submonoids, 
as it allows now arbitrary partitions, and has $Ca(2n)$ elements. 
(Recall that $Ca(k)$ was the $k$th Catalan number.) As before, 
internal components are removed and cells look very familiar 
to the cells of the other planar monoids. For example
$\jcell_{1}$ for $\ppamon[2]$ is:
\begin{gather*}
\xy
(0,0)*{\begin{gathered}
\begin{tabular}{C|C|C}
\arrayrulecolor{tomato}
\cellcolor{mydarkblue!25}
\begin{tikzpicture}[anchorbase]
\draw[usual,dot] (0,0) to (0,0.15);
\draw[usual,dot] (0,0.5) to (0,0.35);
\draw[usual] (0.5,0) to (0.5,0.5);
\end{tikzpicture} & 
\begin{tikzpicture}[anchorbase]
\draw[usual,dot] (0.5,0) to (0.5,0.15);
\draw[usual,dot] (0,0.5) to (0,0.35);
\draw[usual] (0,0) to (0.5,0.5);
\end{tikzpicture} &
\begin{tikzpicture}[anchorbase]
\draw[usual] (0,0) to[out=90,in=180] (0.25,0.2) to[out=0,in=90] (0.5,0);
\draw[usual,dot] (0,0.5) to (0,0.35);
\draw[usual] (0,0) to (0,0.18) to (0.5,0.5);
\end{tikzpicture}
\\
\hline
\begin{tikzpicture}[anchorbase]
\draw[usual,dot] (0,0) to (0,0.15);
\draw[usual,dot] (0.5,0.5) to (0.5,0.35);
\draw[usual] (0.5,0) to (0,0.5);
\end{tikzpicture} & 
\cellcolor{mydarkblue!25}
\begin{tikzpicture}[anchorbase]
\draw[usual,dot] (0.5,0) to (0.5,0.15);
\draw[usual,dot] (0.5,0.5) to (0.5,0.35);
\draw[usual] (0,0) to (0,0.5);
\end{tikzpicture} &
\begin{tikzpicture}[anchorbase]
\draw[usual] (0,0) to[out=90,in=180] (0.25,0.2) to[out=0,in=90] (0.5,0);
\draw[usual,dot] (0.5,0.5) to (0.5,0.35);
\draw[usual] (0,0) to (0,0.5);
\end{tikzpicture}
\\
\hline
\begin{tikzpicture}[anchorbase]
\draw[usual] (0.5,0) to (0,0.32) to (0,0.5);
\draw[usual,dot] (0,0) to (0,0.15);
\draw[usual] (0,0.5) to[out=270,in=180] (0.25,0.3) to[out=0,in=270] (0.5,0.5);
\end{tikzpicture} & 
\begin{tikzpicture}[anchorbase]
\draw[usual] (0,0) to (0,0.5);
\draw[usual,dot] (0.5,0) to (0.5,0.15);
\draw[usual] (0,0.5) to[out=270,in=180] (0.25,0.3) to[out=0,in=270] (0.5,0.5);
\end{tikzpicture} &
\cellcolor{mydarkblue!25}
\begin{tikzpicture}[anchorbase]
\draw[usual] (0,0) to (0,0.5);
\draw[usual] (0,0) to[out=90,in=180] (0.25,0.2) to[out=0,in=90] (0.5,0);
\draw[usual] (0,0.5) to[out=270,in=180] (0.25,0.3) to[out=0,in=270] (0.5,0.5);
\end{tikzpicture}
\end{tabular}
\end{gathered}};
(-25,0)*{\jcell_{1}};
(29,0)*{\hcell(e)\cong\onemon};
\endxy
\quad.
\end{gather*}

In the following we will focus on $\promon$ and $\momon$
as justified by our discussion of the Temperley--Lieb monoid and:

\begin{Lemma}\label{L:TLPPIsKnown}
There is an isomorphism of monoids $\ppamon\cong\tlmon[2n]$.
\end{Lemma}

\begin{proof}
See \cite[(1.5)]{HaRa-partition-algebras}.
\end{proof}

Not surprisingly, the analogs of \autoref{P:TLCells} and \autoref{P:TLSimples2} read almost the same. Below, if not stated otherwise, let $\xmon$ be either $\promon$ or $\momon$.

\begin{Proposition}\label{P:TLCells2}
We have the following.
\begin{enumerate}

\item The left and right cells of $\xmon$
are given by the respective type of diagrams where one fixes 
the bottom respectively top half of the diagram.
The $\leq_{l}$- and the $\leq_{r}$-order increases as the number of through strands 
decreases. Within $\jcell_{k}$ we have
\begin{gather*}
\promon\colon|\lcell|=|\rcell|=\binom{n}{k},
\\
\momon\colon|\lcell|=|\rcell|=\sum_{t=0}^{n}\frac{k+1}{k+t+1}\binom{n}{k+2t}
\binom{k+2t}{t}.
\end{gather*}

\item The $J$-cells $\jcell_{k}$ of $\xmon$
are given by the respective type of diagrams with a fixed number of through strands $k$.
The $\leq_{lr}$-order is a total order and increases as the number of through strands 
decreases. For any $\lcell\subset\jcell_{k}$ we have
\begin{gather*}
\xmon\colon|\jcell_{k}|=|\lcell|^{2}.
\end{gather*}

\item Each $J$-cell of $\xmon$ is idempotent, and $\hcell(e)\cong1$ for all idempotent $H$-cells. We have
\begin{gather*}
\xmon\colon|\hcell|=1.
\end{gather*}

\end{enumerate}	

\end{Proposition}

\begin{proof}
Omitted. See also \cite[Section 3.3]{HaJa-representations-diagram-algebras}. 
Note that the reference gives the dimensions of the simple 
$\promon$- and $\momon$-representations 
in the semisimple case, which are thus the sizes of 
the corresponding cells, see \autoref{P:CellsSemisimple}.
\end{proof}

\begin{Proposition}\label{P:TLSimples2}
The set of apexes for simple $\xmon$-representations can be indexed $1{:}1$ by 
the poset $\Lambda=(\{n,n-1,\dots\},>)$, and
there is precisely one simple $\xmon$-representation of a fixed apex up to $\cong$.
\end{Proposition}

\begin{proof}
Clear by \autoref{P:TLCells2}.
\end{proof}

We can number the simple 
$\xmon$-representations by $\simple[k]$ 
for $k\in\Lambda$.

\begin{Proposition}\label{P:TLSemisimpleDimsOther}
The semisimple dimensions for $\promon$ and $\momon$ are $\ssdimk(\simple[k])=\binom{n}{k}$ and
$\ssdimk(\simple[k])=\sum_{t=0}^{n}\tfrac{k+1}{k+t+1}\binom{n}{k+2t}\binom{k+2t}{t}$, respectively.
\end{Proposition}

\begin{proof}
Directly from \autoref{P:TLCells2} and \autoref{P:TLSimples2}.
\end{proof}

The semisimple dimensions of $\promon$ are given in \autoref{Eq:TlDimTL}. 
(Note that \autoref{Eq:TlDimTL} 
shows the dimensions of the simple $\promon$-representations, but we will see in \autoref{P:TLDimsPRook} below that $\dimk(\simple[k])=\ssdimk(\simple[k])$ holds for $\promon$.)
The semisimple dimensions of $\momon$ behave similarly as the 
semisimple dimensions of $\tlmon$, {\cf} \autoref{Eq:TlDimTL}:
\begin{gather*}
\begin{tikzpicture}[anchorbase]
\node at (0,0) {\includegraphics[height=4.4cm]{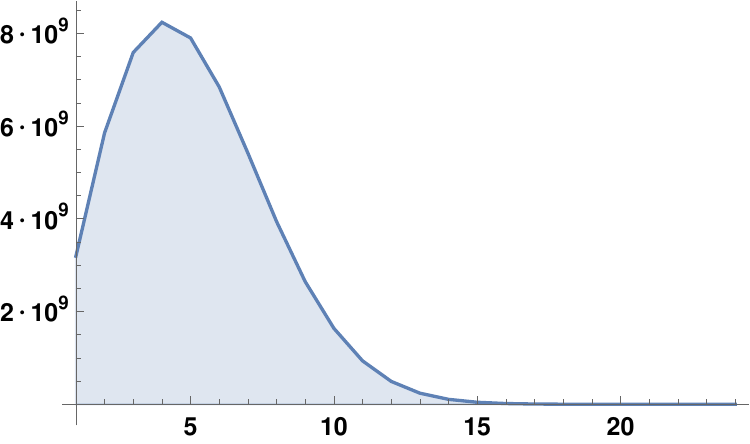}};
\node at (2.25,2) {$\momon[24]/\K$};
\node at (2.25,1.5) {cells increase $\leftarrow$};
\draw[->] (-2.4,2.2) node[above,xshift=0.3cm]{y-axis: ssdim} to (-2.4,2.1) to (-2.8,2.1);
\draw[->] (2.75,-1.5) node[above,yshift=0.75cm]{x-axis:}node[above,yshift=0.3cm]{\# through}node[above,yshift=0.0cm]{strands} to (2.75,-1.8);
\draw[very thick,densely dotted] (-0.7,-2) to (-0.7,1.5) node[above]{$k\approx 2\sqrt{24}$};
\end{tikzpicture}
,
\\
\scalebox{0.7}{$\text{ssdim}\colon
\left(
\begin{gathered}
3192727797,5850674704,7583013474,8234447672,7895719634,6839057544,
\\
5412710842,3938013264,2641866894,1636117512,935163394,492652824,
\\
238637282,105922544,42884259,15742672,5199909,1530144,395922,88504,16674,2552, 299,24,1
\end{gathered}
\right)$}
.
\end{gather*}

The dimensions of simple $\promon$-representations are easy to obtain:

\begin{Proposition}\label{P:TLDimsPRook}
We have $\dimk(\simple[k])=\ssdimk(\simple[k])=\binom{n}{k}$ for $\promon$, and $\promon$ is semisimple.
\end{Proposition}

\begin{proof}
We only need to prove that $\promon$ is semisimple, which implies the 
other results by \autoref{P:CellsSemisimple} and \autoref{P:TLSemisimpleDimsOther}.

To show semisimplicity we use \cite[Theorem 5.21]{St-rep-monoid} 
which says that a finite monoid is semisimple if and only if 
all $J$-cells are idempotent, all idempotent $H$-cells are semisimple 
and the Gram matrices $P(e)$, see \autoref{SS:CellsGram}, for all idempotent $H$-cell $\hcell(e)$ are invertible. 

By \autoref{P:TLCells2} we only need to compute the Gram matrices. 
Since $\hcell(e)\cong\onemon$, the reader familiar with the theory of cellular
algebras will recognize the following calculation.
The Gram matrix for any of the planar monoids discussed in this paper can be computed using analogs of \autoref{L:TLFactor}. Precisely, 
for each $J$-cell there are bottom diagrams $\beta_{1},\dots,\beta_{L}$
and top diagrams $\gamma_{1},\dots,\gamma_{L}$
indexing the rows and columns of the $J$-cell in question. The Gram matrix is then
\begin{gather*}
P(e)_{ij}=
\begin{cases}
1 & \text{if }\beta_{j}\gamma_{i}=1,
\\
0 & \text{else},
\end{cases}
\end{gather*}
where $1$ is the element of $\hcell(e)\cong\onemon$. For example, 
the Gram matrix of $\jcell_{1}$ of $\promon[3]$, see \autoref{Eq:TLCellRook},
takes the form
\begin{gather*}
\xy
(0,0)*{\begin{gathered}
\begin{tabular}{C||C|C|C}
\arrayrulecolor{tomato}
\beta/\gamma& \begin{tikzpicture}[anchorbase]
\draw[usual,dot] (0,0) to (0,0.15);
\draw[usual,dot] (0.5,0) to (0.5,0.15);
\draw[usual] (1,0) to (1,0.25);
\end{tikzpicture} & \begin{tikzpicture}[anchorbase]
\draw[usual,dot] (0,0) to (0,0.15);
\draw[usual,dot] (1,0) to (1,0.15);
\draw[usual] (0.5,0) to (0.5,0.25);
\end{tikzpicture} & \begin{tikzpicture}[anchorbase]
\draw[usual,dot] (1,0) to (1,0.15);
\draw[usual,dot] (0.5,0) to (0.5,0.15);
\draw[usual] (0,0) to (0,0.25);
\end{tikzpicture}
\\
\hline
\hline
\begin{tikzpicture}[anchorbase]
\draw[usual,dot] (0,0) to (0,-0.15);
\draw[usual,dot] (0.5,0) to (0.5,-0.15);
\draw[usual] (1,0) to (1,-0.25);
\end{tikzpicture}	& \begin{tikzpicture}[anchorbase]
\draw[usual,dot] (0,0) to (0,0.15);
\draw[usual,dot] (0.5,0) to (0.5,0.15);
\draw[usual] (1,0) to (1,0.25);
\draw[usual,dot] (0,0) to (0,-0.15);
\draw[usual,dot] (0.5,0) to (0.5,-0.15);
\draw[usual] (1,0) to (1,-0.25);
\end{tikzpicture} & \begin{tikzpicture}[anchorbase]
\draw[usual,dot] (0,0) to (0,0.15);
\draw[usual,dot] (1,0) to (1,0.15);
\draw[usual] (0.5,0) to (0.5,0.25);
\draw[usual,dot] (0,0) to (0,-0.15);
\draw[usual,dot] (0.5,0) to (0.5,-0.15);
\draw[usual] (1,0) to (1,-0.25);
\end{tikzpicture} & \begin{tikzpicture}[anchorbase]
\draw[usual,dot] (1,0) to (1,0.15);
\draw[usual,dot] (0.5,0) to (0.5,0.15);
\draw[usual] (0,0) to (0,0.25);
\draw[usual,dot] (0,0) to (0,-0.15);
\draw[usual,dot] (0.5,0) to (0.5,-0.15);
\draw[usual] (1,0) to (1,-0.25);
\end{tikzpicture}
\\
\hline
\begin{tikzpicture}[anchorbase]
\draw[usual,dot] (0,0) to (0,-0.15);
\draw[usual,dot] (1,0) to (1,-0.15);
\draw[usual] (0.5,0) to (0.5,-0.25);
\end{tikzpicture}	& \begin{tikzpicture}[anchorbase]
\draw[usual,dot] (0,0) to (0,0.15);
\draw[usual,dot] (0.5,0) to (0.5,0.15);
\draw[usual] (1,0) to (1,0.25);
\draw[usual,dot] (0,0) to (0,-0.15);
\draw[usual,dot] (1,0) to (1,-0.15);
\draw[usual] (0.5,0) to (0.5,-0.25);
\end{tikzpicture} & \begin{tikzpicture}[anchorbase]
\draw[usual,dot] (0,0) to (0,0.15);
\draw[usual,dot] (1,0) to (1,0.15);
\draw[usual] (0.5,0) to (0.5,0.25);
\draw[usual,dot] (0,0) to (0,-0.15);
\draw[usual,dot] (1,0) to (1,-0.15);
\draw[usual] (0.5,0) to (0.5,-0.25);
\end{tikzpicture} & \begin{tikzpicture}[anchorbase]
\draw[usual,dot] (1,0) to (1,0.15);
\draw[usual,dot] (0.5,0) to (0.5,0.15);
\draw[usual] (0,0) to (0,0.25);
\draw[usual,dot] (0,0) to (0,-0.15);
\draw[usual,dot] (1,0) to (1,-0.15);
\draw[usual] (0.5,0) to (0.5,-0.25);
\end{tikzpicture}
\\
\hline
\begin{tikzpicture}[anchorbase]
\draw[usual,dot] (1,0) to (1,-0.15);
\draw[usual,dot] (0.5,0) to (0.5,-0.15);
\draw[usual] (0,0) to (0,-0.25);
\end{tikzpicture}	& \begin{tikzpicture}[anchorbase]
\draw[usual,dot] (0,0) to (0,0.15);
\draw[usual,dot] (0.5,0) to (0.5,0.15);
\draw[usual] (1,0) to (1,0.25);
\draw[usual,dot] (1,0) to (1,-0.15);
\draw[usual,dot] (0.5,0) to (0.5,-0.15);
\draw[usual] (0,0) to (0,-0.25);
\end{tikzpicture} & \begin{tikzpicture}[anchorbase]
\draw[usual,dot] (0,0) to (0,0.15);
\draw[usual,dot] (1,0) to (1,0.15);
\draw[usual] (0.5,0) to (0.5,0.25);
\draw[usual,dot] (1,0) to (1,-0.15);
\draw[usual,dot] (0.5,0) to (0.5,-0.15);
\draw[usual] (0,0) to (0,-0.25);
\end{tikzpicture} & \begin{tikzpicture}[anchorbase]
\draw[usual,dot] (1,0) to (1,0.15);
\draw[usual,dot] (0.5,0) to (0.5,0.15);
\draw[usual] (0,0) to (0,0.25);
\draw[usual,dot] (1,0) to (1,-0.15);
\draw[usual,dot] (0.5,0) to (0.5,-0.15);
\draw[usual] (0,0) to (0,-0.25);
\end{tikzpicture}
\end{tabular}
\end{gathered}};
\endxy
\quad
\raisebox{-0.3cm}{$\leftrightsquigarrow
\begin{pmatrix}
1 & 0 & 0
\\
0 & 1 & 0
\\
0 & 0 & 1
\end{pmatrix}.$}
\end{gather*}
This is the identity matrix. In 
fact, $P(e)$ is always a permutation matrix:
any end dot 
needs to hit a start dot in order 
for $\beta_{j}\gamma_{i}$ to keep the 
same number of through strands, 
and there is precisely one $\beta_{j}$ for a fixed 
$\gamma_{i}$ for this to happen. The proof completes.

Alternatively, using an argument from monoid theory, 
$\promon$ is an inverse monoid, namely a submonoid of 
the symmetric inverse monoid that we will meet 
in \autoref{SS:BrauerOther} below. Moreover, by 
\autoref{P:TLCells2} we have $|\hcell|=1$. Thus,
\cite[Corollary 9.4]{St-rep-monoid} 
implies that $\promon$ is semisimple.
\end{proof}

The behavior 
of the dimensions of the simple $\promon$-representations is sketched in \autoref{Eq:CellsDimTL}. 
Sadly, we do not know the dimensions of the simple $\momon$-representations, 
but we have the following.

\begin{Proposition}\label{P:TLMotzkinBound}
Let $\simple[l]^{\tlmon}$ denote the $l$th simple 
$\tlmon$-representation, {\cf} \autoref{SS:TLCell}.
We have $\dimk(\simple[k])\geq\dimk(\simple[k]^{\tlmon})$, if $n-k$ is even,
and $\dimk(\simple[k])\geq\dimk(\simple[k]^{\tlmon[{n-1}]})$, if $n-k$ is odd,
both for $\momon$.
\end{Proposition}

\begin{proof}
Note that $\tlmon$ embeds into $\momon$ by sending 
every element to the element with the same description in $\momon$, {\eg}:
\begin{gather}\label{Eq:TLEmbedding}
\tlmon[3]\ni
\begin{tikzpicture}[anchorbase]
\draw[usual] (0,0) to[out=90,in=180] (0.25,0.2) to[out=0,in=90] (0.5,0);
\draw[usual] (0,0.5) to[out=270,in=180] (0.25,0.3) to[out=0,in=270] (0.5,0.5);
\draw[usual] (1,0) to (1,0.5);
\end{tikzpicture}
\mapsto
\begin{tikzpicture}[anchorbase]
\draw[usual] (0,0) to[out=90,in=180] (0.25,0.2) to[out=0,in=90] (0.5,0);
\draw[usual] (0,0.5) to[out=270,in=180] (0.25,0.3) to[out=0,in=270] (0.5,0.5);
\draw[usual] (1,0) to (1,0.5);
\end{tikzpicture}
\in\momon[3].
\end{gather}
Thus, $\tlmon$ is a submonoid of $\momon$ and \autoref{T:CellsGramSub} 
applies whenever $n-k$ is even since in this case $\jcell_{k}$ 
restricts to an idempotent $J$-cell of $\tlmon$. 

For the odd case we can use the same argument and the embedding 
of semigroups given by adding a pair of an end and a start dot to the right, {\eg}
\begin{gather*}
\tlmon[3]\ni
\begin{tikzpicture}[anchorbase]
\draw[usual] (0,0) to[out=90,in=180] (0.25,0.2) to[out=0,in=90] (0.5,0);
\draw[usual] (0,0.5) to[out=270,in=180] (0.25,0.3) to[out=0,in=270] (0.5,0.5);
\draw[usual] (1,0) to (1,0.5);
\end{tikzpicture}
\mapsto
\begin{tikzpicture}[anchorbase]
\draw[usual] (0,0) to[out=90,in=180] (0.25,0.2) to[out=0,in=90] (0.5,0);
\draw[usual] (0,0.5) to[out=270,in=180] (0.25,0.3) to[out=0,in=270] (0.5,0.5);
\draw[usual] (1,0) to (1,0.5);
\draw[usual,dot] (1.5,0) to (1.5,0.15);
\draw[usual,dot] (1.5,0.5) to (1.5,0.35);
\end{tikzpicture}
\in\momon[4].
\end{gather*}
\autoref{T:CellsGramSub} can be easily extended 
to cover this case as well.
\end{proof}

\begin{Lemma}\label{L:TLAdmissiblePRookMotzkin}
The monoid $\xmon$ is regular.
\end{Lemma}

\begin{proof}
Using \autoref{L:CellsAdmissible} the claim is 
easy to verify: for $\promon[n]$ and $\momon[n]$ symmetric diagrams 
with $k$ end and start dots gives an idempotent in $\jcell_{k}$. See also  \cite[Section 2]{DoEaGr-motzkin-pbrauer} for a proof using regularity.
\end{proof}

This suggests again that we use truncations. Note that $\protru{k}{l}$ 
below is constructed using an honest Rees factor, {\cf} 
\autoref{D:CellsSubquotient}, while $\motru{k}$ is a submonoid of $\momon$.

\begin{Definition}\label{D:TLTruncatedOther}
Define the \emph{$k$-$l$ truncated planar rook monoid} 
for $k\leq l$ and the
\emph{$k$th truncated Motzkin monoid} by
\begin{gather*}
\protru{k}{l}=(\promon[n])_{\geq\jcell_{k}}/(>\jcell_{l})
,\quad
\motru{k}=(\momon[n])_{\geq\jcell_{k}}.
\end{gather*}
\end{Definition}

Let $\xmon$ be either $\promon$ or $\protru{k}{l}$.

\begin{Proposition}\label{P:TLSplitsOther}
Let $\module$ be 
an $\xmon$-representation.
Then any short exact sequence 
\begin{gather*}
0\lra\onebt\lra\module\lra\onebt\lra 0
\end{gather*}
splits.
\end{Proposition} 

\begin{proof}
The monoid $\promon$ is semisimple, see 
\autoref{P:TLDimsPRook}, so \autoref{T:RepGapH1Condition} applies.
The case of $\protru{k}{l}$ follows {\ver} as 
the monoid is also semisimple by the analog of \autoref{P:TLDimsPRook}.
\end{proof} 

\begin{Remark}\label{R:TLSplitsOther}
To prove \autoref{P:TLSplitsOther} for the Motzkin monoid 
and its truncation it suffices to show that they are left-connected: that they are right-connected follows by applying the diagrammatic antiinvolution ${\placeholder}^{\ast}$, that they are null-connected follows from the fact that their $J$-cells are idempotent, the group of units $\group$ is trivial 
which implies $\HH^{1}(\group,\K)\cong 0$, and $\HH^{1}(\momon[n],\K)\cong 0$ as well as its counterpart for $\motru{k}$ follow from the same arguments as in the proof of \autoref{L:TLHTrivial}.
\end{Remark}

The following statement is only about $\promon$, since 
we do not know the dimensions of the simple 
$\momon$-representations.

\begin{Theorem}\label{P:TLRepGapOther}
We have
\begin{gather*}
\gap{\promon}=n,
\quad
\gap{\protru{k}{l}}=
\min\left\{\binom{n}{k},\binom{n}{l-1}\right\}.
\end{gather*}
\end{Theorem}

\begin{proof}
By \autoref{T:RepGapH1Condition}, \autoref{P:TLSplitsOther} and \autoref{P:TLDimsPRook}.
\end{proof}

\begin{Theorem}\label{T:TLIsAGoodExampleMotzkinPRook}
Let $k$ be arbitrary and $l=\floor{2\sqrt{n}}$.
We have the following lower bounds:
\begin{gather*}
\gap[\K]{\protru{l}{n-l}}=\ssgap[\K]{\protru{l}{n-l}}\geq
\binom{n}{\floor{2\sqrt{n}}}
,\\
\faith[\K]{\protru{l}{n-l}}\geq
\sqrt[2\floor{\sqrt{n}}+1]{\binom{n}{\floor{\tfrac{n}{2}}}}
,
\\
\ssgap[\K]{\motru{k}}\geq\ssgap[\LL]{\tltru[{n-1}]{k}}
,\quad
\faith[\K]{\motru{k}}\geq\faith[\K]{\tltru[{n-1}]{k}}
.
\end{gather*}
\end{Theorem}

Because of \autoref{P:TLMotzkinBound}, we also think that $\gap[\K]{\motru{k}}\geq\gap[\K]{\tltru[{n-1}]{k}}$. We can not prove this since we 
would need the analog of \autoref{P:TLSplitsOther} for the Motzkin monoid.

\begin{proof}
\textit{Planar rook.}
We start with $\promon$. The first inequality is 
immediate from \autoref{P:TLDimsPRook} and the behavior 
of binomial coefficients. For the second 
claim we apply \autoref{T:CellsFaithful}. Note that 
$\tltru{l,n-l}$ has $4\floor{\sqrt{n}}+1$ cells, but we 
can restrict to the submonoid with only $2\floor{\sqrt{n}}+1$ 
as in the proof of \autoref{T:TLIsAGoodExample}.

\textit{Motzkin.}
The first claim follows also from \autoref{T:TLIsAGoodExample} 
by identifying the smallest cell of Temperley--Lieb as 
a subcell of the smallest cell of $\motru{k}$.
The final inequality follows then from \autoref{P:TLMotzkinBound},
\autoref{T:CellsFaithful} and \autoref{T:TLIsAGoodExample}.
(Note that using $n-1$ is for convenience 
so that state a closed formula independent of even and odd issues.)
\end{proof}

\begin{Example}\label{E:TLRatioTwo}
As before for the Temperley--Lieb monoid, the various ratios 
are easy to get from the above. For example, $\gap{\protru[16]{6}{10}}\approx 0.34$.
\end{Example}

\begin{Conclusion}\label{C:TLConclusion}
From the viewpoint of linear attacks using 
small representations,
all of the planar monoids $\promon$, $\tlmon$, 
$\momon$ and $\ppamon$, or actually their truncations, 
have only big nontrivial representations.
However, $\tlmon$ 
is our main example: $\promon$ appears to be a bit 
too simple as a monoid to be of use and is semisimple, and $\ppamon$ is just 
$\tlmon[2n]$. The discussion about $\momon$ is unfinished 
and deserves more study.
\end{Conclusion}


\section{Symmetric monoids}\label{S:Brauer}

We still have a fixed field $\K$.


\subsection{Brauer categories and monoids}\label{SS:BrauerMonoid}

We will now recall the definitions of the \emph{Brauer category} $\brcat$, the \emph{Brauer algebra} $\bralg{\delta}$ and explain how to construct set-based versions of these. Brauer categories and algebras are classical topics in representation theory, see {\eg} \cite{Br-brauer-algebra-original} for the original reference. Moreover, the discussion is quite similar to the one in \autoref{S:TL}, so we will be brief.

The crucial difference between $\brcat$ and $\tlcat$ is that the former 
is additionally a symmetric category. The morphisms are 
then called \emph{perfect matchings}. Prototypical examples of these 
perfect matchings are crossingless matchings but also {\eg}:
\begin{gather*}
\begin{tikzpicture}[anchorbase]
\draw[usual] (0.5,0) to[out=90,in=180] (1,0.25) to[out=0,in=90] (1.5,0);
\draw[usual] (0,1) to[out=270,in=180] (0.25,0.75) to[out=0,in=270] (0.5,1);
\draw[usual] (0,0) to (1.5,1);
\draw[usual] (1,0) to (1,1);
\end{tikzpicture}
,\quad
\begin{tikzpicture}[anchorbase]
\draw[usual] (0,0) to[out=90,in=180] (0.75,0.3) to[out=0,in=90] (1.5,0);
\draw[usual] (0.5,0) to (1,1);
\draw[usual] (1,0) to (0.5,1);
\end{tikzpicture}
.
\end{gather*}
The relations on these diagrams are build such that they are the same
if and only if they represent the same perfect matching.
Otherwise the definition of $\brcat$ 
is the same as for $\tlcat$.

Perfect matchings can be numbered by $b(k)=(2k-1)!!$ (the double factorial). 
Letting $P_{m}^{n}$ denote the set of perfect matching 
with $m$ bottom and $n$ top boundary points, we have the following
analog of \autoref{L:TLBasis}:

\begin{Lemma}\label{L:BrauerBasis}
The set $P_{m}^{n}$ is a $\K$-linear basis of $\Hom_{\brcat}(m,n)$.
Hence, the dimension of this space is either zero if $m\not\equiv n\bmod 2$, or otherwise given by $\dimk\Hom_{\brcat}(m,n)=b(\tfrac{m+n}{2})$.
\end{Lemma}

\begin{proof}
Well-known, see {\eg} \cite[Lemma 4.4]{GrLe-cellular}.
\end{proof}

The Brauer algebra on $n$-strands is then 
$\bralg{\delta}=\End_{\brcat}(n)$.

\begin{Remark}\label{R:BrauerSchurWeyl}
Similar as $\tlalg{\delta}$, the algebra $\bralg{\delta}$ originates 
in Schur--Weyl--Brauer duality \cite{Br-brauer-algebra-original}.
See {\eg} \cite[Section 3.4]{AnStTu-semisimple-tilting} for a summary of these dualities.
\end{Remark}

The definition of the 
set theoretical version of these works {\ver} as in \autoref{D:TLSet}.
We then get the \emph{set-theoretic Brauer category} $\brset$ 
and the \emph{Brauer monoid} on $n$-strands is defined by $\brmon=\End_{\brset}(n)$. This monoid has $(2n-1)!!$ elements.

\begin{Lemma}\label{L:BrauerTL}
Sending the $\K$-linear basis from \autoref{L:TLBasis} 
to crossingless matching in $P_{m}^{n}$ from \autoref{L:BrauerBasis} 
defines an embedding of monoids $\tlmon\hookrightarrow\brmon$.
\end{Lemma}

\begin{proof}
Clear by the respective lemmas.
\end{proof}

Note that the symmetric group $\sym$ on $n$-strands 
is isomorphic to the group of units $\group$ of $\brmon$. 
An isomorphism is given by the map
\begin{gather*}
\sym\hookrightarrow\brmon,\quad
(i,i+1)\mapsto
\begin{tikzpicture}[anchorbase]
\draw[usual] (0,0) to (0.5,0.5);
\draw[usual] (0.5,0) to (0,0.5);
\end{tikzpicture}
,
\end{gather*}
where the crossing crosses the $i$th and the $(i+1)$th strand when read from left to right. We will use this to identify $\sym$ with the respective subgroup of $\brmon$ and with the corresponding set of morphisms in $\brset$.

The analog of \autoref{L:TLFactor} now is:

\begin{Lemma}\label{L:BrauerFactor}
For $a\in\Hom_{\brset}(m,n)$ there is unique factorization of the form $a=\gamma\circ\sigma_{k}\circ\beta$ for minimal $k$, and $\beta\in\Hom_{\brset}(m,k)$, $\sigma_{k}\in\sym$ and $\gamma\in\Hom_{\brset}(k,n)$.
\end{Lemma}

\begin{proof}
Very similar to the proof of \autoref{L:TLFactor}.
The picture now is
\begin{gather*}
a=
\begin{tikzpicture}[anchorbase]
\draw[usual] (0,1) to[out=270,in=180] (0.25,0.75) 
to[out=0,in=270] (0.5,1);
\draw[usual] (1,1) to[out=270,in=180] (1.75,0.75) to[out=0,in=270] (2.5,1);
\draw[usual] (0.5,0) to[out=90,in=180] (1,0.25) 
to[out=0,in=90] (1.5,0);
\draw[usual] (2.5,0) to[out=90,in=180] (2.75,0.25) 
to[out=0,in=90] (3,0);
\draw[usual] (0,0) to[out=90,in=180] (0.5,0.5) to (3,0.5) to[out=0,in=270] (3.5,1);
\draw[usual] (1,0) to (1.5,1);
\draw[usual] (2,0) to (2,1);
\end{tikzpicture}
=
\begin{tikzpicture}[anchorbase]
\draw[usual] (0,0) to[out=270,in=180] (0.25,-0.25) 
to[out=0,in=270] (0.5,0);
\draw[usual] (1,0) to[out=270,in=180] (1.75,-0.5) to[out=0,in=270] (2.5,0);
\draw[usual] (1.5,0) to (1.5,-1);
\draw[usual] (2,0) to (2,-1);
\draw[usual] (3,0) to (3,-1);
\end{tikzpicture}
\;\circ\;
\begin{tikzpicture}[anchorbase]
\draw[usual] (0,0) to (1,1);
\draw[usual] (0.5,0) to (0,1);
\draw[usual] (1,0) to (0.5,1);
\end{tikzpicture}
\;\circ\;
\begin{tikzpicture}[anchorbase]
\draw[usual] (-2.5,0) to (-2.5,1);
\draw[usual] (-1.5,0) to (-1.5,1);
\draw[usual] (-2,0) to[out=90,in=180] (-1.5,0.25) 
to[out=0,in=90] (-1,0);
\draw[usual] (-0.5,0) to (-0.5,1);
\draw[usual] (0,0) to[out=90,in=180] (0.25,0.25) 
to[out=0,in=90] (0.5,0);
\end{tikzpicture}
.
\end{gather*}
which one easily generalizes to prove the lemma. 
Note that one can always push crossing in the middle unless they 
have to cross a cap or cup.
\end{proof}

We apply the same terminology as for $\brmon$ regarding 
through strands, bottom half and top half.
As before for $\tlmon$, this notion will give us the cell structure 
of $\brmon$.


\subsection{Cells of the Brauer monoid}\label{SS:BrauerCells}

The cell structure of the Brauer monoid is as follows.

\begin{Remark}\label{R:BrauerCells}
As for all the other monoids 
we have seen, the computation of the cells 
of $\brmon$ is easy and a pleasant exercise. And, as before, there 
are plenty of references on the cell structure, see \cite{Br-algebra-orthogonal} for an early 
reference, \cite{FiGr-canonical-cases-brauer} for a reference from quantum algebra and \cite{KuMaMa-brauer-semigroup} for a reference from monoid theory.
\end{Remark}

The picture for the cell structure of $\brmon$ is now:
\begin{gather*}
\xy
(0,0)*{\begin{gathered}
\begin{tabular}{C|C|C}
\arrayrulecolor{tomato}
\cellcolor{mydarkblue!25}
\begin{tikzpicture}[anchorbase]
\draw[very thick] (0,0) to[out=90,in=180] (0.25,0.2) to[out=0,in=90] (0.5,0);
\draw[very thick] (0,0.5) to[out=270,in=180] (0.25,0.3) to[out=0,in=270] (0.5,0.5);
\draw[very thick] (1,0) to (1,0.5);
\end{tikzpicture} &
\cellcolor{mydarkblue!25}
\begin{tikzpicture}[anchorbase]
\draw[very thick] (0.5,0) to[out=90,in=180] (0.75,0.2) to[out=0,in=90] (1,0);
\draw[very thick] (0,0.5) to[out=270,in=180] (0.25,0.3) to[out=0,in=270] (0.5,0.5);
\draw[very thick] (0,0) to (1,0.5);
\end{tikzpicture} &
\cellcolor{mydarkblue!25}
\begin{tikzpicture}[anchorbase]
\draw[very thick] (0,0) to[out=90,in=180] (0.5,0.2) to[out=0,in=90] (1,0);
\draw[very thick] (0,0.5) to[out=270,in=180] (0.25,0.3) to[out=0,in=270] (0.5,0.5);
\draw[very thick] (0.5,0) to (1,0.5);
\end{tikzpicture}
\\
\hline
\cellcolor{mydarkblue!25}
\begin{tikzpicture}[anchorbase]
\draw[very thick] (0,0) to[out=90,in=180] (0.25,0.2) to[out=0,in=90] (0.5,0);
\draw[very thick] (0.5,0.5) to[out=270,in=180] (0.75,0.3) to[out=0,in=270] (1,0.5);
\draw[very thick] (1,0) to (0,0.5);
\end{tikzpicture} & 
\cellcolor{mydarkblue!25}
\begin{tikzpicture}[anchorbase]
\draw[very thick] (0,0) to (0,0.5);
\draw[very thick] (0.5,0) to[out=90,in=180] (0.75,0.2) to[out=0,in=90] (1,0);
\draw[very thick] (0.5,0.5) to[out=270,in=180] (0.75,0.3) to[out=0,in=270] (1,0.5);
\end{tikzpicture} &
\cellcolor{mydarkblue!25}
\begin{tikzpicture}[anchorbase]
\draw[very thick] (0,0) to[out=90,in=180] (0.5,0.2) to[out=0,in=90] (1,0);
\draw[very thick] (0.5,0.5) to[out=270,in=180] (0.75,0.3) to[out=0,in=270] (1,0.5);
\draw[very thick] (0.5,0) to (0,0.5);
\end{tikzpicture}
\\
\hline
\cellcolor{mydarkblue!25}
\begin{tikzpicture}[anchorbase]
\draw[very thick] (0,0) to[out=90,in=180] (0.25,0.2) to[out=0,in=90] (0.5,0);
\draw[very thick] (0,0.5) to[out=270,in=180] (0.5,0.3) to[out=0,in=270] (1,0.5);
\draw[very thick] (1,0) to (0.5,0.5);
\end{tikzpicture} & 
\cellcolor{mydarkblue!25}
\begin{tikzpicture}[anchorbase]
\draw[very thick] (0.5,0) to[out=90,in=180] (0.75,0.2) to[out=0,in=90] (1,0);
\draw[very thick] (0,0.5) to[out=270,in=180] (0.5,0.3) to[out=0,in=270] (1,0.5);
\draw[very thick] (0,0) to (0.5,0.5);
\end{tikzpicture} &
\cellcolor{mydarkblue!25}
\begin{tikzpicture}[anchorbase]
\draw[very thick] (0.5,0) to (0.5,0.5);
\draw[very thick] (0,0) to[out=90,in=180] (0.5,0.2) to[out=0,in=90] (1,0);
\draw[very thick] (0,0.5) to[out=270,in=180] (0.5,0.3) to[out=0,in=270] (1,0.5);
\end{tikzpicture}
\end{tabular}
\\[1pt]
\begin{tabular}{C}
\arrayrulecolor{tomato}
\cellcolor{mydarkblue!25}
\begin{gathered}
\begin{tikzpicture}[anchorbase]
\draw[very thick] (0,0) to (0,0.5);
\draw[very thick] (0.5,0) to (0.5,0.5);
\draw[very thick] (1,0) to (1,0.5);
\end{tikzpicture}
,
\begin{tikzpicture}[anchorbase]
\draw[very thick] (0,0) to (0.5,0.5);
\draw[very thick] (0.5,0) to (0,0.5);
\draw[very thick] (1,0) to (1,0.5);
\end{tikzpicture}
,
\begin{tikzpicture}[anchorbase]
\draw[very thick] (0,0) to (0,0.5);
\draw[very thick] (0.5,0) to (1,0.5);
\draw[very thick] (1,0) to (0.5,0.5);
\end{tikzpicture}
\\
\begin{tikzpicture}[anchorbase]
\draw[very thick] (0,0) to (1,0.5);
\draw[very thick] (0.5,0) to (0,0.5);
\draw[very thick] (1,0) to (0.5,0.5);
\end{tikzpicture}
,
\begin{tikzpicture}[anchorbase]
\draw[very thick] (0,0) to (0.5,0.5);
\draw[very thick] (0.5,0) to (1,0.5);
\draw[very thick] (1,0) to (0,0.5);
\end{tikzpicture}
,
\begin{tikzpicture}[anchorbase,rounded corners]
\draw[very thick] (0,0) to (1,0.5);
\draw[very thick] (0.5,0) to (0.9,0.25) to (0.5,0.5);
\draw[very thick] (1,0) to (0,0.5);
\end{tikzpicture}
\end{gathered}
\end{tabular}
\end{gathered}};
(-35,7)*{\jcell_{1}};
(-35,-8.5)*{\jcell_{3}};
(35,7)*{\hcell(e)\cong\sym[1]};
(35,-8.5)*{\hcell(e)\cong\sym[3]};
\endxy
\quad.
\end{gather*}
These are the cells of $\brmon[3]$.
Here is another example, where $\hcell(e)\cong\sym[2]$:
\begin{gather*}
\xy
(0,0)*{\begin{gathered}
\begin{tabular}{C|C|C|C|C|C}
\arrayrulecolor{tomato}
\begin{gathered}
\cellcolor{mydarkblue!25}
\begin{tikzpicture}[anchorbase,scale=0.7,tinynodes,yscale=-1]
\draw[usual] (0,0) to[out=90,in=180] (0.25,0.25) to[out=0,in=90] (0.5,0);
\draw[usual] (0,0.75) to[out=270,in=180] (0.25,0.5) to[out=0,in=270] (0.5,0.75);
\draw[usual] (1,0) to (1,0.75);
\draw[usual] (1.5,0) to (1.5,0.75);
\end{tikzpicture}
\\[3pt]
\begin{tikzpicture}[anchorbase,scale=0.7,tinynodes,yscale=-1]
\draw[usual] (0,0) to[out=90,in=180] (0.25,0.25) to[out=0,in=90] (0.5,0);
\draw[usual] (0,0.75) to[out=270,in=180] (0.25,0.5) to[out=0,in=270] (0.5,0.75);
\draw[usual] (1,0) to (1.5,0.75);
\draw[usual] (1.5,0) to (1,0.75);
\end{tikzpicture}
\end{gathered} &
\begin{gathered}
\cellcolor{mydarkblue!25}
\begin{tikzpicture}[anchorbase,scale=0.7,tinynodes,yscale=-1]
\draw[usual] (0,0) to[out=90,in=180] (0.25,0.25) to[out=0,in=90] (0.5,0);
\draw[usual] (0.5,0.75) to[out=270,in=180] (0.75,0.5) to[out=0,in=270] (1,0.75);
\draw[usual] (1,0) to (0,0.75);
\draw[usual] (1.5,0) to (1.5,0.75);
\end{tikzpicture}
\\[3pt]
\begin{tikzpicture}[anchorbase,scale=0.7,tinynodes,yscale=-1]
\draw[usual] (0,0) to[out=90,in=180] (0.25,0.25) to[out=0,in=90] (0.5,0);
\draw[usual] (0.5,0.75) to[out=270,in=180] (0.75,0.5) to[out=0,in=270] (1,0.75);
\draw[usual] (1.5,0) to (0,0.75);
\draw[usual] (1,0) to (1.5,0.75);
\end{tikzpicture}
\end{gathered} &
\begin{gathered}
\begin{tikzpicture}[anchorbase,scale=0.7,tinynodes,yscale=-1]
\draw[usual] (0,0) to[out=90,in=180] (0.25,0.25) to[out=0,in=90] (0.5,0);
\draw[usual] (1,0.75) to[out=270,in=180] (1.25,0.5) to[out=0,in=270] (1.5,0.75);
\draw[usual] (1,0) to (0,0.75);
\draw[usual] (1.5,0) to (0.5,0.75);
\end{tikzpicture}
\\[3pt]
\begin{tikzpicture}[anchorbase,scale=0.7,tinynodes,yscale=-1]
\draw[usual] (0,0) to[out=90,in=180] (0.25,0.25) to[out=0,in=90] (0.5,0);
\draw[usual] (1,0.75) to[out=270,in=180] (1.25,0.5) to[out=0,in=270] (1.5,0.75);
\draw[usual] (1,0) to (0.5,0.75);
\draw[usual] (1.5,0) to (0,0.75);
\end{tikzpicture}
\end{gathered} &
\begin{gathered}
\cellcolor{mydarkblue!25}
\begin{tikzpicture}[anchorbase,scale=0.7,tinynodes,yscale=-1]
\draw[usual] (0,0) to[out=90,in=180] (0.25,0.25) to[out=0,in=90] (0.5,0);
\draw[usual] (0,0.75) to[out=270,in=180] (0.5,0.5) to[out=0,in=270] (1,0.75);
\draw[usual] (1,0) to (0.5,0.75);
\draw[usual] (1.5,0) to (1.5,0.75);
\end{tikzpicture}
\\[3pt]
\begin{tikzpicture}[anchorbase,scale=0.7,tinynodes,yscale=-1]
\draw[usual] (0,0) to[out=90,in=180] (0.25,0.25) to[out=0,in=90] (0.5,0);
\draw[usual] (0,0.75) to[out=270,in=180] (0.5,0.5) to[out=0,in=270] (1,0.75);
\draw[usual] (1,0) to (1.5,0.75);
\draw[usual] (1.5,0) to (0.5,0.75);
\end{tikzpicture}
\end{gathered} &
\begin{gathered}
\cellcolor{mydarkblue!25}
\begin{tikzpicture}[anchorbase,scale=0.7,tinynodes,yscale=-1]
\draw[usual] (0,0) to[out=90,in=180] (0.25,0.25) to[out=0,in=90] (0.5,0);
\draw[usual] (0.5,0.75) to[out=270,in=180] (1,0.5) to[out=0,in=270] (1.5,0.75);
\draw[usual] (1,0) to (0,0.75);
\draw[usual] (1.5,0) to (1,0.75);
\end{tikzpicture}
\\[3pt]
\begin{tikzpicture}[anchorbase,scale=0.7,tinynodes,yscale=1]
\draw[usual] (0.5,0) to[out=90,in=180] (1,0.25) to[out=0,in=90] (1.5,0);
\draw[usual] (0,0.75) to[out=270,in=180] (0.25,0.5) to[out=0,in=270] (0.5,0.75);
\draw[usual] (1,0) to (1,0.75);
\draw[usual] (0,0) to (1.5,0.75);
\end{tikzpicture}
\end{gathered} &
\cellcolor{mydarkblue!25}
\begin{gathered}
\begin{tikzpicture}[anchorbase,scale=0.7,tinynodes,yscale=1]
\draw[usual] (0,0) to[out=90,in=180] (0.75,0.3) to[out=0,in=90] (1.5,0);
\draw[usual] (0,0.75) to[out=270,in=180] (0.25,0.5) to[out=0,in=270] (0.5,0.75);
\draw[usual] (0.5,0) to (1,0.75);
\draw[usual] (1,0) to (1.5,0.75);
\end{tikzpicture}
\\[3pt]
\begin{tikzpicture}[anchorbase,scale=0.7,tinynodes,yscale=1]
\draw[usual] (0,0) to[out=90,in=180] (0.75,0.3) to[out=0,in=90] (1.5,0);
\draw[usual] (0,0.75) to[out=270,in=180] (0.25,0.5) to[out=0,in=270] (0.5,0.75);
\draw[usual] (0.5,0) to (1.5,0.75);
\draw[usual] (1,0) to (1,0.75);
\end{tikzpicture}
\end{gathered}
\\
\arrayrulecolor{tomato}\hline
\begin{gathered}
\cellcolor{mydarkblue!25}
\begin{tikzpicture}[anchorbase,scale=0.7,tinynodes,yscale=-1]
\draw[usual] (0.5,0) to[out=90,in=180] (0.75,0.25) to[out=0,in=90] (1,0);
\draw[usual] (0,0.75) to[out=270,in=180] (0.25,0.5) to[out=0,in=270] (0.5,0.75);
\draw[usual] (0,0) to (1,0.75);
\draw[usual] (1.5,0) to (1.5,0.75);
\end{tikzpicture}
\\[3pt]
\begin{tikzpicture}[anchorbase,scale=0.7,tinynodes,yscale=-1]
\draw[usual] (0.5,0) to[out=90,in=180] (0.75,0.25) to[out=0,in=90] (1,0);
\draw[usual] (0,0.75) to[out=270,in=180] (0.25,0.5) to[out=0,in=270] (0.5,0.75);
\draw[usual] (0,0) to (1.5,0.75);
\draw[usual] (1.5,0) to (1,0.75);
\end{tikzpicture}
\end{gathered} &
\begin{gathered}
\cellcolor{mydarkblue!25}
\begin{tikzpicture}[anchorbase,scale=0.7,tinynodes,yscale=-1]
\draw[usual] (0.5,0) to[out=90,in=180] (0.75,0.25) to[out=0,in=90] (1,0);
\draw[usual] (0.5,0.75) to[out=270,in=180] (0.75,0.5) to[out=0,in=270] (1,0.75);
\draw[usual] (0,0) to (0,0.75);
\draw[usual] (1.5,0) to (1.5,0.75);
\end{tikzpicture}
\\[3pt]
\begin{tikzpicture}[anchorbase,scale=0.7,tinynodes,yscale=-1]
\draw[usual] (0.5,0) to[out=90,in=180] (0.75,0.25) to[out=0,in=90] (1,0);
\draw[usual] (0.5,0.75) to[out=270,in=180] (0.75,0.5) to[out=0,in=270] (1,0.75);
\draw[usual] (0,0) to (1.5,0.75);
\draw[usual] (1.5,0) to (0,0.75);
\end{tikzpicture}
\end{gathered} &
\begin{gathered}
\cellcolor{mydarkblue!25}
\begin{tikzpicture}[anchorbase,scale=0.7,tinynodes,yscale=-1]
\draw[usual] (0.5,0) to[out=90,in=180] (0.75,0.25) to[out=0,in=90] (1,0);
\draw[usual] (1,0.75) to[out=270,in=180] (1.25,0.5) to[out=0,in=270] (1.5,0.75);
\draw[usual] (0,0) to (0,0.75);
\draw[usual] (1.5,0) to (0.5,0.75);
\end{tikzpicture}
\\[3pt]
\begin{tikzpicture}[anchorbase,scale=0.7,tinynodes,yscale=-1]
\draw[usual] (0.5,0) to[out=90,in=180] (0.75,0.25) to[out=0,in=90] (1,0);
\draw[usual] (1,0.75) to[out=270,in=180] (1.25,0.5) to[out=0,in=270] (1.5,0.75);
\draw[usual] (1.5,0) to (0,0.75);
\draw[usual] (0,0) to (0.5,0.75);
\end{tikzpicture}
\end{gathered} & 
\begin{gathered}
\cellcolor{mydarkblue!25}
\begin{tikzpicture}[anchorbase,scale=0.7,tinynodes,yscale=-1]
\draw[usual] (0.5,0) to[out=90,in=180] (0.75,0.25) to[out=0,in=90] (1,0);
\draw[usual] (0,0.75) to[out=270,in=180] (0.5,0.5) to[out=0,in=270] (1,0.75);
\draw[usual] (0,0) to (0.5,0.75);
\draw[usual] (1.5,0) to (1.5,0.75);
\end{tikzpicture}
\\[3pt]
\begin{tikzpicture}[anchorbase,scale=0.7,tinynodes,yscale=-1]
\draw[usual] (0.5,0) to[out=90,in=180] (0.75,0.25) to[out=0,in=90] (1,0);
\draw[usual] (0,0.75) to[out=270,in=180] (0.5,0.5) to[out=0,in=270] (1,0.75);
\draw[usual] (0,0) to (1.5,0.75);
\draw[usual] (1.5,0) to (0.5,0.75);
\end{tikzpicture}
\end{gathered} &
\begin{gathered}
\cellcolor{mydarkblue!25}
\begin{tikzpicture}[anchorbase,scale=0.7,tinynodes,yscale=-1]
\draw[usual] (0.5,0) to[out=90,in=180] (0.75,0.25) to[out=0,in=90] (1,0);
\draw[usual] (0.5,0.75) to[out=270,in=180] (1,0.5) to[out=0,in=270] (1.5,0.75);
\draw[usual] (0,0) to (0,0.75);
\draw[usual] (1.5,0) to (1,0.75);
\end{tikzpicture}
\\[3pt]
\begin{tikzpicture}[anchorbase,scale=0.7,tinynodes,yscale=-1]
\draw[usual] (0.5,0) to[out=90,in=180] (0.75,0.25) to[out=0,in=90] (1,0);
\draw[usual] (0.5,0.75) to[out=270,in=180] (1,0.5) to[out=0,in=270] (1.5,0.75);
\draw[usual] (1.5,0) to (0,0.75);
\draw[usual] (0,0) to (1,0.75);
\end{tikzpicture}
\end{gathered} &
\begin{gathered}
\begin{tikzpicture}[anchorbase,scale=0.7,tinynodes,yscale=1]
\draw[usual] (0,0) to[out=90,in=180] (0.75,0.3) to[out=0,in=90] (1.5,0);
\draw[usual] (0.5,0.75) to[out=270,in=180] (0.75,0.5) to[out=0,in=270] (1,0.75);
\draw[usual] (0.5,0) to (0,0.75);
\draw[usual] (1,0) to (1.5,0.75);
\end{tikzpicture}
\\[3pt]
\begin{tikzpicture}[anchorbase,scale=0.7,tinynodes,yscale=1]
\draw[usual] (0,0) to[out=90,in=180] (0.75,0.3) to[out=0,in=90] (1.5,0);
\draw[usual] (0.5,0.75) to[out=270,in=180] (0.75,0.5) to[out=0,in=270] (1,0.75);
\draw[usual] (0.5,0) to (1.5,0.75);
\draw[usual] (1,0) to (0,0.75);
\end{tikzpicture}
\end{gathered}
\\
\arrayrulecolor{tomato}\hline
\begin{gathered}
\begin{tikzpicture}[anchorbase,scale=0.7,tinynodes,yscale=-1]
\draw[usual] (1,0) to[out=90,in=180] (1.25,0.25) to[out=0,in=90] (1.5,0);
\draw[usual] (0,0.75) to[out=270,in=180] (0.25,0.5) to[out=0,in=270] (0.5,0.75);
\draw[usual] (0,0) to (1,0.75);
\draw[usual] (0.5,0) to (1.5,0.75);
\end{tikzpicture}
\\[3pt]
\begin{tikzpicture}[anchorbase,scale=0.7,tinynodes,yscale=-1]
\draw[usual] (1,0) to[out=90,in=180] (1.25,0.25) to[out=0,in=90] (1.5,0);
\draw[usual] (0,0.75) to[out=270,in=180] (0.25,0.5) to[out=0,in=270] (0.5,0.75);
\draw[usual] (0,0) to (1.5,0.75);
\draw[usual] (0.5,0) to (1,0.75);
\end{tikzpicture}
\end{gathered} &
\begin{gathered}
\cellcolor{mydarkblue!25}
\begin{tikzpicture}[anchorbase,scale=0.7,tinynodes,yscale=-1]
\draw[usual] (1,0) to[out=90,in=180] (1.25,0.25) to[out=0,in=90] (1.5,0);
\draw[usual] (0.5,0.75) to[out=270,in=180] (0.75,0.5) to[out=0,in=270] (1,0.75);
\draw[usual] (0,0) to (0,0.75);
\draw[usual] (0.5,0) to (1.5,0.75);
\end{tikzpicture}
\\[3pt]
\begin{tikzpicture}[anchorbase,scale=0.7,tinynodes,yscale=-1]
\draw[usual] (1,0) to[out=90,in=180] (1.25,0.25) to[out=0,in=90] (1.5,0);
\draw[usual] (0.5,0.75) to[out=270,in=180] (0.75,0.5) to[out=0,in=270] (1,0.75);
\draw[usual] (0,0) to (1.5,0.75);
\draw[usual] (0.5,0) to (0,0.75);
\end{tikzpicture}
\end{gathered} &
\begin{gathered}
\cellcolor{mydarkblue!25}
\begin{tikzpicture}[anchorbase,scale=0.7,tinynodes,yscale=-1]
\draw[usual] (1,0) to[out=90,in=180] (1.25,0.25) to[out=0,in=90] (1.5,0);
\draw[usual] (1,0.75) to[out=270,in=180] (1.25,0.5) to[out=0,in=270] (1.5,0.75);
\draw[usual] (0,0) to (0,0.75);
\draw[usual] (0.5,0) to (0.5,0.75);
\end{tikzpicture}
\\[3pt]
\begin{tikzpicture}[anchorbase,scale=0.7,tinynodes,yscale=-1]
\draw[usual] (1,0) to[out=90,in=180] (1.25,0.25) to[out=0,in=90] (1.5,0);
\draw[usual] (1,0.75) to[out=270,in=180] (1.25,0.5) to[out=0,in=270] (1.5,0.75);
\draw[usual] (0,0) to (0.5,0.75);
\draw[usual] (0.5,0) to (0,0.75);
\end{tikzpicture}
\end{gathered} &
\begin{gathered}
\cellcolor{mydarkblue!25}
\begin{tikzpicture}[anchorbase,scale=0.7,tinynodes,yscale=-1]
\draw[usual] (1,0) to[out=90,in=180] (1.25,0.25) to[out=0,in=90] (1.5,0);
\draw[usual] (0,0.75) to[out=270,in=180] (0.5,0.5) to[out=0,in=270] (1,0.75);
\draw[usual] (0,0) to (0.5,0.75);
\draw[usual] (0.5,0) to (1.5,0.75);
\end{tikzpicture}
\\[3pt]
\begin{tikzpicture}[anchorbase,scale=0.7,tinynodes,yscale=-1]
\draw[usual] (1,0) to[out=90,in=180] (1.25,0.25) to[out=0,in=90] (1.5,0);
\draw[usual] (0,0.75) to[out=270,in=180] (0.5,0.5) to[out=0,in=270] (1,0.75);
\draw[usual] (0.5,0) to (0.5,0.75);
\draw[usual] (0,0) to (1.5,0.75);
\end{tikzpicture}
\end{gathered} &
\begin{gathered}
\cellcolor{mydarkblue!25}
\begin{tikzpicture}[anchorbase,scale=0.7,tinynodes,yscale=-1]
\draw[usual] (1,0) to[out=90,in=180] (1.25,0.25) to[out=0,in=90] (1.5,0);
\draw[usual] (0.5,0.75) to[out=270,in=180] (1,0.5) to[out=0,in=270] (1.5,0.75);
\draw[usual] (0,0) to (0,0.75);
\draw[usual] (0.5,0) to (1,0.75);
\end{tikzpicture}
\\[3pt]
\begin{tikzpicture}[anchorbase,scale=0.7,tinynodes,yscale=-1]
\draw[usual] (1,0) to[out=90,in=180] (1.25,0.25) to[out=0,in=90] (1.5,0);
\draw[usual] (0.5,0.75) to[out=270,in=180] (1,0.5) to[out=0,in=270] (1.5,0.75);
\draw[usual] (0.5,0) to (0,0.75);
\draw[usual] (0,0) to (1,0.75);
\end{tikzpicture}
\end{gathered} &
\cellcolor{mydarkblue!25}
\begin{gathered}
\begin{tikzpicture}[anchorbase,scale=0.7,tinynodes,yscale=1]
\draw[usual] (0,0) to[out=90,in=180] (0.75,0.3) to[out=0,in=90] (1.5,0);
\draw[usual] (1,0.75) to[out=270,in=180] (1.25,0.5) to[out=0,in=270] (1.5,0.75);
\draw[usual] (0.5,0) to (0,0.75);
\draw[usual] (1,0) to (0.5,0.75);
\end{tikzpicture}
\\[3pt]
\begin{tikzpicture}[anchorbase,scale=0.7,tinynodes,yscale=1]
\draw[usual] (0,0) to[out=90,in=180] (0.75,0.3) to[out=0,in=90] (1.5,0);
\draw[usual] (1,0.75) to[out=270,in=180] (1.25,0.5) to[out=0,in=270] (1.5,0.75);
\draw[usual] (0.5,0) to (0.5,0.75);
\draw[usual] (1,0) to (0,0.75);
\end{tikzpicture}
\end{gathered}
\\
\arrayrulecolor{tomato}\hline
\begin{gathered}
\cellcolor{mydarkblue!25}
\begin{tikzpicture}[anchorbase,scale=0.7,tinynodes,yscale=-1]
\draw[usual] (0,0) to[out=90,in=180] (0.5,0.25) to[out=0,in=90] (1,0);
\draw[usual] (0,0.75) to[out=270,in=180] (0.25,0.5) to[out=0,in=270] (0.5,0.75);
\draw[usual] (0.5,0) to (1,0.75);
\draw[usual] (1.5,0) to (1.5,0.75);
\end{tikzpicture}
\\[3pt]
\begin{tikzpicture}[anchorbase,scale=0.7,tinynodes,yscale=-1]
\draw[usual] (0,0) to[out=90,in=180] (0.5,0.25) to[out=0,in=90] (1,0);
\draw[usual] (0,0.75) to[out=270,in=180] (0.25,0.5) to[out=0,in=270] (0.5,0.75);
\draw[usual] (1.5,0) to (1,0.75);
\draw[usual] (0.5,0) to (1.5,0.75);
\end{tikzpicture}
\end{gathered} &
\begin{gathered}
\cellcolor{mydarkblue!25}
\begin{tikzpicture}[anchorbase,scale=0.7,tinynodes,yscale=-1]
\draw[usual] (0,0) to[out=90,in=180] (0.5,0.25) to[out=0,in=90] (1,0);
\draw[usual] (0.5,0.75) to[out=270,in=180] (0.75,0.5) to[out=0,in=270] (1,0.75);
\draw[usual] (0.5,0) to (0,0.75);
\draw[usual] (1.5,0) to (1.5,0.75);
\end{tikzpicture}
\\[3pt]
\begin{tikzpicture}[anchorbase,scale=0.7,tinynodes,yscale=-1]
\draw[usual] (0,0) to[out=90,in=180] (0.5,0.25) to[out=0,in=90] (1,0);
\draw[usual] (0.5,0.75) to[out=270,in=180] (0.75,0.5) to[out=0,in=270] (1,0.75);
\draw[usual] (0.5,0) to (1.5,0.75);
\draw[usual] (1.5,0) to (0,0.75);
\end{tikzpicture}
\end{gathered} &
\begin{gathered}
\cellcolor{mydarkblue!25}
\begin{tikzpicture}[anchorbase,scale=0.7,tinynodes,yscale=-1]
\draw[usual] (0,0) to[out=90,in=180] (0.5,0.25) to[out=0,in=90] (1,0);
\draw[usual] (1,0.75) to[out=270,in=180] (1.25,0.5) to[out=0,in=270] (1.5,0.75);
\draw[usual] (0.5,0) to (0,0.75);
\draw[usual] (1.5,0) to (0.5,0.75);
\end{tikzpicture}
\\[3pt]
\begin{tikzpicture}[anchorbase,scale=0.7,tinynodes,yscale=-1]
\draw[usual] (0,0) to[out=90,in=180] (0.5,0.25) to[out=0,in=90] (1,0);
\draw[usual] (1,0.75) to[out=270,in=180] (1.25,0.5) to[out=0,in=270] (1.5,0.75);
\draw[usual] (0.5,0) to (0.5,0.75);
\draw[usual] (1.5,0) to (0,0.75);
\end{tikzpicture}
\end{gathered} &
\begin{gathered}
\cellcolor{mydarkblue!25}
\begin{tikzpicture}[anchorbase,scale=0.7,tinynodes,yscale=-1]
\draw[usual] (0,0) to[out=90,in=180] (0.5,0.25) to[out=0,in=90] (1,0);
\draw[usual] (0,0.75) to[out=270,in=180] (0.5,0.5) to[out=0,in=270] (1,0.75);
\draw[usual] (0.5,0) to (0.5,0.75);
\draw[usual] (1.5,0) to (1.5,0.75);
\end{tikzpicture}
\\[3pt]
\begin{tikzpicture}[anchorbase,scale=0.7,tinynodes,yscale=-1]
\draw[usual] (0,0) to[out=90,in=180] (0.5,0.25) to[out=0,in=90] (1,0);
\draw[usual] (0,0.75) to[out=270,in=180] (0.5,0.5) to[out=0,in=270] (1,0.75);
\draw[usual] (0.5,0) to (1.5,0.75);
\draw[usual] (1.5,0) to (0.5,0.75);
\end{tikzpicture}
\end{gathered} & 
\begin{gathered}
\begin{tikzpicture}[anchorbase,scale=0.7,tinynodes,yscale=-1]
\draw[usual] (0,0) to[out=90,in=180] (0.5,0.25) to[out=0,in=90] (1,0);
\draw[usual] (0.5,0.75) to[out=270,in=180] (1,0.5) to[out=0,in=270] (1.5,0.75);
\draw[usual] (0.5,0) to (0,0.75);
\draw[usual] (1.5,0) to (1,0.75);
\end{tikzpicture}
\\[3pt]
\begin{tikzpicture}[anchorbase,scale=0.7,tinynodes,yscale=-1]
\draw[usual] (0,0) to[out=90,in=180] (0.5,0.25) to[out=0,in=90] (1,0);
\draw[usual] (0.5,0.75) to[out=270,in=180] (1,0.5) to[out=0,in=270] (1.5,0.75);
\draw[usual] (1.5,0) to (0,0.75);
\draw[usual] (0.5,0) to (1,0.75);
\end{tikzpicture}
\end{gathered} &
\cellcolor{mydarkblue!25}
\begin{gathered}
\begin{tikzpicture}[anchorbase,scale=0.7,tinynodes,yscale=1]
\draw[usual] (0,0) to[out=90,in=180] (0.75,0.3) to[out=0,in=90] (1.5,0);
\draw[usual] (0,0.75) to[out=270,in=180] (0.5,0.5) to[out=0,in=270] (1,0.75);
\draw[usual] (0.5,0) to (0.5,0.75);
\draw[usual] (1,0) to (1.5,0.75);
\end{tikzpicture}
\\[3pt]
\begin{tikzpicture}[anchorbase,scale=0.7,tinynodes,yscale=1]
\draw[usual] (0,0) to[out=90,in=180] (0.75,0.3) to[out=0,in=90] (1.5,0);
\draw[usual] (0,0.75) to[out=270,in=180] (0.5,0.5) to[out=0,in=270] (1,0.75);
\draw[usual] (0.5,0) to (1.5,0.75);
\draw[usual] (1,0) to (0.5,0.75);
\end{tikzpicture}
\end{gathered}
\\
\arrayrulecolor{tomato}\hline
\cellcolor{mydarkblue!25}
\begin{gathered}
\begin{tikzpicture}[anchorbase,scale=0.7,tinynodes,yscale=-1]
\draw[usual] (0.5,0) to[out=90,in=180] (1,0.25) to[out=0,in=90] (1.5,0);
\draw[usual] (0,0.75) to[out=270,in=180] (0.25,0.5) to[out=0,in=270] (0.5,0.75);
\draw[usual] (0,0) to (1,0.75);
\draw[usual] (1,0) to (1.5,0.75);
\end{tikzpicture}
\\[3pt]
\begin{tikzpicture}[anchorbase,scale=0.7,tinynodes,yscale=-1]
\draw[usual] (0.5,0) to[out=90,in=180] (1,0.25) to[out=0,in=90] (1.5,0);
\draw[usual] (0,0.75) to[out=270,in=180] (0.25,0.5) to[out=0,in=270] (0.5,0.75);
\draw[usual] (0,0) to (1.5,0.75);
\draw[usual] (1,0) to (1,0.75);
\end{tikzpicture}
\end{gathered} &
\begin{gathered}
\cellcolor{mydarkblue!25}
\begin{tikzpicture}[anchorbase,scale=0.7,tinynodes,yscale=-1]
\draw[usual] (0.5,0) to[out=90,in=180] (1,0.25) to[out=0,in=90] (1.5,0);
\draw[usual] (0.5,0.75) to[out=270,in=180] (0.75,0.5) to[out=0,in=270] (1,0.75);
\draw[usual] (0,0) to (0,0.75);
\draw[usual] (1,0) to (1.5,0.75);
\end{tikzpicture}
\\[3pt]
\begin{tikzpicture}[anchorbase,scale=0.7,tinynodes,yscale=-1]
\draw[usual] (0.5,0) to[out=90,in=180] (1,0.25) to[out=0,in=90] (1.5,0);
\draw[usual] (0.5,0.75) to[out=270,in=180] (0.75,0.5) to[out=0,in=270] (1,0.75);
\draw[usual] (0,0) to (1.5,0.75);
\draw[usual] (1,0) to (0,0.75);
\end{tikzpicture}
\end{gathered} &
\begin{gathered}
\cellcolor{mydarkblue!25}
\begin{tikzpicture}[anchorbase,scale=0.7,tinynodes,yscale=-1]
\draw[usual] (0.5,0) to[out=90,in=180] (1,0.25) to[out=0,in=90] (1.5,0);
\draw[usual] (1,0.75) to[out=270,in=180] (1.25,0.5) to[out=0,in=270] (1.5,0.75);
\draw[usual] (0,0) to (0,0.75);
\draw[usual] (1,0) to (0.5,0.75);
\end{tikzpicture}
\\[3pt]
\begin{tikzpicture}[anchorbase,scale=0.7,tinynodes,yscale=-1]
\draw[usual] (0.5,0) to[out=90,in=180] (1,0.25) to[out=0,in=90] (1.5,0);
\draw[usual] (1,0.75) to[out=270,in=180] (1.25,0.5) to[out=0,in=270] (1.5,0.75);
\draw[usual] (0,0) to (0.5,0.75);
\draw[usual] (1,0) to (0,0.75);
\end{tikzpicture}
\end{gathered} &
\begin{gathered}
\begin{tikzpicture}[anchorbase,scale=0.7,tinynodes,yscale=-1]
\draw[usual] (0.5,0) to[out=90,in=180] (1,0.25) to[out=0,in=90] (1.5,0);
\draw[usual] (0,0.75) to[out=270,in=180] (0.5,0.5) to[out=0,in=270] (1,0.75);
\draw[usual] (0,0) to (0.5,0.75);
\draw[usual] (1,0) to (1.5,0.75);
\end{tikzpicture}
\\[3pt]
\begin{tikzpicture}[anchorbase,scale=0.7,tinynodes,yscale=-1]
\draw[usual] (0.5,0) to[out=90,in=180] (1,0.25) to[out=0,in=90] (1.5,0);
\draw[usual] (0,0.75) to[out=270,in=180] (0.5,0.5) to[out=0,in=270] (1,0.75);
\draw[usual] (0,0) to (1.5,0.75);
\draw[usual] (1,0) to (0.5,0.75);
\end{tikzpicture}
\end{gathered} &
\begin{gathered}
\cellcolor{mydarkblue!25}
\begin{tikzpicture}[anchorbase,scale=0.7,tinynodes,yscale=-1]
\draw[usual] (0.5,0) to[out=90,in=180] (1,0.25) to[out=0,in=90] (1.5,0);
\draw[usual] (0.5,0.75) to[out=270,in=180] (1,0.5) to[out=0,in=270] (1.5,0.75);
\draw[usual] (0,0) to (0,0.75);
\draw[usual] (1,0) to (1,0.75);
\end{tikzpicture}
\\[3pt]
\begin{tikzpicture}[anchorbase,scale=0.7,tinynodes,yscale=-1]
\draw[usual] (0.5,0) to[out=90,in=180] (1,0.25) to[out=0,in=90] (1.5,0);
\draw[usual] (0.5,0.75) to[out=270,in=180] (1,0.5) to[out=0,in=270] (1.5,0.75);
\draw[usual] (0,0) to (1,0.75);
\draw[usual] (1,0) to (0,0.75);
\end{tikzpicture}
\end{gathered} &
\cellcolor{mydarkblue!25}
\begin{gathered}
\begin{tikzpicture}[anchorbase,scale=0.7,tinynodes,yscale=1]
\draw[usual] (0,0) to[out=90,in=180] (0.75,0.3) to[out=0,in=90] (1.5,0);
\draw[usual] (0.5,0.75) to[out=270,in=180] (1,0.5) to[out=0,in=270] (1.5,0.75);
\draw[usual] (0.5,0) to (0,0.75);
\draw[usual] (1,0) to (1,0.75);
\end{tikzpicture}
\\[3pt]
\begin{tikzpicture}[anchorbase,scale=0.7,tinynodes,yscale=1]
\draw[usual] (0,0) to[out=90,in=180] (0.75,0.3) to[out=0,in=90] (1.5,0);
\draw[usual] (0.5,0.75) to[out=270,in=180] (1,0.5) to[out=0,in=270] (1.5,0.75);
\draw[usual] (0.5,0) to (1,0.75);
\draw[usual] (1,0) to (0,0.75);
\end{tikzpicture}
\end{gathered}
\\
\arrayrulecolor{tomato}\hline
\cellcolor{mydarkblue!25}
\begin{gathered}
\begin{tikzpicture}[anchorbase,scale=0.7,tinynodes,yscale=-1]
\draw[usual] (0,0) to[out=90,in=180] (0.75,0.3) to[out=0,in=90] (1.5,0);
\draw[usual] (0,0.75) to[out=270,in=180] (0.25,0.5) to[out=0,in=270] (0.5,0.75);
\draw[usual] (0.5,0) to (1,0.75);
\draw[usual] (1,0) to (1.5,0.75);
\end{tikzpicture}
\\[3pt]
\begin{tikzpicture}[anchorbase,scale=0.7,tinynodes,yscale=-1]
\draw[usual] (0,0) to[out=90,in=180] (0.75,0.3) to[out=0,in=90] (1.5,0);
\draw[usual] (0,0.75) to[out=270,in=180] (0.25,0.5) to[out=0,in=270] (0.5,0.75);
\draw[usual] (0.5,0) to (1.5,0.75);
\draw[usual] (1,0) to (1,0.75);
\end{tikzpicture}
\end{gathered} &
\begin{gathered}
\begin{tikzpicture}[anchorbase,scale=0.7,tinynodes,yscale=-1]
\draw[usual] (0,0) to[out=90,in=180] (0.75,0.3) to[out=0,in=90] (1.5,0);
\draw[usual] (0.5,0.75) to[out=270,in=180] (0.75,0.5) to[out=0,in=270] (1,0.75);
\draw[usual] (0.5,0) to (0,0.75);
\draw[usual] (1,0) to (1.5,0.75);
\end{tikzpicture}
\\[3pt]
\begin{tikzpicture}[anchorbase,scale=0.7,tinynodes,yscale=-1]
\draw[usual] (0,0) to[out=90,in=180] (0.75,0.3) to[out=0,in=90] (1.5,0);
\draw[usual] (0.5,0.75) to[out=270,in=180] (0.75,0.5) to[out=0,in=270] (1,0.75);
\draw[usual] (0.5,0) to (1.5,0.75);
\draw[usual] (1,0) to (0,0.75);
\end{tikzpicture}
\end{gathered} &
\cellcolor{mydarkblue!25}
\begin{gathered}
\begin{tikzpicture}[anchorbase,scale=0.7,tinynodes,yscale=-1]
\draw[usual] (0,0) to[out=90,in=180] (0.75,0.3) to[out=0,in=90] (1.5,0);
\draw[usual] (1,0.75) to[out=270,in=180] (1.25,0.5) to[out=0,in=270] (1.5,0.75);
\draw[usual] (0.5,0) to (0,0.75);
\draw[usual] (1,0) to (0.5,0.75);
\end{tikzpicture}
\\[3pt]
\begin{tikzpicture}[anchorbase,scale=0.7,tinynodes,yscale=-1]
\draw[usual] (0,0) to[out=90,in=180] (0.75,0.3) to[out=0,in=90] (1.5,0);
\draw[usual] (1,0.75) to[out=270,in=180] (1.25,0.5) to[out=0,in=270] (1.5,0.75);
\draw[usual] (0.5,0) to (0.5,0.75);
\draw[usual] (1,0) to (0,0.75);
\end{tikzpicture}
\end{gathered} &
\cellcolor{mydarkblue!25}
\begin{gathered}
\begin{tikzpicture}[anchorbase,scale=0.7,tinynodes,yscale=-1]
\draw[usual] (0,0) to[out=90,in=180] (0.75,0.3) to[out=0,in=90] (1.5,0);
\draw[usual] (0,0.75) to[out=270,in=180] (0.5,0.5) to[out=0,in=270] (1,0.75);
\draw[usual] (0.5,0) to (0.5,0.75);
\draw[usual] (1,0) to (1.5,0.75);
\end{tikzpicture}
\\[3pt]
\begin{tikzpicture}[anchorbase,scale=0.7,tinynodes,yscale=-1]
\draw[usual] (0,0) to[out=90,in=180] (0.75,0.3) to[out=0,in=90] (1.5,0);
\draw[usual] (0,0.75) to[out=270,in=180] (0.5,0.5) to[out=0,in=270] (1,0.75);
\draw[usual] (0.5,0) to (1.5,0.75);
\draw[usual] (1,0) to (0.5,0.75);
\end{tikzpicture}
\end{gathered} &
\cellcolor{mydarkblue!25}
\begin{gathered}
\begin{tikzpicture}[anchorbase,scale=0.7,tinynodes,yscale=-1]
\draw[usual] (0,0) to[out=90,in=180] (0.75,0.3) to[out=0,in=90] (1.5,0);
\draw[usual] (0.5,0.75) to[out=270,in=180] (1,0.5) to[out=0,in=270] (1.5,0.75);
\draw[usual] (0.5,0) to (0,0.75);
\draw[usual] (1,0) to (1,0.75);
\end{tikzpicture}
\\[3pt]
\begin{tikzpicture}[anchorbase,scale=0.7,tinynodes,yscale=-1]
\draw[usual] (0,0) to[out=90,in=180] (0.75,0.3) to[out=0,in=90] (1.5,0);
\draw[usual] (0.5,0.75) to[out=270,in=180] (1,0.5) to[out=0,in=270] (1.5,0.75);
\draw[usual] (0.5,0) to (1,0.75);
\draw[usual] (1,0) to (0,0.75);
\end{tikzpicture}
\end{gathered} &
\cellcolor{mydarkblue!25}
\begin{gathered}
\begin{tikzpicture}[anchorbase,scale=0.7,tinynodes,yscale=-1]
\draw[usual] (0,0) to[out=90,in=180] (0.75,0.25) to[out=0,in=90] (1.5,0);
\draw[usual] (0,0.75) to[out=270,in=180] (0.75,0.5) to[out=0,in=270] (1.5,0.75);
\draw[usual] (0.5,0) to (0.5,0.75);
\draw[usual] (1,0) to (1,0.75);
\end{tikzpicture}
\\[3pt]
\begin{tikzpicture}[anchorbase,scale=0.7,tinynodes,yscale=-1]
\draw[usual] (0,0) to[out=90,in=180] (0.75,0.25) to[out=0,in=90] (1.5,0);
\draw[usual] (0,0.75) to[out=270,in=180] (0.75,0.5) to[out=0,in=270] (1.5,0.75);
\draw[usual] (0.5,0) to (1,0.75);
\draw[usual] (1,0) to (0.5,0.75);
\end{tikzpicture}
\end{gathered}
\end{tabular}
\end{gathered}};
(-55,0)*{\jcell_{2}};
(58,0)*{\hcell(e)\cong\sym[2]};
\endxy
\quad.
\end{gather*}
This illustrated the cell $\jcell_{2}$ in $\brmon[4]$.

Formally and with contrast to \autoref{P:TLCells}, we have now nontrivial $H$-cells:

\begin{Proposition}\label{P:BrauerCells}
We have the following.
\begin{enumerate}

\item The left and right cells of $\brmon$ 
are given by perfect matchings where one fixes 
the bottom respectively top half of the diagram.
The $\leq_{l}$- and the $\leq_{r}$-order increases as the number of through strands 
decreases. Within $\jcell_{k}$ we have
\begin{gather*}
|\lcell|=|\rcell|=k!\binom{n}{k}(n-k-1)!!.
\end{gather*}
Here $(n-k-1)!!$ denotes the double factorial.

\item The $J$-cells $\jcell_{k}$ of $\brmon$
are given by perfect matchings with a fixed number of through strands $k$.
The $\leq_{lr}$-order is a total order and increases as the number of through strands 
decreases. For any $\lcell\subset\jcell_{k}$ we have
\begin{gather*}
|\jcell_{k}|=\tfrac{1}{k!}|\lcell|^{2}.
\end{gather*}

\item Each $J$-cell of $\brmon$ is idempotent, and $\hcell(e)\cong\sym[k]$ for all idempotent $H$-cells in $\jcell_{k}$. Within $\jcell_{k}$ have
\begin{gather*}
|\hcell|=k!.
\end{gather*}

\end{enumerate}	

\end{Proposition}

\begin{proof}
All of these are known statements. However, the 
cells of $\brmon$ do not correspond to the cells 
coming from the cellular structure of $\bralg{\delta}$, but rather 
from the sandwich cellular structure, {\cf} 
\cite{FiGr-canonical-cases-brauer} or \cite[Section 2D]{TuVa-handlebody}.
\end{proof}

Let $\simple[{\sym[n]}]/{\cong}$ denote the set of simple $\sym[n]$-representations.
For $\cchar=0$ it is well-known that $\simple[{\sym[n]}]/{\cong}$ can be identified 
with partitions 
of $n$. For $\cchar>0$ there is a slightly more involved statement 
of the same kind, see {\eg} \cite[Section 3.4]{Ma-hecke-schur} for an even more general statement.

\begin{Proposition}\label{P:BrauerSimples}
The set of apexes for simple $\brmon$-representations can be indexed $1{:}1$ by 
the poset $\Lambda=(\{n,n-2,\dots\},>)$, and
\begin{gather*}
\{\text{simple $\brmon$-representations of apex $k$}\}/\cong\;
\xleftrightarrow{1{:}1}
\simple[{\sym[k]}]/{\cong}\;
.
\end{gather*}
\end{Proposition}

\begin{proof}
As before by using \autoref{P:CellsSimples} and the 
cell structure in \autoref{P:BrauerCells}.
\end{proof}

By \autoref{P:BrauerSimples}, 
we use the same number scheme and poset as for the Temperley--Lieb monoid
but also keeping track of $\simple[K]\in\simple[{\sym[k]}]/{\cong}$.

\begin{Lemma}\label{L:BrauerCellModules}
Within one $J$-cell all left cell representations $\lmod$ and all right 
cell representations $\rmod$ are isomorphic. We write $\lmod[k]$ respectively $\rmod[k]$ for those in $\jcell_{k}$.

We have $\lmod[k]\cong\rmod[k]$ as $\K$-vector spaces
and $\dimk(\lmod[k])=\dimk(\rmod[k])=\binom{n}{k}(n-k-1)!!$.
\end{Lemma}

\begin{proof}
Using \autoref{P:BrauerCells},
the proof is similar to the Temperley--Lieb case.
\end{proof}

\begin{Proposition}\label{P:BrauerSimplesDim}
Let $K$ be a simple $\sym[k]$-representation, and let $\simple[K]$ denote
its associated simple $\brmon$-representation of apex $\jcell_{k}$. The semisimple dimensions are $\ssdimk(\simple[K])\geq\binom{n}{k}(n-k-1)!!$.
\end{Proposition}

\begin{proof}
As for the Temperley--Lieb case with the extra observation 
that a smallest semisimple dimension is associated to the 
trivial $\sym[k]$-representation.
\end{proof}

\begin{Example}\label{E:BrauerSS}
The lower bound for the semisimple dimensions of $\brmon[24]$
can be illustrated by
\begin{gather*}
\begin{tikzpicture}[anchorbase]
\node at (-0.1,0) {\includegraphics[height=4.4cm]{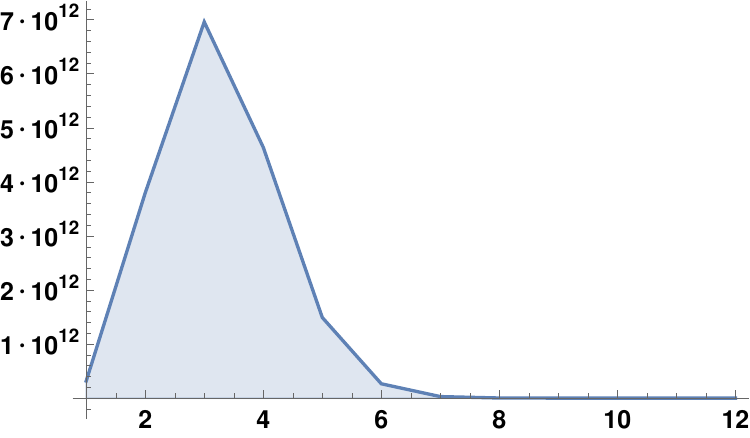}};
\node at (2.25,2) {$\brmon[24]/\K$};
\node at (2.25,1.5) {cells increase $\leftarrow$};
\draw[->] (-2.4,2.2) node[above,xshift=0.3cm]{y-axis: ssdim} to (-2.4,2.1) to (-2.8,2.1);
\draw[->] (2.75,-1.5) node[above,yshift=0.75cm]{x-axis:}node[above,yshift=0.3cm]{\# through}node[above,yshift=-0.15cm]{strands/2} to (2.75,-1.8);
\draw[very thick,densely dotted] (-0.7,-2) to (-0.7,1.5) node[above]{$k\approx 2\sqrt{24}$};
\end{tikzpicture}
,
\\
\scalebox{0.8}{$\text{log${}_{10}$(ssdim)}\colon
\left(
11.5,12.5792,12.8424,12.6663,12.1734,11.4233,10.4489,9.26797,7.88775,6.30512,4.50349,2.44091,0
\right)$}
.
\end{gather*}
For readability, we took the base $10$ log of the actual numbers.
This picture again motivates truncation, and we will do this 
in \autoref{SS:BrauerSubquotient} below.
\end{Example}

Let us discuss the dimensions of simple 
$\brmon$-representations. To the best of our knowledge, the dimensions of simple $\brmon$-representations 
are not known. The best we get is:

\begin{Proposition}\label{P:BrauerBound}
Let $\simple[k]^{\tlmon}$ denote the $k$th simple 
$\tlmon$-representation, {\cf} \autoref{SS:TLCell}. Let 
$K$ be a simple $\sym[k]$-representation and let $\simple[K]$ denote
its associated simple $\brmon$-representation of apex $\jcell_{k}$.
We have $\dimk(\simple[K])\geq\dimk(\simple[k]^{\tlmon})$.
\end{Proposition}

\begin{proof}
The Temperley--Lieb monoid $\tlmon$ embeds into 
$\brmon$ by the evident map that diagrammatically is as the one 
in \autoref{Eq:TLEmbedding}. See also \autoref{L:BrauerTL}.	
The proof is thus essentially the same as for $\momon$, see \autoref{P:TLMotzkinBound}. 
The difference is that we can not use \autoref{T:CellsGramSub} directly, 
but we instead need to argue slightly differently: First, we use the Brauer 
algebra $\bralg{1}$ for circle parameter $1$. The monoid
algebra of $\brmon$ is $\bralg{1}$, hence, finding dimension 
bounds for $\brmon$ or $\bralg{1}$ is the same problem. Working with $\bralg{1}$ 
has the advantage that we can use the cellular structure to split the 
$J$-cells $\jcell_{k}$ further until $H$-cells are of size one, see \cite[Section 4]{GrLe-cellular}. 
This can be achieved by using {\eg} the Kazhdan--Lusztig 
bases of the $\sym[k]$. The sign representation of $\sym[k]$ in its cell structure 
correspond to the bottom cell where there are only through strands. This cell 
for $\bralg{1}$ has then $\tlalg{1}$ inside and the pairing argument 
applies. All other simple $\bralg{1}$-representations associated 
to $\jcell_{k}$ have bigger dimensions, so the proof completes.
\end{proof}


\subsection{Truncating the Brauer monoid}\label{SS:BrauerSubquotient}

We continue with truncation, which is 
almost identical as for the Temperley--Lieb monoid 
in \autoref{SS:TLSubmonoid}.

\begin{Lemma}\label{L:BrauerAdmissible}
The monoid $\brmon$ is regular.
\end{Lemma}

\begin{proof}
The same arguments as in \autoref{L:BrauerAdmissible} work.
In particular, we can use \autoref{L:CellsAdmissible} 
and the well-known fact, that is also easy to prove by hand, 
that $\brmon$ has idempotent $J$-cells, see for example
\cite[Section 3]{KuMaMa-brauer-semigroup}.
\end{proof}

\begin{Definition}\label{D:BrauerTruncated}
Define the \emph{$k$th truncated Brauer monoid} by
\begin{gather*}
\brtru{k}=(\brmon[n])_{\geq\jcell_{k}}.
\end{gather*}
\end{Definition}

Again, let us stress that diagrams in $\brtru{k}$ have at most $k$ 
through strands.
We are almost ready to state our main results, but before 
we need to discuss extensions.


\subsection{Trivial extensions in Brauer monoids}\label{SS:BrauerTrivExtTL} 

The following is the same as for $\tlmon$.

\begin{Lemma}\label{L:BrauerTLWellRounded}
The monoid $\brmon$ is well-connected if $n\geq 5$, and 
the monoid $\brtru{k}$ is well-connected if $n\geq 5$ and $k\leq 3$.
\end{Lemma}

\begin{proof}
This follows from \autoref{L:BrauerFactor} and the respective 
statement about the Temperley--Lieb monoid in \autoref{L:TLWellRounded}.
To see this note that the left-connected condition $ba\approx_{l}a$ implies 
that within on $\approx_{l}$ equivalence class we can focus on the part 
where $\sigma_{k}=1$ since for $a=\gamma\circ\sigma_{k}\circ\beta$ 
we can chose $b=\gamma\circ\sigma_{k}^{-1}\circ\beta^{\ast}$ 
and get $\gamma\circ\id_{k}\circ\beta\approx_{l}\gamma\circ\sigma_{k}\circ\beta$. The same works for right-connected and null-connected.
\end{proof}

We now restrict to a field $\K$ with $\cchar\neq 2$.

\begin{Lemma}\label{L:BrauerSym}
Let $\cchar\neq 2$.
We have $\HH^{1}(\sym,\K)\cong 0$ for all $n\in\N$. 
\end{Lemma}

\begin{proof}
The cases $n=0,1$ are clear, so let $n\geq 2$.
Recall from \autoref{R:RepGapmTrivial} 
that $\HH^{1}(\sym,\K)\cong 0$ is trivial 
if and only if the only homomorphism $\sym\to\K$ is trivial.
To see that this is the case, note that 
any such homomorphism must send the transposition $(i,i+1)$ 
of $\sym=\Aut(\{1,\dots,n\})$ to $k\in\K$ with $2k=0$, which implies $k=0$.
The claim follows since $\sym$ is generated by transpositions.
\end{proof}

Let $\xmon$ be either $\brmon$ or $\brtru{k}$ for $k\geq 3$. 

\begin{Lemma}\label{L:BrauerHTrivial}
Let $\cchar\neq 2$.
We have $\HH^{1}(\xmon,\K)\cong 0$ for all $n\in\N$. 
\end{Lemma}

\begin{proof}
We will use \autoref{L:BrauerSym}.	

\textit{Case $\xmon=\brmon$.}
Similar to 
the proof of \autoref{L:TLHTrivial} with the difference that elements 
in $\jb$ are not generated by idempotents, but rather 
by idempotents and symmetric group generators. 
Idempotents are send to zero, as for the Temperley--Lieb monoid, 
and the symmetric group generators are also send to zero.
These taken together show the claim.
\medskip

\textit{Case $\xmon=\brtru{k}$.} The argument is also similar 
to the proof of \autoref{L:TLHTrivial}. In this case diagrams of width 
$k$ are of the form $a\sigma_{k}b^{\ast}$ where $\sigma_{k}\in\sym$.
Keeping this in mind as well as $\HH^{1}(\sym[k],\K)\cong 0$, the argument 
given in the proof of \autoref{L:TLHTrivial} works {\muta}.
\end{proof}

\begin{Remark}\label{R:BrauerHTrivial}
\autoref{L:BrauerHTrivial} actually works in arbitrary characteristic. To elaborate, in \cite{EaGr-diagram-graphs} the authors 
show that every proper ideal of $\brmon$ is generated as a semigroup
by idempotents. Now, in general, if $\monoid$ is a
monoid containing an ideal $I$ generated by idempotents as a semigroup,
then, for any homomorphism $f\colon\monoid\to\K$, we have
$0=f(I)=f(\monoid I)=f(\monoid)+f(I)=f(\monoid)$. 
This argument really just needs
$\Hom(I,\K)\cong0$, which is a consequence of $I$ being generated by
idempotents as a semigroup, and we get the desired $\HH^{1}(\brmon,\K)\cong 0$ as a special case. The same arguments work for the truncated version.
\end{Remark}

\begin{Proposition}\label{P:BrauerSplits}
Let $\cchar\neq 2$.
Let $\module$ be 
an $\xmon$-representation.
Then any short exact sequence 
\begin{gather*}
0\lra\onebt\lra\module\lra\onebt\lra 0
\end{gather*}
splits.
\end{Proposition} 

\begin{proof}
The proposition follows 
as for the Temperley--Lieb monoid by the above lemmas.
The only difference to \autoref{P:TLSplits} is that the 
group of units is $\group\cong\sym[n]$, but that is taken care of 
in \autoref{L:BrauerSym}.
\end{proof}


\subsection{Representation gap and faithfulness of the Brauer monoid}\label{SS:BrauerMain}

The analog of \autoref{SS:TLMain} is the weaker statement:

\begin{Theorem}\label{T:BrauerGoodExample}
Let $\cchar\neq 2$.
We have the following lower bounds:
\begin{gather*}
\gap{\brtru{k}}\geq\gap[\K]{\tltru[n]{k}}
,\quad
\ssgap{\brtru{k}}\geq\begin{cases*}
\ssgap{\tltru{k}}&\text{always},
\\
(n-1)!!\in\Theta(n^{n/2}e^{n/2})&\text{if $n\gg 0,0\leq k\leq 2\sqrt{n}$},
\end{cases*}
\\
\faith{\brtru{k}}\geq\faith[\K]{\tltru[n]{k}}
.
\end{gather*}
\end{Theorem}

Note that the lower bound $(n-1)!!$ is bigger than 
the one using $\tlmon$ from \autoref{T:TLIsAGoodExample}.

\begin{proof}
Since $\tlmon$ embeds into 
$\brmon$ (see the proof of \autoref{P:BrauerBound}
or \autoref{L:BrauerTL}), \autoref{P:BrauerSplits}
and using the arguments from the proof of 
\autoref{P:BrauerBound}, most of this theorem follows 
from the ones for $\tlmon$ or $\tltru{k}$. 
The exception is the lower bound given by $(n-1)!!$. 
To see that this lower bound holds under the given assumptions, we observe that $\binom{n}{k}(n-k-1)!!$ 
has its minimum at either $k=0$ or $k=\floor{2\sqrt{n}}$. Evaluating at these 
values for $n\gg 0$ (as there are some fluctuations for small $n$) 
shows that the lower bound is achieved at $k=0$.
\end{proof} 


\subsection{Other symmetric monoids}\label{SS:BrauerOther}

We now discuss the remaining symmetric monoids
from \autoref{Eq:IntroDiaMonoids} 
in ascending order (of complexity). We will be brief since almost everything 
follows {\muta} as before.
The basics can be found {\eg} in \cite{HaJa-representations-diagram-algebras}.

\begin{Remark}\label{R:BrauerOther}
Symmetric monoids have the symmetric groups as $H$-cells, 
as we will explain below. As for the planar monoids, 
the cell structure of the symmetric diagram monoids below is 
easy to get. See for example \cite[Chapter 9]{St-rep-monoid}
for many references in the rook monoid case, \cite{HaJa-representations-diagram-algebras} for the rook-Brauer monoid,
and \cite{HaRa-partition-algebras} for 
the partition monoid.
\end{Remark}

The \emph{symmetric group} was discussed in \autoref{E:RepGapSym}, so let 
us start with the \emph{rook monoid} $\romon$. The rook monoid 
is the nonplanar version of $\promon$ and has 
$\sum_{k=0}^{n}k!\binom{n}{k}^{2}$ elements. 
Its $J$-cells are again given by through strands. A typical cell is 
$\jcell_{2}$ for $\romon[3]$:
\begin{gather*}
\xy
(0,0)*{\begin{gathered}
\begin{tabular}{C|C|C}
\arrayrulecolor{tomato}
\begin{gathered}
\cellcolor{mydarkblue!25}
\begin{tikzpicture}[anchorbase]
\draw[usual,dot] (0,0) to (0,0.15);
\draw[usual] (0.5,0) to (0.5,0.5);
\draw[usual] (1,0) to (1,0.5);
\draw[usual,dot] (0,0.5) to (0,0.35);
\end{tikzpicture}
\\[3pt]
\begin{tikzpicture}[anchorbase]
\draw[usual,dot] (0,0) to (0,0.15);
\draw[usual] (0.5,0) to (1,0.5);
\draw[usual] (1,0) to (0.5,0.5);
\draw[usual,dot] (0,0.5) to (0,0.35);
\end{tikzpicture}
\end{gathered} & 
\begin{gathered}
\begin{tikzpicture}[anchorbase]
\draw[usual,dot] (0.5,0) to (0.5,0.15);
\draw[usual] (0,0) to (0.5,0.5);
\draw[usual] (1,0) to (1,0.5);
\draw[usual,dot] (0,0.5) to (0,0.35);
\end{tikzpicture}
\\[3pt]
\begin{tikzpicture}[anchorbase]
\draw[usual,dot] (0.5,0) to (0.5,0.15);
\draw[usual] (0,0) to (1,0.5);
\draw[usual] (1,0) to (0.5,0.5);
\draw[usual,dot] (0,0.5) to (0,0.35);
\end{tikzpicture}
\end{gathered} &
\begin{gathered}
\begin{tikzpicture}[anchorbase]
\draw[usual,dot] (1,0) to (1,0.15);
\draw[usual] (0.5,0) to (0.5,0.5);
\draw[usual] (0,0) to (1,0.5);
\draw[usual,dot] (0,0.5) to (0,0.35);
\end{tikzpicture}
\\[3pt]
\begin{tikzpicture}[anchorbase]
\draw[usual,dot] (1,0) to (1,0.15);
\draw[usual] (0.5,0) to (1,0.5);
\draw[usual] (0,0) to (0.5,0.5);
\draw[usual,dot] (0,0.5) to (0,0.35);
\end{tikzpicture}
\end{gathered}
\\
\hline
\begin{gathered}
\begin{tikzpicture}[anchorbase]
\draw[usual,dot] (0,0) to (0,0.15);
\draw[usual] (0.5,0) to (0,0.5);
\draw[usual] (1,0) to (1,0.5);
\draw[usual,dot] (0.5,0.5) to (0.5,0.35);
\end{tikzpicture}
\\[3pt]
\begin{tikzpicture}[anchorbase]
\draw[usual,dot] (0,0) to (0,0.15);
\draw[usual] (0.5,0) to (1,0.5);
\draw[usual] (1,0) to (0,0.5);
\draw[usual,dot] (0.5,0.5) to (0.5,0.35);
\end{tikzpicture}
\end{gathered} &
\begin{gathered}
\cellcolor{mydarkblue!25}
\begin{tikzpicture}[anchorbase]
\draw[usual,dot] (0.5,0) to (0.5,0.15);
\draw[usual] (0,0) to (0,0.5);
\draw[usual] (1,0) to (1,0.5);
\draw[usual,dot] (0.5,0.5) to (0.5,0.35);
\end{tikzpicture}
\\[3pt]
\begin{tikzpicture}[anchorbase]
\draw[usual,dot] (0.5,0) to (0.5,0.15);
\draw[usual] (0,0) to (1,0.5);
\draw[usual] (1,0) to (0,0.5);
\draw[usual,dot] (0.5,0.5) to (0.5,0.35);
\end{tikzpicture}
\end{gathered} &
\begin{gathered}
\begin{tikzpicture}[anchorbase]
\draw[usual,dot] (1,0) to (1,0.15);
\draw[usual] (0,0) to (0,0.5);
\draw[usual] (0.5,0) to (1,0.5);
\draw[usual,dot] (0.5,0.5) to (0.5,0.35);
\end{tikzpicture}
\\[3pt]
\begin{tikzpicture}[anchorbase]
\draw[usual,dot] (1,0) to (1,0.15);
\draw[usual] (0,0) to (1,0.5);
\draw[usual] (0.5,0) to (0,0.5);
\draw[usual,dot] (0.5,0.5) to (0.5,0.35);
\end{tikzpicture}
\end{gathered}
\\
\hline
\begin{gathered}
\begin{tikzpicture}[anchorbase]
\draw[usual,dot] (0,0) to (0,0.15);
\draw[usual] (1,0) to (0,0.5);
\draw[usual] (0.5,0) to (0.5,0.5);
\draw[usual,dot] (1,0.5) to (1,0.35);
\end{tikzpicture}
\\[3pt]
\begin{tikzpicture}[anchorbase]
\draw[usual,dot] (0,0) to (0,0.15);
\draw[usual] (1,0) to (0.5,0.5);
\draw[usual] (0.5,0) to (0,0.5);
\draw[usual,dot] (1,0.5) to (1,0.35);
\end{tikzpicture}
\end{gathered} & 
\begin{gathered}
\begin{tikzpicture}[anchorbase]
\draw[usual,dot] (0.5,0) to (0.5,0.15);
\draw[usual] (0,0) to (0,0.5);
\draw[usual] (1,0) to (0.5,0.5);
\draw[usual,dot] (1,0.5) to (1,0.35);
\end{tikzpicture}
\\[3pt]
\begin{tikzpicture}[anchorbase]
\draw[usual,dot] (0.5,0) to (0.5,0.15);
\draw[usual] (0,0) to (0.5,0.5);
\draw[usual] (1,0) to (0,0.5);
\draw[usual,dot] (1,0.5) to (1,0.35);
\end{tikzpicture}
\end{gathered} &
\begin{gathered}
\cellcolor{mydarkblue!25}
\begin{tikzpicture}[anchorbase]
\draw[usual,dot] (1,0) to (1,0.15);
\draw[usual] (0,0) to (0,0.5);
\draw[usual] (0.5,0) to (0.5,0.5);
\draw[usual,dot] (1,0.5) to (1,0.35);
\end{tikzpicture}
\\[3pt]
\begin{tikzpicture}[anchorbase]
\draw[usual,dot] (1,0) to (1,0.15);
\draw[usual] (0,0) to (0.5,0.5);
\draw[usual] (0.5,0) to (0,0.5);
\draw[usual,dot] (1,0.5) to (1,0.35);
\end{tikzpicture}
\end{gathered}
\end{tabular}
\end{gathered}};
(-30,0)*{\jcell_{2}};
(35,0)*{\hcell(e)\cong\sym[2]};
\endxy
\quad.
\end{gather*}

The \emph{rook-Brauer monoid} $\robrmon$ is a symmetric version of the Motzkin 
monoid. The rook-Brauer monoid has 
$\sum_{k=0}^{n}k!\big(\sum_{t=0}^{n}\binom{n}{k}\binom{n-k}{2t}(2t-1)!!\big)^{2}$ elements. The $J$-cells are, as usual, indexed by 
through strands. They get huge very fast, so let us just illustrate 
a typing idempotent (and symmetric) $H$-cell:
\begin{gather*}
\xy
(0,0)*{\begin{gathered}
\begin{tabular}{C}
\arrayrulecolor{tomato}
\begin{gathered}
\cellcolor{mydarkblue!25}
\begin{tikzpicture}[anchorbase]
\draw[usual] (0,0) to[out=90,in=180] (0.5,0.2) to[out=0,in=90] (1,0);
\draw[usual] (0,0.5) to[out=270,in=180] (0.5,0.3) to[out=0,in=270] (1,0.5);
\draw[usual] (2,0) to[out=90,in=180] (2.5,0.22) to[out=0,in=90] (3,0);
\draw[usual] (2,0.5) to[out=270,in=180] (2.5,0.28) to[out=0,in=270] (3,0.5);
\draw[usual] (-0.5,0) to (-0.5,0.5);
\draw[usual] (0.5,0) to (0.5,0.5);
\draw[usual] (1.5,0) to (1.5,0.5);
\draw[usual,dot] (2.5,0) to (2.5,0.13);
\draw[usual,dot] (2.5,0.5) to (2.5,0.37);
\end{tikzpicture}
,\;
\begin{tikzpicture}[anchorbase]
\draw[usual] (0,0) to[out=90,in=180] (0.5,0.2) to[out=0,in=90] (1,0);
\draw[usual] (0,0.5) to[out=270,in=180] (0.5,0.3) to[out=0,in=270] (1,0.5);
\draw[usual] (2,0) to[out=90,in=180] (2.5,0.22) to[out=0,in=90] (3,0);
\draw[usual] (2,0.5) to[out=270,in=180] (2.5,0.28) to[out=0,in=270] (3,0.5);
\draw[usual] (-0.5,0) to (0.5,0.5);
\draw[usual] (0.5,0) to (-0.5,0.5);
\draw[usual] (1.5,0) to (1.5,0.5);
\draw[usual,dot] (2.5,0) to (2.5,0.13);
\draw[usual,dot] (2.5,0.5) to (2.5,0.37);
\end{tikzpicture}
\\[3pt]
\begin{tikzpicture}[anchorbase]
\draw[usual] (0,0) to[out=90,in=180] (0.5,0.2) to[out=0,in=90] (1,0);
\draw[usual] (0,0.5) to[out=270,in=180] (0.5,0.3) to[out=0,in=270] (1,0.5);
\draw[usual] (2,0) to[out=90,in=180] (2.5,0.22) to[out=0,in=90] (3,0);
\draw[usual] (2,0.5) to[out=270,in=180] (2.5,0.28) to[out=0,in=270] (3,0.5);
\draw[usual] (-0.5,0) to (-0.5,0.5);
\draw[usual] (0.5,0) to (1.5,0.5);
\draw[usual] (1.5,0) to (0.5,0.5);
\draw[usual,dot] (2.5,0) to (2.5,0.13);
\draw[usual,dot] (2.5,0.5) to (2.5,0.37);
\end{tikzpicture}
,\;
\begin{tikzpicture}[anchorbase]
\draw[usual] (0,0) to[out=90,in=180] (0.5,0.2) to[out=0,in=90] (1,0);
\draw[usual] (0,0.5) to[out=270,in=180] (0.5,0.3) to[out=0,in=270] (1,0.5);
\draw[usual] (2,0) to[out=90,in=180] (2.5,0.22) to[out=0,in=90] (3,0);
\draw[usual] (2,0.5) to[out=270,in=180] (2.5,0.28) to[out=0,in=270] (3,0.5);
\draw[usual] (-0.5,0) to (0.5,0.5);
\draw[usual] (0.5,0) to (1.5,0.5);
\draw[usual] (1.5,0) to (1.5,0.25) to (-0.5,0.5);
\draw[usual,dot] (2.5,0) to (2.5,0.13);
\draw[usual,dot] (2.5,0.5) to (2.5,0.37);
\end{tikzpicture}
\\[3pt]
\begin{tikzpicture}[anchorbase]
\draw[usual] (0,0) to[out=90,in=180] (0.5,0.2) to[out=0,in=90] (1,0);
\draw[usual] (0,0.5) to[out=270,in=180] (0.5,0.3) to[out=0,in=270] (1,0.5);
\draw[usual] (2,0) to[out=90,in=180] (2.5,0.22) to[out=0,in=90] (3,0);
\draw[usual] (2,0.5) to[out=270,in=180] (2.5,0.28) to[out=0,in=270] (3,0.5);
\draw[usual] (-0.5,0) to (-0.5,0.25) to (1.5,0.5);
\draw[usual] (0.5,0) to (-0.5,0.5);
\draw[usual] (1.5,0) to (0.5,0.5);
\draw[usual,dot] (2.5,0) to (2.5,0.13);
\draw[usual,dot] (2.5,0.5) to (2.5,0.37);
\end{tikzpicture}
,\;
\begin{tikzpicture}[anchorbase]
\draw[usual] (0,0) to[out=90,in=180] (0.5,0.2) to[out=0,in=90] (1,0);
\draw[usual] (0,0.5) to[out=270,in=180] (0.5,0.3) to[out=0,in=270] (1,0.5);
\draw[usual] (2,0) to[out=90,in=180] (2.5,0.22) to[out=0,in=90] (3,0);
\draw[usual] (2,0.5) to[out=270,in=180] (2.5,0.28) to[out=0,in=270] (3,0.5);
\draw[usual] (-0.5,0) to (-0.5,0.2) to (1.5,0.5);
\draw[usual] (0.5,0) to (0.5,0.5);
\draw[usual] (1.5,0) to (-0.5,0.3) to (-0.5,0.5);
\draw[usual,dot] (2.5,0) to (2.5,0.13);
\draw[usual,dot] (2.5,0.5) to (2.5,0.37);
\end{tikzpicture}
\end{gathered}
\end{tabular}
\end{gathered}};
(-55,0)*{\text{Part of }\jcell_{3}};
(55,0)*{\hcell(e)\cong\sym[3]};
\endxy
\quad.
\end{gather*}

The \emph{partition monoid} $\pamon$ contains all the other 
planar and symmetric monoids as submonoid. It has $Be(2n)$ elements, where 
$Be$ denotes the Bell number.
The $J$-cells are still given by through strands. As for $\robrmon$, 
the sizes of the cells are very large, so we only illustrate an idempotent 
(and symmetric) $H$-cell:
\begin{gather*}
\xy
(0,0)*{\begin{gathered}
\begin{tabular}{C}
\arrayrulecolor{tomato}
\begin{gathered}
\cellcolor{mydarkblue!25}
\begin{tikzpicture}[anchorbase]
\draw[usual] (0,0) to[out=90,in=180] (0.5,0.2) to[out=0,in=90] (1,0);
\draw[usual] (0,0.5) to[out=270,in=180] (0.5,0.3) to[out=0,in=270] (1,0.5);
\draw[usual] (1,0) to[out=90,in=180] (1.5,0.2) to[out=0,in=90] (2,0);
\draw[usual] (1,0.5) to[out=270,in=180] (1.5,0.3) to[out=0,in=270] (2,0.5);
\draw[usual] (2,0) to[out=90,in=180] (2.25,0.2) to[out=0,in=90] (2.5,0);
\draw[usual] (2,0.5) to[out=270,in=180] (2.25,0.3) to[out=0,in=270] (2.5,0.5);
\draw[usual] (0,0) to (0,0.5);
\draw[usual] (0.5,0) to (0.5,0.5);
\draw[usual] (1.5,0) to (1.5,0.5);
\end{tikzpicture}
\\[3pt]
\begin{tikzpicture}[anchorbase]
\draw[usual] (0,0) to[out=90,in=180] (0.5,0.2) to[out=0,in=90] (1,0);
\draw[usual] (0,0.5) to[out=270,in=180] (0.5,0.3) to[out=0,in=270] (1,0.5);
\draw[usual] (1,0) to[out=90,in=180] (1.5,0.2) to[out=0,in=90] (2,0);
\draw[usual] (1,0.5) to[out=270,in=180] (1.5,0.3) to[out=0,in=270] (2,0.5);
\draw[usual] (2,0) to[out=90,in=180] (2.25,0.2) to[out=0,in=90] (2.5,0);
\draw[usual] (2,0.5) to[out=270,in=180] (2.25,0.3) to[out=0,in=270] (2.5,0.5);
\draw[usual] (0,0) to (0,0.5);
\draw[usual] (0.5,0) to (1.5,0.5);
\draw[usual] (1.5,0) to (0.5,0.5);
\end{tikzpicture}
\end{gathered}
\end{tabular}
\end{gathered}};
(-35,0)*{\text{Part of }\jcell_{2}};
(35,0)*{\hcell(e)\cong\sym[2]};
\endxy
\quad.
\end{gather*}

Below, if not stated otherwise, let $\xmon$ be $\romon$, $\robrmon$ or $\pamon$.

\begin{Proposition}\label{P:BrauerCellsOther}
We have the following.
\begin{enumerate}

\item The left and right cells of $\xmon$ 
are given by the respective type of diagrams where one fixes 
the bottom respectively top half of the diagram.
The $\leq_{l}$- and the $\leq_{r}$-order increases as the number of through strands 
decreases. Within $\jcell_{k}$ we have
\begin{gather*}
\romon\colon|\lcell|=|\rcell|=k!\binom{n}{k},
\\
\robrmon\colon|\lcell|=|\rcell|=k!\sum_{t=0}^{n}\binom{n}{k}\binom{n-k}{2t}(2t-1)!!,
\\
\pamon\colon|\lcell|=|\rcell|=k!\sum_{t=0}^{n}\begin{Bmatrix}n\\t\end{Bmatrix}\binom{t}{k}.
\end{gather*}
Here $\begin{Bsmallmatrix}n\\t\end{Bsmallmatrix}$ denotes the Stirling 
number of the second kind.

\item The $J$-cells $\jcell_{k}$ of $\xmon$
are given by the respective type of diagrams with a fixed number of through strands $k$.
The $\leq_{lr}$-order is a total order and increases as the number of through strands 
decreases. For any $\lcell\subset\jcell_{k}$ we have
\begin{gather*}
\xmon\colon|\jcell_{k}|=\tfrac{1}{k!}|\lcell|^{2}.
\end{gather*}

\item Each $J$-cell of $\xmon$ is idempotent, 
and $\hcell(e)\cong\sym[k]$ for all idempotent $H$-cells in $\jcell_{k}$. Within $\jcell_{k}$ have
\begin{gather*}
\xmon\colon|\hcell|=k!.
\end{gather*}

\end{enumerate}	

\end{Proposition}

\begin{proof}
Easy and omitted, see also \cite[Section 3.3]{HaJa-representations-diagram-algebras}.
\end{proof}

\begin{Proposition}\label{P:BrauerSimplesOther}
The set of apexes for simple $\xmon$-representations can be indexed $1{:}1$ by 
the poset $\Lambda=(\{n,n-1,\dots\},>)$, and
\begin{gather*}
\{\text{simple $\xmon$-representations of apex $k$}\}/\cong\;
\xleftrightarrow{1{:}1}
\simple[{\sym[k]}]/\cong\;
.
\end{gather*}
\end{Proposition}

\begin{proof}
By \autoref{P:BrauerCellsOther}.
\end{proof}

\begin{Proposition}\label{P:BrauerSimplesDimOther}
Let $K$ be a simple $\sym[k]$-representation, and let $\simple[K]$ denote
its associated simple $\xmon$-representation of apex $\jcell_{k}$. 
The semisimple dimensions are $\ssdimk(\simple[K])\geq\binom{n}{k}$, 
$\ssdimk(\simple[K])\geq\sum_{t=0}^{n}\binom{n}{k}\binom{n-k}{2t}(2t-1)!!$ respectively 
$\ssdimk(\simple[K])\geq\sum_{t=0}^{n}\begin{Bsmallmatrix}n\\t\end{Bsmallmatrix}\binom{t}{k}$.
\end{Proposition}

\begin{proof}
This follows {\ver} as \autoref{P:BrauerSimplesDim}.
\end{proof}

\begin{Example}\label{E:BrauerOther}
As before, let us illustrate the lower bound for the 
semisimple dimensions:
\begin{gather*}
\begin{tikzpicture}[anchorbase]
\node at (-0.1,0) {\includegraphics[height=4.4cm]{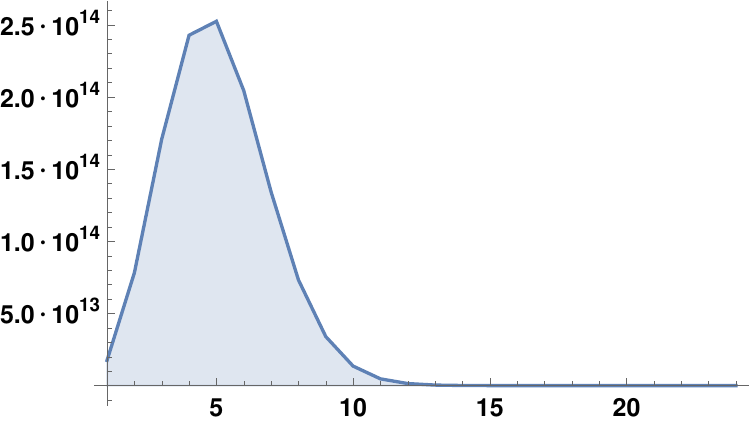}};
\node at (2.25,2) {$\robrmon[24]/\K$};
\node at (2.25,1.5) {cells increase $\leftarrow$};
\draw[->] (-2.4,2.2) node[above,xshift=0.3cm]{y-axis: ssdim} to (-2.4,2.1) to (-2.8,2.1);
\draw[->] (2.75,-1.5) node[above,yshift=0.75cm]{x-axis:}node[above,yshift=0.3cm]{\# through}node[above,yshift=0.0cm]{strands} to (2.75,-1.8);
\draw[very thick,densely dotted] (-0.7,-2) to (-0.7,1.5) node[above]{$k\approx 2\sqrt{24}$};
\end{tikzpicture}
,
\begin{tikzpicture}[anchorbase]
\node at (-0.1,0) {\includegraphics[height=4.4cm]{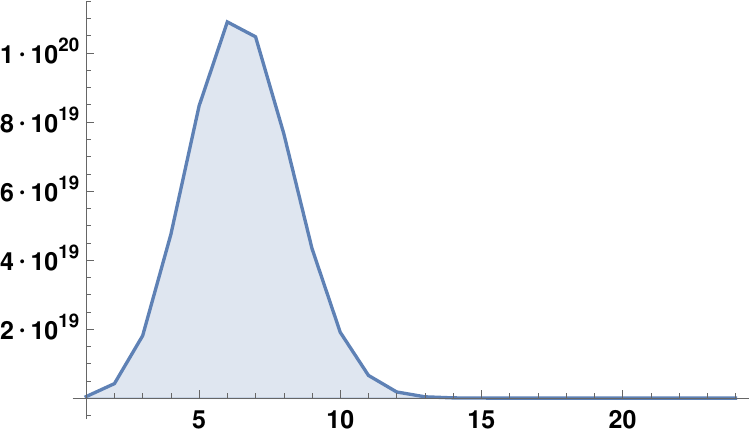}};
\node at (2.25,2) {$\pamon[24]/\K$};
\node at (2.25,1.5) {cells increase $\leftarrow$};
\draw[->] (-2.4,2.2) node[above,xshift=0.3cm]{y-axis: ssdim} to (-2.4,2.1) to (-2.8,2.1);
\draw[->] (2.75,-1.5) node[above,yshift=0.75cm]{x-axis:}node[above,yshift=0.3cm]{\# through}node[above,yshift=0.0cm]{strands} to (2.75,-1.8);
\draw[very thick,densely dotted,xshift=0.3cm] (-0.7,-2) to (-0.7,1.5) node[above]{$k\approx 2\sqrt{24}$};
\end{tikzpicture}
,
\\
\scalebox{0.7}{$\text{RoBr log${}_{10}$(ssdim)}\colon
\left(
\begin{gathered}
13.2428,13.8931,14.2325,14.3852,14.4022,14.3105,14.1279,13.8651,13.5313,13.1314,12.671,12.152,
\\
11.5786,10.9499,10.2701,9.53474,8.74966,7.90469,7.00985,6.0434,5.02637,3.90827,2.74194,1.38021,0
\end{gathered}
\right)$}
,
\\
\scalebox{0.7}{$\text{Pa log${}_{10}$(ssdim)}\colon
\left(
\begin{gathered}
17.6493,18.6225,19.2572,19.6761,19.9277,20.0373,20.0198,19.8843,19.6367,19.2804,18.8176,18.2495,
\\
17.5764,16.7982,15.9143,14.9234,13.8234,12.6115,11.2832,9.83189,8.24761,6.51485,4.60773,2.47712,0
\end{gathered}
\right)$}
,
\end{gather*}
where we again took the base $10$ log.
We have not illustrated the situation for $\romon$ 
as it is the same as for $\promon$ with semisimple 
dimensions given by binomial coefficients, {\cf} \autoref{Eq:CellsDimTL}.
\end{Example}

Let $\xmon[Y]$ denote $\promon$, $\momon$ or $\ppamon$ associated 
to their respective $\xmon$.

\begin{Lemma}\label{L:BrauerEmbedding}
The monoid $\xmon$
contains $\xmon[Y]$ as a submonoid, by the analog 
of \autoref{Eq:TLEmbedding}.
\end{Lemma}

\begin{proof}
Clear, see also \autoref{L:BrauerTL}.
\end{proof}

\begin{Proposition}\label{P:BrauerBoundOther}
Let $\simple[k]^{\xmon[Y]}$ denote the $k$th simple 
$\xmon[Y]$-representation. Let 
$K$ denote a simple $\sym[k]$-representation and let $\simple[K]$ denote
its associated simple $\xmon$-representation of apex $\jcell_{k}$.
We have $\dimk(\simple[K]^{\xmon})\geq\dimk(\simple[k]^{\xmon[Y]})$.
\end{Proposition}

\begin{proof}
This follows again by observing that the planar version embed into their 
nonplanar counterparts, see \autoref{L:BrauerEmbedding}.
\end{proof}

\begin{Proposition}\label{P:BrauerDimsRook}
Let $\cchar\nmid n!$, including $\cchar=0$.
We have $\dimk(\simple[k])=\ssdimk(\simple[k])=\binom{n}{k}$ for $\romon$, and $\romon$ is semisimple.
\end{Proposition}

\begin{proof}
The argument is the same as in \autoref{P:TLDimsPRook} with the additional 
caveat of the symmetric groups $\sym[k]$ for $0\leq k\leq n$ 
appearing as idempotent $H$-cells 
which forces the condition $\cchar\nmid n!$.
\end{proof}

\begin{Lemma}\label{L:BrauerAdmissibleOther}
The monoid $\xmon$ is regular.
\end{Lemma}

\begin{proof}
This follows as before from \autoref{L:CellsAdmissible} 
and a construction of an idempotent for each $J$-cell. 
The latter is easy and omitted, but also well-known, see {\eg} the references in \autoref{R:BrauerOther}.
\end{proof}

\begin{Definition}\label{D:BrauerTruncatedOther}
Define the \emph{$k$-$l$ truncated rook monoid} 
for $k\leq l$ and the
\emph{$k$th truncated rook-Brauer monoid} respectively \emph{$k$ truncated partition monoid} by
\begin{gather*}
\rotru{k}{l}=(\romon[n])_{\geq\jcell_{k}}/(>\jcell_{l})
,\quad
\robrtru{k}=(\robrmon[n])_{\geq\jcell_{k}},\quad
\patru{k}=(\pamon[n])_{\geq\jcell_{k}}.
\end{gather*}
\end{Definition}

Let $\xmon$ be either of the above monoids or their truncations.
For the following theorem, note that the \emph{$k$th 
truncated planar partition monoid} $\ppatru{k}$ 
can be defined in the evident way.

\begin{Theorem}\label{T:BrauerIsAGoodExampleOther}
Let $\cchar\nmid n!$, including $\cchar=0$, and let $\LL$ 
be an arbitrary field.
We have the following lower bounds:
\begin{gather*}
\gap[\K]{\rotru{k}{l}}\geq\gap[\K]{\protru{k}{l}}
\\
\ssgap[\LL]{\rotru{k}{l}}=\ssgap[\LL]{\protru{k}{l}}
,\quad
\faith[\LL]{\rotru{k}{l}}\geq\faith[\LL]{\protru{k}{l}}
,
\\
\ssgap[\LL]{\robrtru{k}}\geq
\begin{cases*}
\ssgap[\LL]{\motru{k}}&\text{always},
\\
\sum_{t=0}^{n}\binom{n}{2t}(2t-1)!!&\text{if $n\gg 0,0\leq k\leq 2\sqrt{n}$},
\end{cases*}
\\
\faith[\LL]{\robrtru{k}}\geq\faith[\LL]{\motru{k}}
\\
\ssgap[\LL]{\patru{k}}\geq
\begin{cases*}
\ssgap[\LL]{\ppatru{k}}&\text{always},
\\
\sum_{t=0}^{n}\begin{Bsmallmatrix}n\\t\end{Bsmallmatrix}&\text{if $n\gg 0,0\leq k\leq 2\sqrt{n}$},
\end{cases*}
\\
\faith[\LL]{\patru{k}}\geq\faith[\LL]{\ppatru{k}}
.
\end{gather*}
\end{Theorem}

Note that the above lower bounds in the cases $n\gg 0,0\leq k\leq 2\sqrt{n}$ 
are bigger than the ones coming from the embeddings. We also expect 
that
\begin{gather*}
\gap{\robrtru{k}}\geq\gap{\motru{k}}
,\quad
\gap{\patru{k}}\geq\gap{\ppatru{k}}
,
\end{gather*}
but we were not able to prove this since there might be extensions.

\begin{proof}
All lower bounds except the first follow directly by
using the embedding in \autoref{L:BrauerEmbedding}. 
The first uses 
additionally \autoref{P:BrauerDimsRook} which also holds for the truncation.

The equality $\ssgap{\rotru{k}{l}}=\ssgap{\protru{k}{l}}$ is clear by \autoref{P:BrauerCellsOther}. For the semisimple gaps of 
$\robrmon$ and $\patru{k}$ one can use the same arguments as 
in \autoref{T:BrauerGoodExample}.
\end{proof}

\begin{Remark}\label{R:BrauerIsAGoodExampleOther}
The exact value for $\faith[\LL]{\rotru{k}{l}}$ can be computed 
using \cite[Theorems 15 and 17]{MaSt-effective-dimension-semigroups}.
The methods from that paper together with \cite{EaMiRuTo-congruence-lattices} may also be used to compute 
the faithfulness of other diagram monoids and their truncations.
\end{Remark}	

\begin{Conclusion}\label{C:BrauerConclusion}
As with the planar monoids, 
all of the symmetric monoids $\romon$, $\brmon$, 
$\robrmon$ and $\pamon$ appear to have  big nontrivial 
representations. However, it is 
not clear why they should be preferable over 
their planar counterparts since they are, roughly speaking, their
planar version inflated by the symmetric group $\sym[k]$. 
In fact most of our arguments above use the planar versions 
to derive bounds.
\end{Conclusion}


\newcommand{\etalchar}[1]{$^{#1}$}

\end{document}